\documentclass[a4paper]{amsart}
\usepackage[hypertex]{hyperref}
\usepackage{amsmath,amssymb}
\usepackage{bbold}
\usepackage{psfrag}
\usepackage{graphicx}
\usepackage{subfigure}
\usepackage{leftidx}
\usepackage[all]{xy}
\usepackage[latin1]{inputenc}

\newtheorem{thm}{Theorem}[section]
\newtheorem{cor}[thm]{Corollary}
\newtheorem{lem}[thm]{Lemma}
\newtheorem{prop}[thm]{Proposition}

\theoremstyle{definition}

\newtheorem{rem}[thm]{Remark}
\newtheorem{exa}[thm]{Example}

\newcommand{\co}{\colon}
\newcommand{\id}{\mathrm{id}}
\newcommand{\kk}{\Bbbk}
\newcommand{\kt}{$\Bbbk$\nobreakdash-\hspace{0pt}}

\newcommand{\opp}{\mathrm{op}}
\newcommand{\cop}{\mathrm{cop}}
\newcommand{\mop}{\mathrm{o}}

\newcommand{\cc}{\mathcal{C}}
\newcommand{\bb}{\mathcal{B}}
\newcommand{\bbb}{\overline{\bb}}
\newcommand{\dd}{\mathcal{D}}
\newcommand{\aaa}{\mathcal{A}}
\newcommand{\ee}{\mathcal{E}}
\newcommand{\ff}{\mathcal{F}}

\newcommand{\uu}{\mathcal{U}}

\newcommand{\zz}{\mathcal{Z}}

\newcommand{\rmod}[2]{\leftidx{}{#1}{_{#2}}}
\newcommand{\lmod}[2]{\leftidx{_{#2}}{#1}{}}

\newcommand{\dinato}{\to}
\newcommand{\iso}{\stackrel{\sim}{\longrightarrow}}

\newcommand{\trait}{\nobreakdash-\hspace{0pt}}
\newcommand{\Rt}{$\mathrm{R}$\nobreakdash-\hspace{0pt}}
\newcommand{\ti}{\mbox{-}\,}
\newcommand{\un}{\mathbb{1}}
\newcommand{\cp}{\rtimes}

\newcommand{\Ob}{\mathrm{Ob}}

\newcommand{\End}{\mathrm{End}}
\newcommand{\Hom}{\mathrm{Hom}}
\newcommand{\Nat}{{\textsc{Hom}}}

\newcommand{\Fun}{{\textsc{Fun}}}

\newcommand{\Dinat}{{\textsc{Dinat}}}
\newcommand{\Vect}{{\mathrm{Vect}}}
\newcommand{\vect}{\mathrm{vect}}

\newcommand{\rank}{\mathrm{rank}}
\newcommand{\lev}{\mathrm{ev}}
\newcommand{\rev}{\widetilde{\mathrm{ev}}}
\newcommand{\lcoev}{\mathrm{coev}}
\newcommand{\rcoev}{\widetilde{\mathrm{coev}}}

\newcommand{\ldual}[1]{\leftidx{^\vee}{\!#1}{}}
\newcommand{\rdual}[1]{{#1}^\vee}

\newcommand{\lldual}[1]{\leftidx{^{\vee\vee}}{\!#1}{}}

\newcommand{\rexcla}[1]{\leftidx{}{#1}{^!}}
\newcommand{\adjunct}[2]{\!\!\raisebox{.6ex}{\xymatrix{\ar@/^.4pc/[r]^{#1}  &  \ar@/^.4pc/[l]^{#2}}}\!\!}
\newcommand{\rsdraw}[3]{\raisebox{-#1\height}{\scalebox{#2}{\includegraphics{#3.eps}}}}
\providecommand{\bysame}{\leavevmode\hbox to3em{\hrulefill}\thinspace}

\begin{document}
\title{The double of a Hopf monad}
\author[A. Brugui\`eres]{Alain Brugui\`eres}
\author[A. Virelizier]{Alexis Virelizier}
\email{bruguier@math.univ-montp2.fr \and virelizi@math.univ-montp2.fr}

\date{\today}

\begin{abstract} The center $\zz(\cc)$ of an autonomous category $\cc$ is monadic over $\cc$ (if certain coends exist in $\cc$).
The notion of Hopf monad naturally arises if one tries to reconstruct the structure of  $\zz(\cc)$
in terms of its monad $Z$: we show that $Z$ is a quasitriangular Hopf monad on $\cc$ and $\zz(\cc)$ is  isomorphic to the braided category $Z\ti\cc$ of $Z$\ti modules.
More generally, let $T$ be a Hopf monad on an autonomous category $\cc$.
We construct a Hopf monad $Z_T$ on~$\cc$, the \emph{centralizer of $T$}, and a canonical  distributive law $\Omega\co TZ_T \to Z_T T$. By Beck's theory, this has two consequences. On one hand, $D_T=Z_T \circ_\Omega T$ is a quasitriangular Hopf monad on $\cc$, called the \emph{double of~$T$}, and $\zz(T\ti\cc)\simeq D_T \ti \cc$ as braided categories.
As an illustration, we
define the double $D(A)$ of a Hopf algebra $A$ in a braided autonomous category in such a way that the center of the category of $A$\ti modules is the braided category of $D(A)$\ti modules (generalizing the Drinfeld double).
On the other hand, the canonical distributive law $\Omega$ also lifts
$Z_T$ to a Hopf monad $\Tilde{Z}_T^\Omega$ on~$T\ti\cc$, and $\Tilde{Z}_T^\Omega(\un,T_0)$ is the coend of~$T\ti\cc$.
For $T=Z$, this gives   an explicit description of the Hopf algebra structure of the coend of $\zz(\cc)$ in terms of the structural morphisms of $\cc$. Such a description is useful in quantum topology, especially when
$\cc$ is a spherical fusion category, as  $\zz(\cc)$ is then modular.
\end{abstract}
\maketitle

\setcounter{tocdepth}{1} \tableofcontents

\section*{Introduction}
The center $\zz(\cc)$ of an autonomous category $\cc$, introduced by Drinfeld, is a braided autonomous category.
This construction establishes a bridge between the non-braided world and the braided world. It is useful, in particular, for comparing quantum invariants of $3$ manifolds. Indeed,  the center $\zz(\cc)$ of spherical fusion category $\cc$ is modular (see \cite{Mueg}) and it is conjectured  that the Turaev-Viro  invariant $\mathrm{TV}_\cc$ (as revisited in \cite{BW}) is equal to the Reshetikhin-Turaev invariant $\mathrm{RT}_{\zz(\cc)}$ (see \cite{Tur2}).

Let $\cc$ be an autonomous category. If the coend:
$$
Z(X)=\int^{Y \in \cc} \ldual{Y} \otimes X \otimes Y
$$
exists for all object $X$ of $\cc$, then Day and Street \cite{DayStreet} showed that $Z$ is a monad on $\cc$ and the center $\zz(\cc)$ is isomorphic to the category $Z\ti\cc$ of $Z$\ti modules in $\cc$ (also called $Z$\ti algebras).
By Tannaka reconstruction, we endow the monad $Z$ with a quasitriangular Hopf monad structure which reflects the braided autonomous structure of $\zz(\cc)$ in the sense that $\zz(\cc)\simeq Z\ti\cc$ as braided categories. The notion of Hopf monad, which generalizes Hopf algebras to the non-braided (and non-linear) setting, was introduced in \cite{BV2} for this very purpose.

The Reshetikhin-Turaev invariant can be expressed in terms of the simple objects of the category (as in
Reshetikhin and Turaev's original construction) or in terms of the coend of the category (following Lyubashenko, see~\cite{Lyu2,BV1}).
In order to compute $\mathrm{RT}_{\zz(\cc)}$, the first approach is not practicable for lack of a workable
description of the simple objects of the center. And so we need to provide an  explicit description of the coend of $\zz(\cc)$ and its algebraic structure. To fulfill this objective, we extend the previous construction of $Z$ to a more general situation.
Let $T$ be a Hopf monad on an autonomous category $\cc$. We denote $T\ti\cc$ the category of $T$\ti modules (also called $T$\ti algebras), which is autonomous. Assume $T$ is \emph{centralizable},  meaning that the coend:
$$
Z_T(X)=\int^{Y \in \cc} \ldual{T}(Y) \otimes X \otimes Y.
$$
exists for every object $X$ of $\cc$. We endow $Z_T$ with a structure of a Hopf monad on~$\cc$ and call $Z_T$ the \emph{centralizer of $T$}.  In particular,
$Z_{1_\cc}=Z$ as Hopf monads. Note that the coend of $\cc$ is $Z(\un)=Z_{1_\cc}(\un)$, and so the coend of $\zz(\cc)$ is $Z_{1_{\zz(\cc)}}(\un)=Z_{1_{Z\ti\cc}}(\un)$.


Using adjunction and exactness properties of Hopf monads, we show that
$1_{T\ti\cc}$ is centralizable and:
\begin{equation*}
U_TZ_{1_{T\ti\cc}}=Z_TU_T.
\end{equation*}
Note that this implies that, as an object of $\cc$, the coend of the category $T\ti\cc$ is $Z_T(\un)$  and, in particular, the coend of $\zz(\cc)=Z\ti\cc$ is $Z_Z(\un)$.
Now $U_TZ_{1_{T\ti\cc}}=Z_TU_T$ means in fact that the Hopf monad $Z_{1_{T\ti\cc}}$ is a lift to $T\ti\cc$ of the Hopf monad~$Z_T$. Extending Beck's theory of distributive laws to Hopf monads, we show that such a lift is encoded by an invertible comonoidal distributive law $\Omega \co T Z_T \to Z_T T$, called the \emph{canonical distributive law of $T$}. The coend of $T\ti\cc$ is therefore $(Z_T(\un),Z_T(T_0) \Omega_\un)$. When $T$ is quasitriangular, this coend has a structure of a Hopf algebra in the braided autonomous category $T\ti\cc$, which we elucidate in terms of $T$. Hence, for $T=Z$, an explicit description of the coend of $\zz(\cc)$.  The case of fusion categories is treated in detail.

The canonical distributive law $\Omega$ also endows the composition of $Z_T$ and $T$ with a Hopf monad structure,
denoted $D_T=Z_T \circ_\Omega T$ and called the \emph{double of $T$}. We prove that $D_T$ is quasitriangular and give a braided isomorphism:
$$D_T \ti \cc \simeq \zz(T\ti\cc).$$

This construction, which holds for any centralizable Hopf monad on an autonomous category, generalizes the Drinfeld
double in an non-braided setting. As an illustration, we apply this to Hopf monads associated with Hopf algebras.
This leads to the double $D(A)$ of a Hopf algebra $A$ in a braided autonomous category $\bb$. More precisely,  the endofunctor $? \otimes A$  is a Hopf monad on $\bb$. Assuming $\bb$ admits a coend $C$, the Hopf monad $?\otimes A$  is centralizable, and its centralizer is of the form $? \otimes Z(A)$, where $Z(A)=\ldual{A}\otimes C$ is a Hopf algebra in $\bb$.
The canonical distributive law of $? \otimes A$ is of the form $\id_{1_\bb} \otimes \Omega$, where $\Omega \co Z(A) \otimes A \to A \otimes Z(A)$ is
a distributive law of Hopf algebras. Then $D(A)=A \otimes_\Omega Z(A)$ is a quasitriangular Hopf algebra in $\bb$,
such that:
\begin{equation*}
\zz(\rmod{\bb}{A}) \simeq \rmod{\bb}{{D(A)}}\simeq\lmod{\bb}{{D(A)}} \simeq \overline{\zz(\lmod{\bb}{A})}.
\end{equation*}
as braided categories, where $\lmod{\bb}{A}$ and $\rmod{\bb}{A}$ denote the categories of left and right modules over $A$. Note that a Hopf algebra $B$ in $\bb$ is quasitriangular when it is endowed with a \Rt matrix. In this context,
we define \Rt matrices to be morphisms $r\co C \otimes C \to B \otimes B$ which encode braidings on $\rmod{\bb}{B}$ (or equivalently $\lmod{\bb}{B}$).
 When $\bb$ is the category of finite-dimensional vector spaces over a field $\kk$, we recover the usual definition of \Rt matrices and the Drinfeld double of a Hopf algebra $H$. Indeed, in that case: $C=\kk$, $Z(H)={H^*}^{\cop}$, and $D(H)=H \otimes_\Omega {H^*}^{\cop}$.

The canonical distributive law of a Hopf monad
is in fact naturally defined in a more general setting.
Let $T$ be a Hopf monad on an autonomous category $\cc$ and $Q$ be a Hopf monad on $T\ti\cc$. Their cross product $Q\cp T=U_TQF_T$  is a Hopf monad on~$\cc$. If $Q\cp T$ is centralizable, then so is~$Q$ and the Hopf monad $Z_Q$ is a lift to $T\ti\cc$ of the Hopf monad $Z_{Q\cp T}$:
\begin{equation*}
U_TZ_Q=Z_{Q\cp T}U_T.
\end{equation*}
Hence a canonical distributive law $\Omega\co T Z_{Q\cp T} \to Z_{Q\cp T} T$ and a Hopf monad  $D_{Q,T}=Z_{Q\cp T} \circ_\Omega T$ on $\cc$. Moreover, we show: $$D_{Q,T}\ti\cc\simeq \zz_Q(T\ti\cc),$$ where $\zz_Q(T\ti\cc)$ is the center of $T\ti\cc$ relative to $Q$. When  $Q=\id_{T\ti\cc}$, we obtain the previous results.

This paper is organized as follows. In Section~\ref{sect-prelims}, we review several facts about monoidal categories and Hopf algebras in braided categories. Section~\ref{sect-HopfMon} recalls the definition and elementary properties of Hopf monads.  Section~\ref{sect-hopfmonadj} deals with monoidal adjunctions, exactness properties, and cross products of Hopf monads. In Section~\ref{sect-distlaw},
 we briefly recall the basic results of  Beck's theory of distributive laws and extend them to the Hopf monad setting.
In Section~\ref{sect-centHM}, we define the  centralizer $Z_T$ of a Hopf monad $T$ on $\cc$ and relate it to the center
$\zz_T(\cc)$ of $\cc$ relative to $T$. In Section~\ref{sect-double}, we define the canonical distributive law $\Omega$ of $T$ over $Z_T$ and the double $D_T=Z_T\circ_\Omega T$,  and state their categorical properties.
In Section~\ref{sect-cent-mod}, we study the centralizer $Z_Q$  of a Hopf monad $Q$ on $T\ti\cc$ and construct the canonical distributive law of $T$ over $Z_{Q\cp T}$. Section~\ref{sect-cas-hopf-ds-braided} is devoted to Hopf monads on a braided category $\bb$. In particular, we define the double $D(A)$ of a Hopf algebra $A$ in a braided autonomous category. In Section~\ref{sect-fusion-cat}, we treat the case of the center of a fusion category.

\section{Preliminaries and notations}\label{sect-prelims}

\subsection{Categories} Unless otherwise specified, categories are small, and monoidal categories
are strict.

If $\cc$ is a category, we denote $\Ob(\cc)$ the set of objects of $\cc$ and $\Hom_\cc(X,Y)$ the set of morphisms in
$\cc$ from an object $X$ to an object $Y$. The identity functor of~$\cc$ will be denoted by $1_\cc$. We denote
$\cc^\opp$ the \emph{opposite category} (where arrows are reversed).

Let $\cc$, $\dd$ be two categories. Functors from $\cc$ to $\dd$ are
the objects of a category $\Fun(\cc,\dd)$. Given two functors $F,
G\co \cc \to \dd$, a morphism $\alpha \co F \to G$ is a family
$\{\alpha_X\co F(X) \to G(X)\}_{X \in \Ob(\cc)}$ of morphisms in
$\dd$ satisfying the following functoriality condition: $\alpha_Y
F(f)=G(f)\alpha_X$ for every morphism $f\co X \to Y$ in $\cc$. Such
a morphism is called a \emph{natural transformation}. We denote
$\Nat(F,G)$ the set $\Hom_{\Fun (\cc,\dd)}(F,G)$ of natural
transformations from $F$ to $G$,   and $\id_F$ the identity natural
transformation of a functor $F$.

If $\cc$, $\cc'$ are two categories, we denote $\sigma_{\cc,\cc'}$ the flip functor $\cc \times \cc' \to \cc' \times
\cc$ defined by $(X, X') \mapsto (X', X)$.

\subsection{Monoidal categories}
Let $\cc$ be monoidal category with monoidal product $\otimes \co
\cc \times \cc \to \cc$ and unit object $\un$. For $n\geq 0$, we denote $\otimes_n$ the functor:
\begin{equation*}
\otimes_n\co\cc^n=\underbrace{\cc \times \cdots \times \cc}_{\text{$n$ times}} \to \cc, \quad (X_1,\dots,X_n)  \mapsto  X_1 \otimes \dots \otimes X_n.
\end{equation*}
Note that $\otimes_0$ is the constant functor equal to $\un$,
 $\otimes_1=\id_{1_\cc}$,
and $\otimes_2=\otimes$.

For a family of functors $\{F_i\co \aaa_i \to \cc\}_{1\leq i\leq
n}$, set:
\begin{equation*}
F_1 \otimes \cdots \otimes F_n =
\otimes_n\circ(F_1,\cdots,F_n)\co \aaa_1 \times \cdots \times \aaa_n \to \cc.
\end{equation*}
For a functor $F\co \aaa \to \cc$, set:
\begin{equation*}
F^{\otimes n}=\underbrace{F \otimes \cdots \otimes F}_{\text{$n$ times}}\co \aaa^n \to \cc.
\end{equation*}

If $\cc$ is a monoidal category,  we denote $\cc^{\otimes \opp}$
the monoidal category with opposite monoidal product $\otimes^\opp$ defined by $X \otimes^\opp Y=Y \otimes X$.

\subsection{Monoidal functors}\label{sect-monofunctor}
Let $(\cc,\otimes,\un)$ and $(\dd, \otimes, \un)$ be two monoidal categories.
A \emph{monoidal functor} from $\cc$ to $\dd$ is a triple
$(F,F_2,F_0)$, where $F\co \cc \to \dd$ is a functor, $F_2\co
F\otimes F \to F \otimes$ is a morphism of functors, and $F_0\co\un
\to F(\un)$ is a morphism in~$\dd$, such that:
\begin{align*}
& F_2(X,Y \otimes Z) (\id_{F(X)} \otimes F_2(Y,Z))= F_2(X \otimes Y, Z)(F_2(X,Y) \otimes \id_{F(Z)}) ;\\
& F_2(X,\un)(\id_{F(X)} \otimes F_0)=\id_{F(X)}=F_2(\un,X)(F_0
\otimes \id_{F(X)}) ;
\end{align*}
for all objects $X,Y,Z$ of $\cc$.

A monoidal functor $(F,F_2,F_0)$ is said to be \emph{strong} (resp.\@ \emph{strict}) if $F_2$ and $F_0$ are
isomorphisms (resp.\@ identities).

By a \emph{monoidal isomorphism}, we mean a strong monoidal functor which is an isomorphism.

\subsection{Monoidal natural transformations}
Let $F\co\cc \to \dd$ and $G\co\cc \to \dd$ be two monoidal functors.
A natural transformation $\varphi\co F \to G$ is \emph{monoidal}
if it satisfies:
\begin{equation*}
\varphi_{X \otimes Y} F_2(X,Y)= G_2(X,Y) (\varphi_X \otimes
\varphi_Y) \quad \text{and} \quad G_0=\varphi_\un F_0.
\end{equation*}

\subsection{Comonoidal functors}\label{sect-comonofunctor}
Let $(\cc,\otimes,\un)$ and $(\dd, \otimes, \un)$ be two monoidal categories.
A \emph{comonoidal functor}\footnote{Comonoidal functors are also called \emph{opmonoidal functors}} from $\cc$ to
$\dd$ is a triple $(F,F_2,F_0)$, where $F\co \cc \to \dd$ is a functor, $F_2\co F \otimes \to F\otimes F$ is a natural
transformation, and $F_0\co F(\un) \to \un$ is a morphism in $\dd$, such that:
\begin{align*}
& \bigl(\id_{F(X)} \otimes F_2(Y,Z)\bigr) F_2(X,Y \otimes Z)= \bigl(F_2(X,Y) \otimes \id_{F(Z)}\bigr) F_2(X \otimes Y, Z) ;\\
& (\id_{F(X)} \otimes F_0) F_2(X,\un)=\id_{F(X)}=(F_0 \otimes \id_{F(X)}) F_2(\un,X) ;
\end{align*}
for all objects $X,Y,Z$ of $\cc$.

It is convenient to denote $F_3\co F\otimes_3 \to F^{\otimes 3}$ the natural transformation defined by
$F_3(X,Y,Z)=(\id_{F(X)} \otimes F_2(Y,Z)) F_2(X,Y \otimes Z)= (F_2(X,Y) \otimes \id_{F(Z)}) F_2(X \otimes Y, Z)$.

A comonoidal functor $(F,F_2,F_0)$ is said to be \emph{strong} (resp.\@ \emph{strict}) if $F_2$ and $F_0$ are
isomorphisms (resp.\@ identities). In that case, $(F,F^{-1}_2,F^{-1}_0)$ is a strong (resp. strict) monoidal functor.

\subsection{Comonoidal natural transformations}
Let $F\co\cc \to \dd$ and $G\co\cc \to \dd$ be two comonoidal functors.
A natural transformation $\varphi\co F \to G$
is \emph{comonoidal} if it satisfies:
\begin{equation*}
G_2(X,Y) \varphi_{X \otimes Y}= (\varphi_X \otimes \varphi_Y) F_2(X,Y)\quad \text{and} \quad G_0 \varphi_\un= F_0.
\end{equation*}

\subsection{Autonomous categories}
Recall that a \emph{duality} in a monoidal category $\cc$ is a quadruple $(X,Y,e,d)$, where $X$, $Y$ are objects of
$\cc$, $e\co X \otimes Y \to \un$ (the \emph{evaluation}) and $d \co \un \to Y \otimes X$ (the \emph{coevaluation}) are
morphisms in $\cc$, such that:
\begin{equation*}
(e \otimes \id_X)(\id_X \otimes d)=\id_X \quad \text{and} \quad (\id_Y \otimes e)(d \otimes \id_Y)=\id_Y.
\end{equation*}
Then $(X,e,d)$ is a \emph{left dual} of $Y$ and $(Y,e,d)$ is a \emph{right dual} of $X$.

If $D=(X,Y,e,d)$ and $D'=(X',Y',e',d')$ are two dualities, two morphisms $f\co X \to X'$ and $g\co Y' \to Y$ are
\emph{in duality with respect to $D$ and $D'$} if
\begin{equation*}
e'(f \otimes \id_{Y'})=e(\id_X \otimes g) \quad \bigl (\text{or, equivalently, $(\id_{Y'} \otimes f)d=(g \otimes \id_X)
d'$}\bigr).
\end{equation*}
In that case we write $f= \ldual{g}_{D,D'}$ and $g = \rdual{f}_{D,D'}$, or simply $f=\ldual{g}$ and $g=\rdual{f}$. Note that this defines a bijection between $\Hom_\cc(X,X')$ and
$\Hom_\cc(Y',Y)$.

Left and right duals, if they exist, are essentially unique: if $(Y,e,d)$ and $(Y',e',d')$ are right duals of some
object $X$, then there exists a unique isomorphism $u\co Y \to Y'$ such that $e'=e(\id_X \otimes u^{-1})$ and $d'=(u
\otimes \id_X)d$.

A \emph{left autonomous} (resp.\@ \emph{right autonomous}, resp.\@ \emph{autonomous}) category is a monoidal category
for which every object admits a left dual (resp.\@ a right dual, resp.\@ both a left and a right dual).

Assume $\cc$ is a left autonomous category and, for each object $X$, pick a left dual $(\ldual{X},\lev_X,\lcoev_X)$.
This data defines a strong monoidal functor $\ldual{?}\co \cc^{\opp,\otimes \opp} \to \cc$.

Likewise, if $\cc$ is a right autonomous category, picking a right dual $(\rdual{X},\rev_X,\rcoev_X)$ for each object
$X$ defines a strong monoidal functor $\rdual{?}\co \cc^{\opp,\otimes \opp}\to \cc$.

Subsequently, when dealing with left or right autonomous categories, we shall always assume tacitly that left duals or
right duals have been chosen. Moreover, in formulae, we abstain from writing down the following
canonical isomorphisms:
\begin{align*}
& \ldual{?}_2(X,Y)\co \ldual{Y} \otimes \ldual{X} \to \ldual{(X \otimes Y)}, & \ldual{?}_0 \co \un \to \ldual{\un}, \\
& \rdual{?}_2(X,Y)\co \rdual{Y} \otimes \rdual{X} \to \rdual{(X \otimes Y)}, & \rdual{?}_0 \co \un \to \rdual{\un},
\end{align*} and
\begin{align*}
& (\rev_X \otimes \id_{\ldual{(\rdual{X})}})(\id_X \otimes \lcoev_{\rdual{X}}) \co X \to \ldual{(\rdual{X})},\\
&(\id_{\rdual{(\ldual{X})}} \otimes \lev_X)(\rcoev_{\ldual{X}} \otimes \id_X) \co X \to \rdual{(\ldual{X})}.
\end{align*}

\subsection{Braided categories}
Recall that a \emph{braiding} on a monoidal category $\cc$ is a natural isomorphism $\tau\co \otimes \to \otimes \sigma_{\cc,\cc}$ such that:
\begin{equation*}
\tau_{X, Y\otimes Z}=(\id_Y \otimes \tau_{X, Z})(\tau_{X, Y} \otimes \id_Z) \quad \text{and} \quad
\tau_{X\otimes Y,Z}=(\tau_{X,Z} \otimes \id_Y)(\id_X \otimes \tau_{Y,Z}).
\end{equation*}
A \emph{braided category} is a monoidal category endowed with a braiding.

The \emph{mirror of a braiding $\tau$} is the braiding $\overline{\tau}$ defined by
$\overline{\tau}_{X,Y}=\tau^{-1}_{Y,X}$.  The \emph{mirror of a braided category $\bb$} is the braided category $\bbb$ which coincides with $\bb$ as a monoidal category but is endowed with the mirror braiding.

If $\cc$ is braided with braiding $\tau$, then $\cc^{\otimes\opp}$ is braided with braiding $\tau^\mop$ defined by $\tau^\mop_{X,Y}=\tau_{Y,X}$. Note that $\tau \mapsto \tau^\mop$ is a bijection between braidings on $\cc$ and  braidings on $\cc^{\otimes\opp}$.

\subsection{Braided functors}
A \emph{braided functor} between two braided categories $\bb$ and~$\bb'$ is a strong monoidal functor $F\co \bb \to \bb'$ such that:
\begin{equation*}
F(\tau_{X,Y})F_2(X,Y)=F_2(Y,X)\tau'_{F(X),F(Y)}
\end{equation*}
for all objects $X,Y$ of $\bb$, where $\tau$ and $\tau'$ are the braidings of $\bb$ and $\bb'$.

\begin{exa}\label{rem-isobraidedop} If $\bb$ is a braided category with braiding $\tau$, the monoidal functor $(1_\bb,\tau, \id_\un) \co \bb^{\otimes \opp} \to \bb$ is a braided isomorphism.
\end{exa}

\subsection{The center of a monoidal category}\label{sect-centerusual}
Let $\cc$ be a monoidal category. A \emph{left half braiding} of $\cc$ is a pair $(M,\sigma)$, where $M$ is an object of $\cc$ and $\sigma\co M \otimes 1_\cc
\to 1_\cc \otimes M$ is a natural transformation such that:
\begin{enumerate}
\renewcommand{\labelenumi}{{\rm (\roman{enumi})}}
\item $\sigma_{Y \otimes Z}=(\id_Y \otimes \sigma_Z)(\sigma_Y \otimes \id_Z)$ for all $Y,Z \in \Ob(\cc)$;
\item $\sigma_\un=\id_M$;
\item $\sigma$ is an isomorphism.
\end{enumerate}
Note that if $\cc$ is autonomous, (iii) is a consequence of (i) and (ii).

The \emph{center of $\cc$} is the braided category $\zz(\cc)$ defined as follows. Its objects are left half braidings of $\cc$. A morphism in $\zz(\cc)$ from $(M,\sigma)$ to $(M',\sigma')$ is a morphism $f \co M \to M'$ in $\cc$ such that:
$(\id_{1_\cc} \otimes f)\sigma=\sigma'(f \otimes \id_{1_\cc})$. The monoidal product and braiding $\tau$ are:
\begin{equation*}
(M,\sigma) \otimes (N,\gamma)=\bigl(M \otimes N,(\sigma \otimes \id_N)(\id_M \otimes \gamma) \bigr) \quad \text{and} \quad \tau_{(M,\sigma),(N,\gamma)}=\sigma_{N}.
\end{equation*}
Note that if $\cc$ is autonomous, so is $\zz(\cc)$.

\begin{rem}\label{rem-alter-center}
Likewise, define a \emph{right half braiding} of a monoidal category~$\cc$ to be a pair $(M,\sigma)$, where $M$ is an object of $\cc$ and $\sigma\co 1_\cc \otimes M
\to  M \otimes 1_\cc$ is a natural transformation satisfying analogous axioms. Right half braidings form a braided category $\zz'(\cc)$, with braiding: $\tau'_{(M,\sigma),(N,\gamma)}=\gamma_{M}$.
We have:  $$\zz'(\cc)=\zz(\cc^{\otimes \opp})^{\otimes \opp}. $$
The braided category $\zz'(\cc)$ is isomorphic to the mirror of $\zz(\cc)$ via the braided isomorphism given by $(M,\sigma) \mapsto (M,\sigma^{-1})$.
\end{rem}

\subsection{Algebras, bialgebras, and Hopf algebras in categories}
Let $\cc$ be a monoidal category. An \emph{algebra in $\cc$} is an object $A$ of $\cc$ endowed with morphisms $m\co A \otimes A \to A$ (the product) and $u\co \un \to A$ (the unit) such that:
\begin{equation*}
m(m \otimes \id_A)=m(\id_A \otimes m) \quad \text{and} \quad m(\id_A \otimes u)=\id_A=m(u \otimes \id_A).
\end{equation*}
A \emph{coalgebra in $\cc$} is an object $C$ of $\cc$ endowed with morphisms $\Delta\co C \to C \otimes C$ (the coproduct) and $\varepsilon\co C \to \un$ (the counit) such that:
\begin{equation*}
(\Delta \otimes \id_C)\Delta=(\id_C \otimes \Delta)\Delta \quad \text{and} \quad (\id_C \otimes \varepsilon)\Delta=\id_C=(\varepsilon \otimes \id_C)\Delta.
\end{equation*}

Let $\bb$ be a braided category, with braiding
$\tau$. A \emph{bialgebra in $\bb$} is an object $A$ of~$\bb$ endowed with an algebra structure $(m,u)$ and a coalgebra structure $(\Delta,\varepsilon)$ in~$\bb$ satisfying:
\begin{align*}
\Delta m&=(m \otimes m)(\id_A \otimes \tau_{A,A} \otimes \id_A)(\Delta \otimes \Delta),  & \Delta u&=u \otimes u, \\  \varepsilon m&=\varepsilon \otimes \varepsilon, & \varepsilon u&=\id_\un.
\end{align*}
Let $A$ be a bialgebra in $\bb$.
Set:
\begin{equation*}
m^\opp=m\tau^{-1}_{A,A} \quad \text{and} \quad \Delta^\cop= \tau^{-1}_{A,A}\Delta.
\end{equation*}
Then $(A,m^\opp,u,\Delta,\varepsilon)$ is a bialgebra in the mirror $\bbb$ of $\bb$, called the \emph{opposite of~$A$}, and denoted
$A^\opp$. Similarly $(A,m,u,\Delta^\cop,\varepsilon)$ is a bialgebra in $\bbb$, called the \emph{co-opposite of $A$},
and denoted $A^\cop$. Consequently $A^{\cop,\opp}=(A^\cop)^\opp$ is a bialgebra in $\bb$ (with product $m\tau_{A,A}$ and coproduct $\tau^{-1}_{A,A}\Delta$).

An \emph{antipode} for a bialgebra $A$ is a morphism $S\co A \to A$ in $\bb$ such that:
\begin{equation*}
m(S \otimes \id_A)\Delta=u \varepsilon=m(\id_A \otimes S)\Delta.
\end{equation*}
If it exists, an antipode is unique, and it is a morphism of bialgebras $A\to A^{\cop,\opp}$. A \emph{Hopf algebra in $\bb$} is a bialgebra in $\bb$ which admits an invertible antipode.

If $A$ is a Hopf algebra in $\bb$, with antipode $S$, then $A^\opp$ and $A^\cop$ are Hopf algebras in the mirror $\bbb$ of $\bb$, with antipode $S^{-1}$.

\subsection{Modules in categories}\label{sect-modulcatclassicHA}
Let $(A,m,u)$ be an algebra in a monoidal category~$\cc$.
A \emph{left $A$\ti module} (in $\cc$) is a pair~$(M,r)$, where $M$ is an object of $\cc$ and $r\co A \otimes M \to M$ is a morphism in $\cc$, such that:
\begin{equation*}
r(m \otimes \id_M)=r(\id_A \otimes r) \quad \text{and} \quad r(u \otimes \id_M)=\id_M.
\end{equation*}
An \emph{$A$\ti linear morphism} between two left $A$\ti modules $(M,r)$ and $(N,s)$ is a morphism $f\co M \to N$ such that $fr=s(\id_A \otimes f)$. Hence the category $\lmod{\cc}{A}$  of left $A$\ti modules. Likewise, one defines the category $\rmod{\cc}{A}$ of right $A$\ti modules.

Let $A$ be a bialgebra in a braided category $\bb$.  Then the category $\lmod{\bb}{A}$ is monoidal, with unit object $(\un,\varepsilon)$ and monoidal product:
\begin{equation*}
(M,r) \otimes (N,s)= (r \otimes s)(\id_A \otimes \tau_{A,M} \otimes \id_N)(\Delta \otimes \id_{M \otimes N}),
\end{equation*}
where $\Delta$ and $\epsilon$ are the coproduct and counit of $A$, and $\tau$ is the braiding of $\bb$.
Likewise the category $\rmod{\cc}{A}$ is monoidal, with unit object $(\un,\varepsilon)$ and monoidal product:
\begin{equation*}
(M,r) \otimes (N,s)= (r \otimes s)(\id_M \otimes \tau_{N,A} \otimes \id_A)(\Delta \otimes \id_{M \otimes N}).
\end{equation*}

Assume $\bb$ is autonomous. Then $\lmod{\bb}{A}$ is autonomous if and only if $\rmod{\bb}{A}$ is autonomous, if and only if $A$ is a Hopf algebra.
If $A$ is a Hopf algebra, with antipode~$S$, then the duals of a left $A$-module $(M,r)$ are:
\begin{align*}
&\ldual{(M,r)}=\bigl(\ldual{M}, (\lev_M \otimes \id_{\ldual{M}}) (\id_{\ldual{M}} \otimes r(S \otimes \id_M)\otimes \id_{\ldual{M}})(\tau_{A,\ldual{M}} \otimes \lcoev_M)\bigr),\\
&\rdual{(M,r)}=\bigl(\rdual{M},  (\id_{\rdual{M}} \otimes \rev_M)(\id_{\rdual{M}} \otimes r\tau^{-1}_{A,M} \otimes \id_{\rdual{M}})(\rcoev_M \otimes S^{-1} \otimes \id_{\rdual{M}})\bigr),
\end{align*}
and the duals of a right $A$-module $(M,r)$ are:
\begin{align*}
&\ldual{(M,r)}=\bigl(\ldual{M},  (\lev_M \otimes \id_{\ldual{M}})(\id_{\ldual{M}} \otimes r\tau^{-1}_{M,A} \otimes \id_{\ldual{M}})(\id_{\ldual{M}} \otimes S^{-1} \otimes \lcoev_M) \bigr),\\
&\rdual{(M,r)}=\bigl(\rdual{M}, (\id_{\rdual{M}} \otimes \rev_M) (\id_{\rdual{M}} \otimes r(\id_M \otimes S) \otimes \id_{\rdual{M}})(\rcoev_M \otimes \tau_{\rdual{M},A}) \bigr).
\end{align*}

\begin{rem}\label{rem-braidleftright}
Let $A$ be a Hopf algebra in a braided category $\bb$, with braiding~$\tau$. The functor $F_A\co \lmod{\bb}{A} \to \rmod{\bb}{A}$, defined by $F_A(M,r)=\bigl(M,r \tau_{M,A}(\id_M \otimes S)\bigr)$ and $F_A(f)=f$, gives rise to a monoidal isomorphism of categories:
$$F_A=(F_A,\tau, \un) \co (\lmod{\bb}{A})^{\otimes\opp} \to \rmod{\bb}{A}.$$
Therefore braidings on $\lmod{\bb}{A}$ are in bijection with braidings on $\rmod{\bb}{A}$. More precisely, if $c$ is a braiding on  $\rmod{\bb}{A}$, then:
$$c'_{(M,r),(N,s)}=\tau_{M,N}\,c_{F_A(N,s),F_A(M,r)}\,\tau^{-1}_{N,M}$$ is a braiding on $\lmod{\bb}{A}$ (making $F_A$ braided), and the correspondence $c \mapsto c'$ is bijective.
\end{rem}

\subsection{Penrose graphical calculus} We represent morphisms in a category by diagrams to be read from bottom to top.   Thus we draw the identity $\id_X$ of an object $X$, a morphism $f\co X \to Y$, and its composition with a morphism $g\co Y \to Z$  as follows:
\begin{center}
\psfrag{X}[Bc][Bc]{\scalebox{.8}{$X$}} \psfrag{Y}[Bc][Bc]{\scalebox{.8}{$Y$}} \psfrag{h}[Bc][Bc]{\scalebox{.8}{$f$}} \psfrag{g}[Bc][Bc]{\scalebox{.8}{$g$}}
\psfrag{Z}[Bc][Bc]{\scalebox{.8}{$Z$}} $\id_X=$ \rsdraw{.45}{.9}{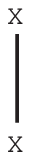}\,,\quad $f=$ \rsdraw{.45}{.9}{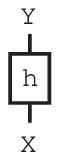} ,\quad \text{and} \quad $gf=$ \rsdraw{.45}{.9}{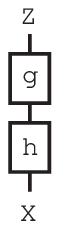}\,.
\end{center}
In a monoidal category, we represent the monoidal product of two morphisms $f\co X \to Y$ and $g \co U \to V$ by juxtaposition:
\begin{center}
\psfrag{X}[Bc][Bc]{\scalebox{.8}{$X$}} \psfrag{h}[Bc][Bc]{\scalebox{.8}{$f$}}
\psfrag{Y}[Bc][Bc]{\scalebox{.8}{$Y$}}  $f\otimes g=$ \rsdraw{.45}{.9}{morphism} \psfrag{X}[Bc][Bc]{\scalebox{.8}{$U$}} \psfrag{g}[Bc][Bc]{\scalebox{.8}{$g$}}
\psfrag{Y}[Bc][Bc]{\scalebox{.8}{$V$}} \rsdraw{.45}{.9}{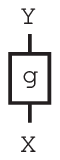}\,.
\end{center}
The duality morphisms of an autonomous category are depicted as:
\begin{center}
\psfrag{C}[Bc][Bc]{\scalebox{.8}{$X$}} \psfrag{A}[Bc][Bc]{\scalebox{.8}{$\ldual{X}$}} $\lev_X=$ \rsdraw{.45}{.9}{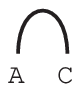}\,,\quad
\psfrag{A}[Bc][Bc]{\scalebox{.8}{$X$}} \psfrag{C}[Bc][Bc]{\scalebox{.8}{$\ldual{X}$}} $\lcoev_X=$ \rsdraw{.45}{.9}{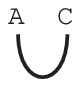}\,,\quad
\psfrag{A}[Bc][Bc]{\scalebox{.8}{$X$}} \psfrag{C}[Bc][Bc]{\scalebox{.8}{$\rdual{X}$}} $\rev_X=$ \rsdraw{.45}{.9}{eval}\,,\quad \text{and} \quad
\psfrag{C}[Bc][Bc]{\scalebox{.8}{$X$}} \psfrag{A}[Bl][Bl]{\scalebox{.8}{$\rdual{X}$}} $\rcoev_X=$ \rsdraw{.45}{.9}{coeval}\,.
\end{center}
The braiding $\tau$ of a braided category, and its inverse, are depicted as:
\begin{center}
\psfrag{X}[Bc][Bc]{\scalebox{.8}{$X$}} \psfrag{Y}[Bc][Bc]{\scalebox{.8}{$Y$}} $\tau_{X,Y}=\,$\rsdraw{.45}{.9}{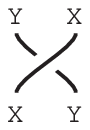} \quad \text{and} \quad $\tau^{-1}_{Y,X}=\,$\rsdraw{.45}{.9}{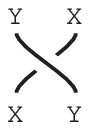}.
\end{center}

Given a Hopf algebra $A$ in a braided category, we depict its product $m$, unit $u$,
coproduct $\Delta$, counit $\varepsilon$, antipode $S$, and $S^{-1}$ as follows:
\begin{center}
\psfrag{A}[Bc][Bc]{\scalebox{.8}{$A$}} $m=$\rsdraw{.45}{.9}{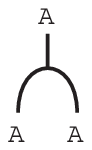}, \quad $u=$\rsdraw{.45}{.9}{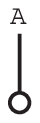}, \quad
$\Delta=$\rsdraw{.45}{.9}{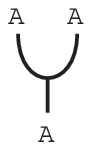}, \quad $\varepsilon=$\rsdraw{.45}{.9}{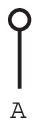}, \quad $S=\,\rsdraw{.45}{.9}{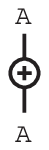}$\,,
\quad $S^{-1}=\,\rsdraw{.45}{.9}{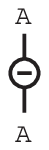}$\,.
\end{center}

\section{Hopf monads and their modules}\label{sect-HopfMon}
In this section, we review the notion of a Hopf monad. For a general treatment, we refer to \cite{BV2}.

\subsection{Monads}\label{monad}
Let $\cc$ be a category. Recall that the category $\End(\cc)$ of endofunctors of $\cc$ is strict monoidal with
composition for monoidal product and identity functor $1_\cc$ for unit object. A \emph{monad} on $\cc$ is an algebra in
$\End(\cc)$, that is, a triple $(T,\mu,\eta)$, where $T\co \cc \to \cc$ is a functor, $\mu\co T^2 \to T$ and $\eta\co
1_\cc \to T$ are natural transformations, such that:
\begin{equation*}
\mu_X T(\mu_X)=\mu_X\mu_{T(X)} \quad \text{and} \quad
\mu_X\eta_{T(X)}=\id_{T(X)}=\mu_X T(\eta_X)
\end{equation*}
for any object $X$ of $\cc$.

\begin{exa} \label{exa-MfromA} Let $A$ be an algebra in a monoidal category $\cc$, with product $m$ and unit $u$. Then the endofunctor $? \otimes A$ of $\cc$, defined by $X \mapsto X \otimes A$, has a structure
of a monad on $\cc$, with product $\mu=\id_{1_\cc} \otimes m$ and unit $\eta=\id_{1_\cc} \otimes u$.
Similarly, the endofunctor $A \otimes \,?$ is a monad on $\cc$ with product
$m \otimes \id_{1_\cc}$ and unit $u \otimes \id_{1_\cc}$.
\end{exa}

\subsection{Bimonads}
A \emph{bimonad}\footnote{Bimonads were introduced in \cite{Moer} under the name `Hopf monads', which we prefer to
reserve for bimonads with antipodes by analogy with Hopf algebras.} on a monoidal category~$\cc$ is a monad
$(T,\mu,\eta)$ on $\cc$ such that the functor $T\co \cc \to \cc$ is comonoidal and the natural transformations $\mu\co
T^2 \to T$ and $\eta\co 1_\cc \to T$ are comonoidal. In other words, $T$ is endowed with a natural
transformation $T_2\co T \otimes \to T\otimes T$ and a morphism $T_0\co T(\un) \to \un$ in $\cc$ such that:
\begin{align*}
& \bigl(\id_{T(X)} \otimes T_2(Y,Z)\bigr) T_2(X,Y \otimes Z)= \bigl(T_2(X,Y) \otimes \id_{T(Z)}\bigr) T_2(X \otimes Y, Z) ,\\
& (\id_{T(X)} \otimes T_0) T_2(X,\un)=\id_{T(X)}=(T_0 \otimes \id_{T(X)}) T_2(\un,X),
\end{align*}
and
\begin{align*}
& T_2(X,Y) \mu_{X \otimes Y}= (\mu_X \otimes \mu_Y) T_2(T(X),T(Y)) T(T_2(X,Y)),\\
& T_0 \mu_\un= T_0 T(T_0) , \qquad T_2(X,Y) \eta_{X \otimes Y}= (\eta_X \otimes \eta_Y) , \qquad T_0 \eta_\un= \id_\un.
\end{align*}

\begin{rem}\label{bimonad-alg}
A bimonad on a monoidal category $\cc$ is nothing but an algebra in
the strict monoidal category
of comonoidal endofunctors of $\cc$ (with monoidal product $\circ$ and unit object $1_\cc$).
\end{rem}

\begin{rem}\label{remoppositemonad}
A bimonad $T$ on a monoidal category $\cc=(\cc,\otimes,\un)$ may be viewed as a bimonad $T^\mop$ on the
monoidal category $\cc^{\otimes \opp}=(\cc,\otimes^\opp,\un)$, with comonoidal structure $T_2^\mop=T_2\sigma_{\cc,\cc}$ and $T_0^\mop=T_0$. The bimonad $T^\mop$  is called  the \emph{opposite} of the bimonad $T$. We have: $T^\mop
\ti\cc^{\otimes\opp}=(T\ti \cc)^{\otimes \opp}$.
\end{rem}

\subsection{Antipodes}\label{sect-antipod}
Right and left antipodes of a Hopf monad generalize the antipode of a Hopf algebra and its inverse. Let $(T,\mu,\eta)$ be a bimonad on a monoidal category~$\cc$.

Assume $\cc$ is left autonomous.  A \emph{left antipode for $T$} is
a natural transformation $s^l=\{s^l_X\co T(\ldual{T(X)}) \to
\ldual{X}\}_{X \in \Ob(\cc)}$ satisfying:
\begin{align*}
& T_0 T(\lev_X)T(\ldual{\eta_X} \otimes \id_X)=\lev_{T(X)}(s^l_{T(X)}T(\ldual{\mu}_X) \otimes
\id_{T(X)})T_2(\ldual{T(X)},X); \\
& (\eta_X \otimes \id_{\ldual{X}})\lcoev_X T_0=(\mu_X \otimes s^l_X) T_2(T(X),\ldual{T(X)})T(\lcoev_{T(X)});
\end{align*}
for every object $X$ of $\cc$. By \cite[Theorem 3.7]{BV2}, a left antipode $s^l$ is `anti-(co)multi\-pli\-ca\-ti\-ve':
for all objects $X,Y$ of $\cc$,
\begin{align*}
&s^l_X \mu_{\ldual{T(X)}}=s^l_X T(s^l_{T(X)}) T^2(\ldual{\mu_X}); &
&s^l_X\eta_{\ldual{T(X)}}=\ldual{\eta_X};\\
&s^l_{X\otimes Y}T(\ldual{T_2(X,Y)})=(s^l_Y \otimes s^l_X)T_2(\ldual{T(Y)},\ldual{T(X)});& &s^l_\un T(\ldual{T_0})=T_0.
\end{align*}

Assume $\cc$ is right autonomous.  A \emph{right antipode for $T$}
is a natural transformation $s^r=\{s^r_X\co T(\rdual{T(X)}) \to
\rdual X\}_{X \in \Ob(\cc)}$ satisfying:
\begin{align*}
& T_0 T(\rev_X)T(\id_X \otimes \eta_X^\vee)=\rev_{T(X)}(\id_{T(X)}
\otimes s^r_{T(X)}T(\mu_X^\vee))T_2(X,\rdual{T(X)});\\
& (\id_{X^\vee}\otimes \eta_X )\rcoev_X T_0=(s^r_X \otimes\mu_X) T_2(\rdual{T(X)},T(X))T(\rcoev_{T(X)});
\end{align*}
for every object $X$ of $\cc$. By \cite[Theorem 3.7]{BV2}, a right antipode $s^r$ is also  `anti-(co)multi\-pli\-ca\-ti\-ve': for all objects $X,Y$ of
$\cc$,
\begin{align*}
&s^r_X \mu_{\rdual{T(X)}}=s^r_X T(s^r_{T(X)}) T^2(\rdual{\mu_X}); &
&s^r_X \eta_{\rdual{T(X)}}=\rdual{\eta_X};\\
&s^r_{X \otimes Y}T(\rdual{T_2(X,Y)})=(s^r_Y \otimes s^r_X)T_2(\rdual{T(Y)},\rdual{T(X)}); & &s^r_\un
T(\rdual{T_0})=T_0.
\end{align*}

Note that if a left (resp.\@ right) antipode exists, then it is unique. Furthermore, when both exist, the left antipode
$s^l$ and the right antipode $s^r$ are `inverse' to each other in the sense that $ \id_{T(X)}=s^r_{\ldual{T(X)}}
T(\rdual{(s_X^l)})=s^l_{\rdual{T(X)}} T(\ldual{(s_X^r)})$ for any object $X$ of $\cc$.

\subsection{Hopf monads}\label{sect-hopfmon}
A \emph{Hopf monad} is a bimonad on an autonomous category which has a left antipode and a right antipode.

Hopf monads generalize Hopf algebras in a non-braided setting. In particular, finite-dimensional Hopf algebras and
several generalizations (Hopf algebras in braided autonomous categories, bialgebroids, etc...) provide examples
of Hopf monads. If fact, any monoidal adjunction between autonomous categories gives rise to a Hopf monad (see
Theorem~\ref{thm-gen-ex}). It turns out that much of the theory of finite-dimensional Hopf algebras (such as the
decomposition of Hopf modules, the existence of integrals, Maschke's criterium of semisimplicity, etc...) extends to
Hopf monads, see \cite{BV2}.

\begin{exa}[Hopf monads associated with Hopf algebras]\label{exa-HMfromHA}  Let $A$ be a Hopf algebra in a braided autonomous category $\bb$, with braiding
$\tau$.
According to \cite{BV2}, the endofunctor $? \otimes A$ of $\bb$
has a structure of a Hopf monad
on $\bb$, with product  $\mu=\id_{1_\bb} \otimes m$, unit
$\eta=\id_{1_\bb}\otimes u$, comonoidal structure given by:
\begin{equation*}
(? \otimes A)_2(X,Y)=(\id_X \otimes \tau_{Y,A} \otimes \id_A)(\id_{X \otimes Y} \otimes \Delta) \quad \text{and} \quad
(?\otimes A)_0=\varepsilon,
\end{equation*}
and left and right antipodes:
\begin{align*}
&s^l_X=( \lev_A \otimes  \id_{\ldual{X}})( \id_{\ldual{A}} \otimes \tau_{\ldual{X},A} )( \id_{\ldual{A}\otimes\ldual{X}} \otimes S^{-1} ),\\
&s^r_X=( \rev_A\otimes \id_{\rdual{X}})\tau_{\rdual{A}\otimes\rdual{X},A}( \id_{\rdual{A}\otimes\rdual{X}} \otimes S ).
\end{align*}
Pictorially, the structural morphisms of $? \otimes A$  are:
\begin{center}
\psfrag{A}[Bc][Bc]{\scalebox{.8}{$A$}} \psfrag{X}[Bc][Bc]{\scalebox{.8}{$X$}} \psfrag{Y}[Bc][Bc]{\scalebox{.8}{$Y$}}
  $\mu_X=\rsdraw{.45}{.9}{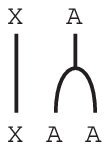}$, \qquad $\eta_X=\rsdraw{.45}{.9}{uA}$, \qquad $(?\otimes A)_2(X,Y)=\rsdraw{.45}{.9}{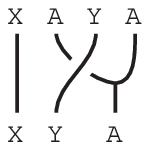}$,
\\[.6em]
$(?\otimes A)_0=\rsdraw{.45}{.9}{epsA}$,\qquad
   \psfrag{B}[Bl][Bl]{\scalebox{.8}{$\ldual{A}$}}\psfrag{X}[Bc][Bc]{\scalebox{.8}{$\ldual{X}$}} $s^l_X=\rsdraw{.45}{.9}{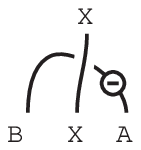}$\,,
\qquad \psfrag{B}[Bc][Bc]{\scalebox{.8}{$\rdual{A}$}}\psfrag{X}[Bl][Bl]{\scalebox{.8}{$\rdual{X}$}}
$s^r_X=\rsdraw{.45}{.9}{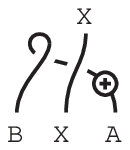}$\,.
\end{center}
Similarly, the endofunctor
$A \otimes ?$ of $\bb$ has a structure of
a Hopf monad on $\bb$, with product $\mu= m \otimes\id_{1_\bb}$,
unit $\eta=u\otimes \id_{1_\bb}$, comonoidal structure:
\begin{equation*}
(A \otimes ?)_2(X,Y)=(\id_A \otimes \tau_{A,X} \otimes \id_Y)(\Delta \otimes \id_{X \otimes Y}) \quad \text{and} \quad
(A \otimes ?)_0=\varepsilon,
\end{equation*}
and left and right antipodes:
\begin{align*}
&s^l_X=(\id_{\ldual{X}}\otimes \lev_A)\tau_{A,\ldual{X}\otimes\ldual{A}}(S \otimes \id_{\ldual{X}\otimes\ldual{A}}),\\ &s^r_X=(\id_{\rdual{X}} \otimes \rev_A)(\tau_{A,\rdual{X}} \otimes \id_{\ldual{A}})(S^{-1} \otimes \id_{\rdual{X}\otimes\rdual{A}}).
\end{align*}
Pictorially, the structural morphisms of $A \otimes ?$ are:
\begin{center}
\psfrag{A}[Bc][Bc]{\scalebox{.8}{$A$}} \psfrag{X}[Bc][Bc]{\scalebox{.8}{$X$}} \psfrag{Y}[Bc][Bc]{\scalebox{.8}{$Y$}}
  $\mu_X=\rsdraw{.45}{.9}{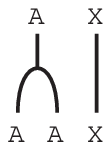}$, \qquad $\eta_X=\rsdraw{.45}{.9}{uA}$, \qquad $(A \otimes ?)_2(X,Y)=\rsdraw{.45}{.9}{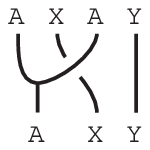}$,
\\[.6em]
$(A \otimes ?)_0=\rsdraw{.45}{.9}{epsA}$,\qquad
   \psfrag{B}[Bl][Bl]{\scalebox{.8}{$\ldual{A}$}}\psfrag{X}[Bc][Bc]{\scalebox{.8}{$\ldual{X}$}} $s^l_X=\rsdraw{.45}{.9}{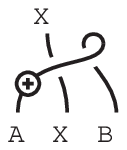}$\,,
\qquad \psfrag{B}[Bc][Bc]{\scalebox{.8}{$\rdual{A}$}}\psfrag{X}[Bl][Bl]{\scalebox{.8}{$\rdual{X}$}}
$s^r_X=\rsdraw{.45}{.9}{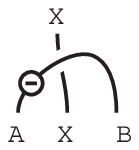}$.
\end{center}\end{exa}

\begin{exa} \label{exa-HMfromHAcent} The previous example can be extended to the non-braided setting as follows.
Let $\cc$ be a autonomous category  and $(A,\sigma)$ be a Hopf algebra in the center $\zz(\cc)$ of $\cc$ (see Section~\ref{sect-centerusual}). Denote
$m$, $u$, $\Delta$, $\varepsilon$, $S$ the product, unit, coproduct, counit, and antipode of $(A,\sigma)$. Observe that $(A,m,u)$ is an algebra in $\cc$.
Then the endofunctor $A \otimes ?$ of $\cc$
has a structure of a Hopf monad
on $\cc$,  denoted $A \otimes_\sigma ?$,  with product $\mu= m \otimes\id_{1_\cc}$,
unit $\eta=u\otimes \id_{1_\cc}$, comonoidal structure:
\begin{equation*}
(A \otimes_\sigma ?)_2(X,Y)=(\id_A \otimes \sigma_X \otimes \id_Y)(\Delta \otimes \id_{X \otimes Y})
\quad \text{and} \quad
(A \otimes_\sigma ?)_0=\varepsilon,
\end{equation*}
and left and right antipodes:
\begin{align*}
&s^l_X=(\id_{\ldual{X}}\otimes \lev_A)\sigma_{\ldual{X}\otimes\ldual{A}}(S \otimes \id_{\ldual{X}\otimes\ldual{A}}),\\ &s^r_X=(\id_{\rdual{X}} \otimes \rev_A)(\sigma_{\rdual{X}} \otimes \id_{\ldual{A}})(S^{-1} \otimes \id_{\rdual{X}\otimes\rdual{A}}).
\end{align*}
Likewise, if $(A,\sigma)$ is a Hopf algebra in $\zz'(\cc)$ (see Remark~\ref{rem-alter-center}), then
the endofunctor $? \otimes A$ of $\cc$ has a structure of a Hopf monad
on $\cc$,  denoted $? \otimes_\sigma A$,  with product  $\mu=\id_{1_\cc} \otimes m$, unit
$\eta=\id_{1_\cc}\otimes u$, comonoidal structure given by:
\begin{equation*}
(? \otimes_\sigma A)_2(X,Y)=(\id_X \otimes \sigma_Y \otimes \id_A)(\id_{X \otimes Y} \otimes \Delta) \quad \text{and} \quad
(?\otimes_\sigma A)_0=\varepsilon,
\end{equation*}
and left and right antipodes:
\begin{align*}
&s^l_X=( \lev_A \otimes  \id_{\ldual{X}})( \id_{\ldual{A}} \otimes \sigma_{\ldual{X}} )( \id_{\ldual{A}\otimes\ldual{X}} \otimes S^{-1} ),\\
&s^r_X=(\rev_A \otimes \id_{\rdual{X}})\sigma_{\rdual{A}\otimes\rdual{X}}( \id_{\rdual{A}\otimes\rdual{X}} \otimes S ).
\end{align*}
Note that if $A$ is a Hopf algebra in an autonomous braided category $\bb$ with braiding~$\tau$, then $(A,\tau_{A,-})$ is a Hopf algebra in $\zz(\bb)$, $(A,\tau_{-,A})$ is a Hopf algebra in $\zz'(\bb)$, and we have
$A \otimes ?= A \otimes_{\tau_{A,-}} ?$ and $? \otimes A= ? \otimes_{\tau_{-,A}} A$ as Hopf monads  on $\bb$.
\end{exa}

\subsection{Modules over a monad}
Let $(T,\mu,\eta)$ be a monad on a category $\cc$. An \emph{action} of $T$ on an object $M$ of $\cc$ is a morphism
$r\co T(M) \to M$ in $\cc$ such that:
\begin{equation*}
r T(r)= r \mu_M \quad \text{and} \quad r \eta_M= \id_M.
\end{equation*}
The pair $(M,r)$ is then called a \emph{$T$\ti module in $\cc$}, or just a \emph{$T$-module}\footnote{Pairs $(M,r)$ are
usually called $T$-algebras in the literature (see \cite{ML1}). However, throughout this paper, pairs $(M,r)$ are
considered as the analogues of modules over an algebra, and so the term `algebra' would be awkward in this context.}.

Given two $T$-modules $(M,r)$ and $(N,s)$ in $\cc$, a \emph{morphism of $T$\ti modules} from $(M,r)$ to $(N,r)$ is a
morphism $f\in \Hom_\cc(M,N)$ which is \emph{$T$-linear}, that is, such that $f r=s T(f)$. This gives rise to the
\emph{category $T\ti\cc$ of $T$-modules (in $\cc$)}, with composition inherited from $\cc$. We denote by $U_T\co
T\ti\cc \to \cc$ the \emph{forgetful functor of $T$} defined by $U_T(M,r)=M$ for any $T$-module $(M,r)$ and $U_T(f)=f$
for any $T$-linear morphism~$f$.

\begin{exa}\label{exa-modulesA}  Let $A$ be an algebra in
a monoidal category $\cc$ and consider the monads $?\otimes A$ and
$A \otimes ?$ of Example~\ref{exa-MfromA}. Then the category of $(?
\otimes A)$\ti modules (resp.\@ of $(A \otimes ?)$\ti modu\-les)
coincides with the category $\rmod{\cc}{A}$ of right $A$-modules
in~$\cc$ (resp.\@~with the category $\lmod{\cc}{A}$ of left
$A$-modules in $\cc$):
\begin{equation*}
(? \otimes A)\ti \cc =\rmod{\cc}{A} \quad \text{and} \quad (A \otimes ?)\ti \cc =\lmod{\cc}{A}.
\end{equation*}
\end{exa}

\subsection{Tannaka dictionary}
Structures of bimonad and Hopf monad on a monad~$T$ have natural interpretations in terms
of
the category of $T$\ti modules:
\begin{thm}[\cite{BV2}]\label{thm-biHopfmon}
Let $T$ be a monad on a monoidal category $\cc$ and $T\ti\cc$ be the category of $T$-modules.
\begin{enumerate}
\renewcommand{\labelenumi}{{\rm (\alph{enumi})}}
\item If~$T$ is a bimonad, then the category $T\ti\cc$ is
monoidal by setting:
\begin{equation*}
(M,r) \otimes (N,s)=\bigl(M \otimes N, (r \otimes s) T_2(M,N)\bigr) \quad \text{and} \quad \un_{T\ti\cc}=(\un,T_0).
\end{equation*}
Moreover this gives a bijective correspondence between bimonad structures on the monad $T$ and monoidal structures on $T\ti\cc$ such that the forgetful functor $U_T \co T\ti\cc \to \cc$ is strict monoidal.
\item Assume $T$ is a bimonad and $\cc$ is left autonomous (resp.\@ right) autonomous. Then $T$ has a left
(resp.\@ right) antipode if and only if $T\ti\cc$ is left (resp.\@ right) autonomous. If $s^l$ is a left antipode
for~$T$, left duals in $T\ti\cc$ are given by:
\begin{equation*}
\ldual{(M,r)}=(\ldual{M}, s^l_M T(\ldual{r})),\quad \lev_{(M,r)}=\lev_M, \quad\lcoev_{(M,r)}=\lcoev_M,
\end{equation*}
and if $s^r$ is a right antipode for $T$, right duals in $T\ti\cc$ are given by:
\begin{equation*}
\rdual{(M,r)}=(\rdual{M}, s^r_M T(\rdual{r})),\quad \rev_{(M,r)}=\rev_M, \quad\rcoev_{(M,r)}=\rcoev_M.
\end{equation*}
\item Assume $T$ is a bimonad and $\cc$ is autonomous. Then $T$ is a
Hopf monad if and only if $T\ti\cc$ is autonomous.
\end{enumerate}
\end{thm}

\begin{exa}\label{exa-modulesHA}  Let $A$ be a Hopf algebra in
a braided autonomous category $\bb$ and consider the Hopf monads $?\otimes A$ and $A \otimes ?$ of
Example~\ref{exa-HMfromHA}. Then:
\begin{equation*}
(? \otimes A)\ti \bb =\rmod{\bb}{A} \quad \text{and} \quad (A \otimes ?)\ti \bb =\lmod{\bb}{A}
\end{equation*}
as monoidal categories.
\end{exa}

\begin{exa}\label{exa-catmodAsigma}
More generally, let $\cc$ be a monoidal category and  $(A,\sigma)$ be a Hopf algebra in
the braided category $\zz(\cc)$.  Then $\lmod{\cc}{A}$ coincides with the category of modules over the  Hopf monad $A \otimes_\sigma ?$ on $\cc$ defined in Example~\ref{exa-HMfromHAcent}.  Hence $\lmod{\cc}{A}$ is autonomous,
with unit object $(\un,\varepsilon)$ and monoidal product:
\begin{equation*}
(M,r) \otimes (N,s)= (r \otimes s)(\id_A \otimes \sigma_M \otimes \id_N)(\Delta \otimes \id_{M \otimes N}).
\end{equation*}
Likewise, if $(A,\sigma)$ is a Hopf algebra in the braided category $\zz'(\cc)$  (see Remark~\ref{rem-alter-center}), then $\rmod{\cc}{A}$ coincides with the category of modules over the Hopf monad $? \otimes_\sigma A$,  and so is autonomous, with unit object $(\un,\varepsilon)$ and monoidal product:
\begin{equation*}
(M,r) \otimes (N,s)= (r \otimes s)(\id_M \otimes \sigma_N \otimes \id_A)(\id_{M \otimes N} \otimes \Delta).
\end{equation*}
\end{exa}

\subsection{Quasitriangular Hopf monads}\label{sect-QTHM}
A \emph{\Rt matrix} for a Hopf monad $(T,\mu,\eta)$ on an autonomous category $\cc$ is a natural transformation
\begin{equation*}
R=\{R_{X,Y}\co X\otimes Y \to T(Y) \otimes T(X)\}_{X,Y \in \Ob(\cc)}
\end{equation*}
such that,  for all objects $X,Y,Z$ of $\cc$,
\begin{align*}
& (\mu_Y \otimes \mu_X)R_{T(X),T(Y)}T_2(X,Y)=(\mu_Y \otimes \mu_X)T_2(T(Y),T(X))T(R_{X,Y}); \\
\begin{split}
& (\id_{T(Z)} \otimes T_2(X,Y))R_{X \otimes Y,Z} \\
& \phantom{XXXXXX}=(\mu_Z \otimes \id_{T(X) \otimes T(Y)}) (R_{X,T(Z)} \otimes \id_{T(Y)} )(\id_X \otimes R_{Y,Z}) ;
\end{split}\\
\begin{split}
& (T_2(Y,Z) \otimes \id_{T(X)})R_{X,Y \otimes Z} \\
& \phantom{XXXXXX}=( \id_{T(Y) \otimes T(Z)} \otimes \mu_X) ( \id_{T(Y)} \otimes R_{T(X),Z} ) (R_{X,Y} \otimes \id_Z);\\
& (\id_{T(X)} \otimes T_0) R_{\un,X}=\eta_X=(T_0 \otimes \id_{T(X)}) R_{X,\un}.
\end{split}
\end{align*}
A \emph{quasitriangular Hopf monad} is a Hopf monad equipped with an \Rt matrix.

\begin{rem}
For a bimonad, an \Rt matrix is also required to be $*$-invertible (see \cite[Section 8.2]{BV2}). For a Hopf monad $T$,
this condition is automatic and we have:
\begin{align*}
R^{*-1}_{X,Y} &= \bigl (\id_{T(X) \otimes T(Y)}\otimes\lev_X (s^l_X \otimes \id_X)\bigr )\\
& \phantom{XXX} (\id_{T(X)} \otimes
R_{\ldual{T(X)},Y} \otimes \id_X ) (\lcoev_{T(X)} \otimes\id_{Y \otimes X});\\
& =\bigl(\rev_Y (\id_Y \otimes s^r_Y ) \otimes\id_{T(X) \otimes T(Y)}\bigr )\\
& \phantom{XXX} (\id_Y \otimes R_{X,\rdual{T(Y)}} \otimes \id_{T(Y)} ) (\id_{Y \otimes X} \otimes\rcoev_{T(Y)});
\end{align*}
where $s^l$ and $s^r$ are the left and right antipodes of $T$.
\end{rem}

There is a natural interpretation of \Rt matrices for a Hopf monad $T$ in terms of braidings  on
the category of $T$\ti modules:
\begin{thm}[\cite{BV2}]\label{thm-quasitrigHopfmon}
Let $T$ be a Hopf monad on an autonomous category $\cc$.
Then any \Rt matrix $R$ for $T$ defines a braiding $\tau$ on the category $T\ti\cc$ as follows:
\begin{equation*}
\tau_{(M,r),(N,s)}=(s \otimes t)R_{M,N}\co (M,r) \otimes (N,s) \to (N,s) \otimes (M,r).
\end{equation*}
This assignment is a bijection between \Rt matrices for~$T$ and braidings on $T\ti\cc$.
\end{thm}

\begin{rem}\label{rem-quasitrigHA}
In Section~\ref{sect-RmatrixA}, we define \Rt matrices for a Hopf algebra $A$ in a braided autonomous category~$\bb$ admitting a coend $C$. These \Rt matrices are morphisms $\mathfrak{r} \co C \otimes C \to A \otimes A$, which encode \Rt matrices for the Hopf monads $? \otimes A$ and $A \otimes ?$. They generalize usual \Rt matrices
for finite-dimensional Hopf algebras.
\end{rem}

\subsection{Morphisms of Hopf monads} \label{sect-mor-mon}
A \emph{morphism of monads} between two monads $(T,\mu,\eta)$ and $(T',\mu',\eta')$ on a category $\cc$ is a
natural
transformation $f \co T \to T'$ such that, for every object $X$ of $\cc$,
\begin{equation*}
f_X \mu_X=\mu'_X f_{T'(X)}T(f_X) \quad \text{and} \quad f_X \eta_X=\eta'_X.
\end{equation*}

According to \cite[Lemma~1.7]{BV2}, a morphism of monads $f\co T \to T'$ yields a functor $f^* \co T'\ti\cc \to
T\ti\cc$ defined by $f^*(M,r)=(M, r f_M)$. Moreover, the mapping $f\mapsto f^*$ is a bijective correspondence between:
(i) morphisms of monads $f\co T \to T'$, and (ii) functors $F\co T'\ti\cc \to T\ti\cc$ such that $U_T F=U_{T'}$.

A \emph{morphism of bimonads} between two bimonads $T$ and $T'$ on a monoidal category $\cc$ is a morphism of monads
$f\co T \to T'$ which is comonoidal, that is:
\begin{equation*}
T'_2(X,Y)f_{X\otimes Y}=(f_X\otimes f_Y)T_2(X,Y) \quad \text{and} \quad T'_0f_\un =T_0.
\end{equation*}
According to~\cite[Lemma~2.9]{BV2}, the associated functor $f^*\co T'\ti\cc\to\cc \to
T\ti\cc$ is then monoidal strict. Moreover, the mapping $f\mapsto f^*$ is a bijective correspondence between: (i)
morphisms of bimonads $f\co T \to T'$, and (ii) monoidal functors $F\co T'\ti\cc \to T\ti\cc$ such that $U_T
F=U_{T'}$ as monoidal functors.

A \emph{morphism of Hopf monads} is a morphism of bimonads between Hopf monads.

\begin{exa}\label{exa-morphismHA}  Let $A$ be a Hopf algebra in a braided autonomous category $\bb$, with braiding
$\tau$.  Recall that $A^\opp$ is a Hopf algebra in the mirror $\bbb$ of $\bb$. The Hopf monad~$? \otimes A^\opp$
on $\bbb$ may be seen as a Hopf monad on $\bb$. Then
\begin{equation*}
\tau_{A, ?} \co A \otimes ? \to \, ? \otimes A^\opp
\end{equation*}
is an isomorphism of Hopf monads and
\begin{equation*}
 (\tau_{A, ?})^* \co \bbb_{A^\opp}=(? \otimes A^\opp)\ti \bb
 \to (A \otimes ?)\ti \bb=\lmod{\bb}{A}
\end{equation*}
is an isomorphism of monoidal categories. Likewise, since
$(A^\opp)^\opp=A$ as Hopf algebras in $\bb$,  $\tau_{?, A}$ induces
isomorphisms  $? \otimes A \to
A^\opp \otimes ?$ and
$\lmod{\bbb}{{A^\opp}}\to
\rmod{\bb}{A}$.
\end{exa}

\section{Hopf monads, monoidal adjunctions, and coends}\label{sect-hopfmonadj}
Monads and adjunctions are closely related. This relationship
extends naturally to Hopf monads and monoidal adjunctions between
autonomous categories.  We show that the forgetful functor of a Hopf monad creates and preserves coends. Lastly, we define the pushforward of a Hopf monad
under an adjunction and, as a special case, the cross product of
Hopf monads.

\subsection{Adjunctions}\label{sect-HM-adj} Let $\cc$ and
$\dd$ be categories. Recall that an \emph{adjunction} is a pair of functors $(F\co \cc \to \dd, U\co \dd \to \cc)$
endowed with a bijection:
\begin{equation*}
\Hom_\dd\bigl (F(X),Y\bigr ) \simeq \Hom_\cc\bigl(X,U(Y)\bigr)
\end{equation*}
which is natural in both $X \in \Ob(\cc)$ and $Y\in \Ob(\dd)$. The functor $F$ is then called \emph{left adjoint of
$U$} and the functor $U$ \emph{right adjoint of $F$}. Note that a left (resp.\@ right) adjoint of a given functor, if
it exists, is unique up to unique natural isomorphism.

An adjunction $(F,U)$ is entirely determined by two natural transformations $\eta \co 1_\cc \to UF$ and $\varepsilon\co
FU \to 1_\dd$ satisfying:
\begin{equation*}
U(\varepsilon)\,\eta_U=\id_U \quad \text{and} \quad \varepsilon_F F(\eta)=\id_F.
\end{equation*}
These transformations $\eta$ and $\varepsilon$ are respectively called the \emph{unit} and \emph{counit} of the
adjunction, and collectively the \emph{adjunction morphisms}.

Adjunctions may be composed: given two adjunctions  $(F\co \cc \to \dd, U\co \dd \to \cc)$ and $(F'\co \dd \to \ee,
U'\co \ee \to \dd)$, the pair $(F'F\co \cc \to \ee, UU'\co \ee \to \cc)$ is an adjunction called the \emph{composite}
of $(F,U)$ and $(F',U')$.

Adjunctions and monads are closely related.
Indeed if $T$ is a monad on a category~$\cc$, then the forgetful functor $U_T\co T\ti\cc \to \cc$ has a left adjoint
$F_T \co \cc \to T\ti\cc$, defined by $F_T(X)=(T(X),\mu_X)$ for any object $X$ of $\cc$ and $F_T(f)=T(f)$ for any
morphism $f$ in $\cc$. The unit of the adjunction $(F_T,U_T)$ is the unit $\eta\co 1_\cc \to T=U_TF_T$
of the monad $T$, and the counit $\varepsilon\co F_T U_T \to 1_{T\ti\cc}$ of  $(F_T,U_T)$ is the $T$\ti action, that is,  $\varepsilon_{(M,r)}=r$ for any
$T$-module $(M,r)$.

Moreover if $(F\co \cc \to \dd, U\co \dd \to \cc)$ is a pair of adjoint functors, with adjunction morphisms $\eta \co
1_\cc \to UF$ and $\varepsilon \co FU \to 1_\dd$, then $T= UF$ is a monad on~$\cc$, with product $\mu=U (\varepsilon_F)
\co T^2 \to T$ and unit $\eta$. The monad $(T,\mu,\eta)$ is the called the \emph{monad of the adjunction $(F,U)$}. In
addition there exists a unique functor $K \co \dd \to T\ti\cc$ such that $U_T K = U$ and $KF=F_T$. The functor $K$ is
called the \emph{comparison functor} and is given by $K(D)= \bigl(U(D), U(\varepsilon_D)\bigr)$ for any object $D$ of
$\dd$.

Note that if $T$ is a monad on $\cc$, then $T$ is the monad of the adjunction $(F_T,U_T)$ and the comparison functor is
the identity functor. In general, however, the comparison functor of an adjunction need not be an equivalence.

\subsection{Monadic adjunctions} An adjunction is
\emph{mona\-dic} if its comparison functor (see Section~\ref{sect-HM-adj}) is an equivalence. Remark that the composite
adjunction of two monadic adjunctions need not be monadic.

A functor $U\co \dd \to \cc$ is \emph{monadic} if it admits a left adjoint $F\co \cc \to \dd$ and the
adjunction $(F,U)$ is monadic. If such is the case, the monad $T=UF$ of the adjunction $(F,U)$ is called \emph{the
monad of $U$}. It is well-defined up to unique isomorphism of monads (as the left adjoint $F$ is unique up to unique
natural isomorphism).

For example, if $T$ is a monad on a category~$\cc$, the forgetful functor $U_T \co T\ti \cc \to \cc$ is monadic with
monad $T$.

\begin{rem}\label{rem-def-monadic}
Let $U \co \dd \to \cc$ be a functor. If there exist a monad $T$ on $\cc$ and an isomorphism of categories $K\co\dd \to T\ti\cc$ such that $U=U_TK$, then $F=K^{-1}F_T$ is left adjoint to
$U$ and the adjunction $(F,U)$ is monadic with monad $T$ and comparison functor~$K$.
\end{rem}

\subsection{Hopf monads and monoidal adjunctions} Let $\cc$ and $\dd$ be monoidal categories. An adjunction
$(F\co \cc \to \dd, U\co \dd \to \cc)$ is said to be \emph{monoidal} if the right adjoint $U \co \dd \to \cc$ is strong
monoidal. For example, if $T$ is a bimonad on a monoidal category $\cc$, then the adjunction $(F_T,U_T)$ is monoidal.

The monad of a monoidal adjunction between monoidal categories
(resp.\@ autonomous categories) is a bimonad (resp. a Hopf monad). More precisely:
\begin{thm}[\cite{BV2}]\label{thm-gen-ex}
Let $(F\co \cc \to \dd, U\co \dd \to \cc)$ be a monoidal adjunction between monoidal categories. Denote $T=UF$ the
monad of this adjunction. Then the functor $F$ is comonoidal and $T$ is a bimonad on $\cc$. The comparison
functor $K \co \dd \to T\ti \cc$ is strong monoidal, satisfies $U_T K=U$ as monoidal functors, and $KF=F_T$ as
comonoidal functors. If the categories $\cc$ and~$\dd$ are furthermore autonomous, then the
bimonad $T$ is a Hopf monad.
\end{thm}

\begin{rem}\label{rem-explicit-ex-gen}
Let $(F, U)$ be a monoidal adjunction between autonomous categories, with unit $\eta$ and counit $\varepsilon$. Let $T=UF$ be the Hopf monad associated with this monoidal adjunction (see Theorem~\ref{thm-gen-ex}). Then the comonoidal structure and antipodes of $T$ are:
\begin{align*}
&T_2(X,Y)=U_2(F(X),F(Y))^{-1} \, U(\varepsilon_{F(X) \otimes F(Y)}) \, UF\bigl (U_2(F(X),F(Y))(\eta_X \otimes \eta_Y) \bigr),\\
& T_0=U_0^{-1} \, U(\varepsilon_\un) \, UF(U_0),\\
& s^l_X=\ldual{\eta}_X  \, U_1^l(F(X))^{-1}\, U(\varepsilon_{\ldual{F}(X)}) \, UF\bigl(U^l_1(F(X))\bigr),\\
& s^r_X=\eta_X^\vee \, U_1^r(F(X))^{-1}\, U(\varepsilon_{F(X)^\vee}) \, UF\bigl(U^r_1(F(X))\bigr),
\end{align*}
where $U^l_1(Y)\co\ldual{U}(Y) \to U(\ldual{Y})$ and $U^r_1(Y)\co U(Y)^\vee \to U(\rdual{Y})$ are the compatibility isomorphisms of $U$ with duals (see \cite[Section 3.2]{BV2}).
\end{rem}

\subsection{Hopf monads and right adjoints}

If $F\co \cc \to\dd$ is a functor between autonomous categories, denote $\rexcla{F}\co \cc \to \dd$ the functor defined
by: $\rexcla{F}(X)=F(\ldual{X}\rdual{)}$ and $\rexcla{F}(f)=F(\ldual{f}\rdual{)}$ for all object $X$ and morphism $f$
in $\cc$.

\begin{lem}\label{lem-right-adj}
Let $U \co \dd \to \cc$ be a strong monoidal functor between autonomous categories. If $F\co \cc \to \dd$ is a left
adjoint for $U$, then $\rexcla{F}$ is a right adjoint for $U$.
\end{lem}
\begin{proof}  Since $U$ is strong monoidal, we have  $U(\ldual{X}) \simeq \ldual{U(X)}$  for any objet $X$ of $\cc$.
Hence  the following isomorphisms:
\begin{align*}
\Hom_\cc\bigl(U&(X),Y\bigr) \simeq \Hom_\cc\bigl(\ldual{Y},\ldual{U}(X)\bigr)\simeq \Hom_\cc\bigl(\ldual{Y},U(\ldual{X})\bigr) \\
& \simeq \Hom_\dd\bigl(F(\ldual{Y}),\ldual{X}\bigr) \simeq \Hom_\dd\bigl(X, F( \ldual{Y}\rdual{)} \bigr) =
\Hom_\dd\bigl(X, \rexcla{F} (Y)\bigr)
\end{align*}
which are natural in both $X \in \Ob(\cc)$ and $Y \in \Ob(\dd)$.
\end{proof}

\begin{prop}\label{prop-right-adj}
Let $T$ be a Hopf monad on an autonomous category $\cc$. Then:
\begin{enumerate}
\renewcommand{\labelenumi}{{\rm (\alph{enumi})}}
\item The endofunctor $\rexcla{T}$ of $\cc$ is a right adjoint of $T$.
\item The
functor $F^!_T\co \cc \to T\ti\cc$ is a right adjoint of the forgetful functor $U_T$.
\end{enumerate}
\end{prop}
\begin{proof}
Part (a) is \cite[Corollary~3.12]{BV2}. Part (b) is Lemma~\ref{lem-right-adj} applied to the monoidal adjunction
$(F_T,U_T)$.
\end{proof}

\begin{rem}\label{radj-T}
If $T$ is a Hopf monad on an autonomous category $\cc$, then the adjunction morphisms $e \co T\rexcla{T} \to 1_\cc$ and
$h \co 1_\cc \to \rexcla{T}\,T$ are given by $e_X = s^r_{\ldual{X}}$ and $h_X =\rdual{(s^l_X)}$, where $s^l$ and $s^r$
denote the left and right antipodes of $T$.
\end{rem}

Recall that a functor $G\co \dd \to \cc$ \emph{preserves colimits} if the image under $G$ of a colimit in $\dd$ is a
colimit in $\cc$. A functor $G\co \dd \to \cc$ \emph{creates colimits} if, for any functor $F\co I \to \dd$ such that
$GF\co I \to \cc$ has a colimit, this colimit lifts uniquely to a colimit of $F$. See \cite{ML1} for more precise
definitions.

Since the forgetful functor of a monad which preserves colimits creates colimits (by \cite[Proposition~4.3.2]{Bor2}), Proposition~\ref{prop-right-adj} admits the following corollary:
\begin{cor}\label{thm-creat-colim}
If $T$ is a Hopf monad on an autonomous category $\cc$, then $T$ preserves colimits and the forgetful functor $U_T\co T \ti \cc \to \cc$ creates
and preserves colimits.
\end{cor}

\subsection{Coends and Hopf monads}\label{sect-coendHM}\label{sect-rappelcolimit}\label{sect-rappelcoend}\label{sect-coendHopfmon}
Let $\cc$ and $\dd$ be categories and $F\co \cc^\opp \times \cc \to \dd$ be a functor. A \emph{dinatural
transformation} $d \co F \dinato Z$ from $F$ to an object $Z$ of~$\dd$ is family $d=\{d_X \co F(X,X) \to
Z\}_{X \in \Ob(\cc)}$ of morphisms in $\dd$ satisfying the dinaturality condition:
\begin{equation*}
d_{Y}F(\id_Y,f)=F(f,\id_X)d_X
\end{equation*}
for every morphism $f\co X \to Y$ in $\cc$. We denote $\Dinat(F,Z)$ the set of dinatural transformations from $F$ to
$Z$.

A \emph{coend} of a functor $F\co \cc^\opp \times \cc \to \dd$ consists of an object $C$ of $\cc$ and a dinatural
transformation $i\co F \dinato C$ which is universal in the sense that, for every dinatural transformation $d\co F
\dinato Z$, there exists a unique morphism $r\co C \to Z$ such that $d_X=r \circ i_X$ for all $X \in \Ob(\cc)$. In
other words, the map:
\begin{equation*} \left\{\begin{array}{rcl}
\Hom_\dd(C,Z) &\to &\Dinat(F,Z)\\
r &\mapsto & r i
\end{array}\right.
\end{equation*}
is a bijection. The dinatural transformation $i$ is then called a \emph{universal dinatural transformation} for $F$. A
coend of $F$, if it exists, is unique up to unique isomorphism. Following \cite{ML1}, we  denote it $\int^{X \in \cc}
F(X,X)$.

Coends are well-behaved under adjunction:
\begin{lem}\label{lem-coend-adj}
Let $\cc$, $\dd$, $\ee$ be categories, $(F\co \cc \to \dd, U\co \dd \to \cc)$ be an adjunction,  and $G\co \dd^\opp
\times \cc \to \ee$ be a functor. We have:
\begin{equation*}
\int^{X \in \cc}\!\!\!\! \!\! \!\! G(F(X),X)\,\simeq\, \int^{Y \in \dd}\!\!\!\! \!\! \!\! G(Y,U(Y)),
\end{equation*}
meaning that if either coend exists, then both exist and they are naturally isomorphic.
\end{lem}
\begin{proof}
Denote $\eta\co 1_\cc \to UF$ and $\varepsilon\co FU \to 1_\dd$ the adjunction morphisms. The lemma results from the
existence of a bijection:
\begin{equation*}
\psi\co \Dinat\bigl(G(F\times 1_\cc), E\bigr) \to \Dinat\bigl(G(1_{\dd^\opp} \times U), E\bigr)
\end{equation*}
which is natural in $E\in \Ob(\ee)$. It is defined by $\psi(d)=d_{U}G(\varepsilon, \id_{U})$, and its inverse by $\psi^{-1}(t)=t_{F}G(\id_{F}, \eta)$.
\end{proof}

Coends are special cases of colimits (see \cite{ML1}), and particular, a functor which preserves (resp.\@ creates)
colimits preserves (resp.\@ creates) coends. Hence, by Corollary~\ref{thm-creat-colim}:

\begin{prop}\label{HM-creat-colim}
Let $T$ be a Hopf monad on an autonomous category $\cc$ and $F\co \dd^\opp \times \dd \to T\ti\cc$ be a functor. Then
the coend $C=\int^{Y \in \dd} U_TF(Y,Y)$ exists if and only if the coend $\int^{Y \in \dd} F(Y,Y)$ exists. Moreover,
given a coend $C=\int^{Y \in \dd} U_TF(Y,Y)$ with universal dinatural transformation $i_Y\co U_TF(Y,Y) \to C$,  there
exists a unique action $r\co T(C) \to C$ of $T$ on $C$ such that $i_Y\co F(Y,Y) \to (C,r)$ is $T$\ti linear. We have
then  $(C,r)=\int^{Y \in \dd} F(Y,Y)$, with universal dinatural transformation $i$. The morphism $r\co T(C) \to C$ is
characterized by
\begin{equation*}
rT(i_Y)= i_Y \alpha_Y \quad \text{where} \quad F(Y,Y)=(U_TF(Y,Y), \alpha_Y),
\end{equation*}
as $T(i)$ is a universal dinatural transformation.
\end{prop}

\subsection{Pushforward of a monad under an
adjunction}\label{sect-dir-img} Let $(F\co \cc \to \dd, U\co
\dd \to \cc)$ be an adjunction and $Q$ be an endofunctor of $\dd$.
The endofunctor $UQF$ of $\cc$ is called  the \emph{pushforward of
$Q$ under the adjunction $(F,U)$} and is denoted by $(F,U)_*Q$.

If $Q$ is a monad, then $(F,U)_*Q$ is a monad: it is the monad of the composite adjunction $(F_QF,UU_Q)$ of $(F,U)$ and
$(F_Q,U_Q)$.

If $Q$ is comonoidal and $(F,U)$ is monoidal, then $(F,U)_*Q$ is comonoidal with comonoidal structure the composition
of the comonoidal structures of $U_T$, $Q$, and~$F_T$.

By Theorem~\ref{thm-gen-ex}, if the adjunction $(F,U)$ is monoidal and $Q$ is a bimonad, then $(F,U)_*Q$ is a bimonad
(since the composite of monoidal adjunctions is a monoidal adjunction).

Finally, if $\cc$ and $\dd$ are autonomous, $(F,U)$ is monoidal, and $Q$ is a Hopf monad, then $(F,U)_*Q$ is
a Hopf monad.

\begin{rem}
The structural morphisms of  $(F,U)_*Q$ can be expressed using those of $Q$ and the adjunction morphisms of $(F,U)$ (by applying Remark~\ref{rem-explicit-ex-gen}).
\end{rem}

\subsection{Cross products}\label{sect-crossprod}  Let $T$ be a monad on a category $\cc$ and $Q$ be an
endofunctor of $T\ti\cc$. Denote $\eta$ and $\varepsilon$ the unit and counit of $(U_T,F_T)$. The pushforward of $Q$ under the adjunction $(F_T,U_T)$ is called the \emph{cross product of
$Q$ by $T$} and denoted by $Q\cp  T$. Recall: $Q\cp  T=U_TQF_T$ as an endofunctor of $\cc$.

If $(Q,q,v)$ is a monad on~$T\ti\cc$, then $Q\cp  T$ is a monad on $\cc$ with product~$p$ and unit~$e$ given by:
\begin{equation*}
p=q_{F_T} Q(\varepsilon_{Q F_T})  \quad \text{and}\quad e=v_{F_T}\eta.
\end{equation*}

If $T$ is a bimonad and $Q$ is comonoidal, then $Q\cp  T$ is a comonoidal with comonoidal structure given by:
\begin{align*}
& (Q\cp  T)_2(X,Y)=Q_2\bigl(F_T(X),F_T(Y)\bigr)\, Q\bigl(\varepsilon_{F_T(X)\otimes F_T(Y)} F_T(\eta_X\otimes \eta_Y)\bigr),\\
& (Q\cp  T)_0=Q_0\, Q(\varepsilon_{(\un,T_0)}).
\end{align*}

If $T$ and $Q$ are bimonads, then $Q\cp  T$ is a
bimonad. If $T$ and $Q$ are Hopf monads, then $Q\cp  T$ is a Hopf monad, with left and right antipodes given by:
\begin{equation*}
a^l_X=\ldual{\eta}_X  S^l_{\ldual{F}_T(X)}  Q(\varepsilon_{\ldual{Q}F_T(X)}) \quad \text{and} \quad
a^r_X=\eta^\vee_X  S^r_{F^\vee_T(X)}  Q(\varepsilon_{QF_T(X)^\vee}),
\end{equation*}
where $S^l$ and $S^r$ are the antipodes of $Q$.

\begin{exa}
Let $H$ be a bialgebra over a field $\kk$ and $A$ be a $H$-module algebra, that is, an algebra in the monoidal category $\lmod{\mathrm{Vect}_\kk}{H}$ of left $H$\ti modules. In this situation, we may form the cross product $A \rtimes H$, which is a $\kk$-algebra (see \cite{Maj2}). Recall $H \otimes ?$ is a monad on $\Vect_\kk$ and $A \otimes ?$ is a monad on $\lmod{\mathrm{Vect}_\kk}{H}$. Then:
\begin{equation*}
(A \otimes ?)\cp (H \otimes ?)=(A\rtimes H) \otimes ?
\end{equation*}
as monads. Moreover, if $H$ is a quasitriangular bialgebra and $A$ is a $H$-module bialgebra, that is, a bialgebra in the braided category $\lmod{\mathrm{Vect}_\kk}{H}$, then $A\rtimes H$ is a $\kk$-bialgebra, and
$(A \otimes ?)\cp (H \otimes ?)=(A\rtimes H) \otimes ?$ as bimonads.
\end{exa}

\section{Distributive laws and liftings}\label{sect-distlaw}

Given two monads $P$ and $T$ on a category $\cc$,  when is the composition $PT$ a monad? How can one lift $P$ to a
monad on the category $T\ti\cc$? Beck's theory of distributive law \cite{Beck1} provides an answer for these questions. In this section, we recall the basic results of this theory and extend them to Hopf monads.

\subsection{Distributive laws between algebras} \label{sect-cartier}
Let $(A,m,u)$ and $(B,\mu,\eta)$ be two algebras in a monoidal
category $\cc$. Given a morphism $\Omega \co B \otimes A \to A \otimes B$ in $\cc$, set:
\begin{equation*}
p=(m \otimes \mu)(\id_A \otimes \Omega \otimes \id_B) \co (A \otimes B) \otimes (A \otimes B) \to (A \otimes B).
\end{equation*}
Then $(A \otimes B, p, u\otimes \eta)$ is an algebra in $\cc$ if and only if $\Omega$ satisfies:
\begin{align*}
& \Omega (\id_B \otimes m)=(m \otimes \id_B)(\id_A \otimes \Omega)(\Omega \otimes \id_A) ; &
& \Omega(\id_B \otimes u)=u \otimes \id_B ;\\
& \Omega (\mu \otimes \id_A)=(\id_A \otimes \mu)(\Omega \otimes \id_B)(\id_B \otimes \Omega) ; &
& \Omega(\eta \otimes \id_A)=\id_A \otimes \eta.
\end{align*}
If such is the case, we say that $\Omega$ is a \emph{distributive law of $B$ over $A$}. The algebra $(A \otimes B, p,
u\otimes \eta)$ is then denoted $A \otimes_\Omega
B$. Note that $i=(\id_A \otimes \eta)\co A \to A \otimes_\Omega B$ and $j=(u \otimes \id_B)\co B \to A\otimes_\Omega B$
are algebra morphisms, and the \emph{middle unitary law} holds:
\begin{equation*}
p (\id_A \otimes \eta \otimes u \otimes \id_A)= \id_{A \otimes B}.
\end{equation*}
In other words, we have $p(i \otimes j)=\id_{A \otimes B}$.

\begin{rem} \label{rem-cartier}
Let $(C,p,e)$ be an algebra in $\cc$ and $i \co A \to C$, $j\co B \to C$ be two algebra morphisms such that $\Theta=p(i
\otimes j)\co A \otimes B \to C$ is an isomorphism in $\cc$. Then there exists a unique distributive law
$\Omega$ of $B$ over $A$ such that $\Theta$ is an algebra isomorphism from $A
\otimes_\Omega B$ to $C$. Moreover:
\begin{equation*}
\Omega= \Theta^{-1}p(j \otimes i),\quad i=\Theta(\id_A \otimes \eta)\quad\mbox{and}\quad j=\Theta(u \otimes \id_B).
\end{equation*}
\end{rem}

\begin{rem}\label{rem-cartieriso}
If a distributive law $\Omega\co B \otimes A \to A \otimes B$  of $B$ over $A$ is an
isomorphism, then $\Omega^{-1}$ is a distributive law of $A$ over $B$ and $\Omega\co B \otimes_{\Omega^{-1}} A \to  A
\otimes_\Omega B $ is an isomorphism of algebras.
\end{rem}

\begin{exa}\label{exa-cartierHA}
Let $A$ and $B$ be bialgebras in a braided
category~$\bb$. A
distributive law of $B$ over $A$ is \emph{comultiplicative} if it satisfies:
\begin{align*}
& (\id_A \otimes \tau_{A,B} \otimes \id_B)(\Delta_A \otimes
\Delta_B) \Omega
=(\Omega \otimes \Omega) (\id_B \otimes \tau_{B,A} \otimes \id_A)(\Delta_B \otimes \Delta_A),\\
& (\varepsilon_A \otimes \varepsilon_B) \Omega=\varepsilon_B \otimes
\varepsilon_A,
\end{align*}
where $\tau$ is the braiding of $\bb$. A comultiplicative distributive law is nothing but a distributive law
between algebras in the monoidal category of
coalgebras in $\bb$. Let $\Omega$ be a comultiplicative distributive law
of $B$ over $A$. Then  $A\otimes_\Omega B$ is a bialgebra in $\bb$. Furthermore, if $A$ and $B$ are Hopf algebras, then $A\otimes_\Omega B$ is a Hopf algebra with structural morphisms:
\begin{align*}
&m_{A\otimes_\Omega B}=(m_A \otimes m_B)(\id_A \otimes \Omega \otimes \id_B), &&u_{A\otimes_\Omega B}=u_A \otimes u_B,\\
&\Delta_{A\otimes_\Omega B}=(\id_A \otimes \tau_{A,B} \otimes
\id_B)(\Delta_A \otimes \Delta_B), &
&\varepsilon_{A\otimes_\Omega B}= \varepsilon_A \otimes \varepsilon_B,\\
&S_{A\otimes_\Omega B}=S_A \otimes S_B,
\end{align*}
where $m_C$, $u_C$, $\Delta_C$, $\varepsilon_C$, $S_C$ denote
respectively the product, unit, coproduct, counit, and antipode of a
Hopf algebra $C$.
\end{exa}

\subsection{Lifting monads and bimonads} \label{sect-def-lift}
Let $(P,m,u)$ be a monad on a category $\cc$ and $U\co\dd \to \cc$ be a functor. A \emph{lift of the monad $P$ to
$\dd$} is a monad $(\Tilde{P},\Tilde{m},\Tilde{u})$ on~$\dd$ such that $PU=U\Tilde{P}$, $m_{U}=U(\Tilde{m})$, and
$u_{U}=U(\Tilde{u})$.

Let $P$ be a bimonad on a monoidal category $\cc$ and $U \co \dd \to \cc$ be a strong monoidal functor. A \emph{lift of
the bimonad $P$ to $\dd$} is bimonad $\Tilde{P}$ on $\dd$ which is a lift of the monad $P$ to $\dd$ such that $U
\Tilde{P}=P U$ as comonoidal functors.

\subsection{Distributive laws between monads}\label{sect-dist-law-monads}
Let $(T,\mu,\eta)$ and $(P,m,u)$ be monads on a category $\cc$. Following Beck \cite{Beck1}, a \emph{distributive law
of $T$ over $P$} is a natural transformation $\Omega\co TP \to PT$ verifying:
\begin{align*}
& \Omega_X T(m_X)=m_{T(X)} P(\Omega_X) \Omega_{P(X)} ; &
& \Omega_X T(u_X)=u_{T(X)} ; \\
& \Omega_X \mu_{P(X)}=P(\mu_X) \Omega_{T(X)} T(\Omega_X) ; &
& \Omega_X \eta_{P(X)}=P(\eta_X) ;
\end{align*}
for all object $X$ of $\cc$.

\begin{rem}
Viewing the monads $T$ and $P$ as algebras in the monoidal category of endofunctors of $\cc$ (with monoidal product
$\circ$ and unit object $1_\cc$), the above definition of a distributive law agrees with that given in
Section~\ref{sect-cartier}.
\end{rem}

Let $\Omega$ be a distributive law of $T$ over $P$. Firstly, $\Omega$ defines a monad structure on the endofunctor
$PT$ of $\cc$, with product $p$ and unit $e$ given by:
\begin{equation*}
p_X=m_{T(X)}P^2(\mu_X)P(\Omega_{T(X)}) \quad \text{and} \quad e_X=u_{T(X)}\eta_X.
\end{equation*}
The monad $(PT, p,e)$ is denoted $P \circ_\Omega T$. Secondly $\Omega$ defines a lift $(\Tilde{P}^\Omega, \Tilde{m},
\Tilde{u})$ of the monad $P$ to the category $T\ti\cc$ as follows:
\begin{equation*}
\Tilde{P}^\Omega(M,r)=\bigl(P(M), P(r)\Omega_M\bigr), \quad \Tilde{m}_{(M,r)}=m_M, \quad \text{and} \quad \Tilde{u}_{(M,r)}=u_M.
\end{equation*}
Furthermore,  there is a canonical isomorphism of categories:
\begin{equation*}
K\co\left\{\begin{array}{ccc}
\Tilde{P}^\Omega\ti(T\ti\cc) & \longrightarrow & (P\circ_\Omega T)\ti\cc \\
\bigl ((M,r),s \bigr) & \longmapsto & \bigl(M,U_T(s)P(r)\bigr)
\end{array}\right.
\end{equation*}
with inverse: \begin{equation*}
K^{-1}\co\left\{\begin{array}{ccl}
(P\circ_\Omega T)\ti\cc & \longrightarrow & \Tilde{P}^\Omega\ti(T\ti\cc) \\
(A,\alpha)& \longmapsto & \bigl ((A,\alpha \,u_{T(A)}), \alpha P(\eta_A) \bigr)
\end{array}\right..
\end{equation*}
In fact $K$ is the comparison functor of the composite adjunction:
\begin{equation*}
\Tilde{P}^\Omega\ti(T\ti\cc) \adjunct{U_{\Tilde{P}^\Omega}}{F_{\Tilde{P}^\Omega}}T\ti\cc\adjunct{U_T}{F_T} \cc
\end{equation*}
Hence this composite adjunction is monadic with monad $P\circ_\Omega T$.

The assignments $\Omega \mapsto P \circ_\Omega T$ and $\Omega \mapsto \Tilde{P}^\Omega$ are one-to-one in the following
sense:

\begin{thm}[\cite{Beck1}] \label{propdistrib}
Let $(T,\mu,\eta)$ and $(P,m,u)$ be monads on a category $\cc$.  We have bijective correspondences between:
\begin{enumerate}
\renewcommand{\labelenumi}{{\rm (\roman{enumi})}}
\item Distributive laws $\Omega\co TP \to PT$ of $T$ over $P$;
\item Products $p\co PTPT \to PT$ for which:
\begin{enumerate}
\renewcommand{\labelenumi}{{\rm (\alpha{enumi})}}
\item $(PT,p,u_T\eta)$ is a monad on $\cc$;
\item $u_T \co T \to PT$ and $P(\eta) \co P \to PT$ are morphisms of monads;
\item the middle unitary law $p_XP(\eta_{PT(X)} u_{T(X)})=\id_{PT(X)}$ holds;
\end{enumerate}
\item Lifts of the monad $P$ on $\cc$ to a monad $\Tilde{P}$ on $T\ti\cc$.
\end{enumerate}
\end{thm}

\subsection{Distributive laws between bimonads}\label{sect-dist-law-bimonads}
Let $T$ and $P$ be bimonads on a monoidal category $\cc$. Recall that $TP$ and $PT$ are comonoidal endofunctors of
$\cc$. A distributive law $\Omega\co TP \to PT$ of $T$ over $P$ is \emph{comonoidal} if it is comonoidal as a natural
transformation, that is, if it satisfies:
\begin{equation*}
(PT)_2(X,Y) \Omega_{X \otimes Y}=(\Omega_X \otimes \Omega_Y) (TP)_2(X,Y)
\quad \text{and} \quad (PT)_0 \Omega_\un=(TP)_0.
\end{equation*}

\begin{rem}
Viewing the bimonads $T$ and $P$ as algebras in the monoidal category of comonoidal endofunctors of $\cc$ (see
Remark~\ref{bimonad-alg}), a comonoidal distributive law is a distributive law in the sense of
Section~\ref{sect-cartier}.
\end{rem}

Beck's Theorem \ref{propdistrib} was generalized by Street~\cite{Street0} to monads in a 2\ti category. Applying this theorem to the case of the 2-category of monoidal categories and comonoidal functors, we obtain:

\begin{thm}\label{propdistribbimon}
Let $(T,\mu,\eta)$ and $(P,m,u)$ be bimonads on a monoidal category~$\cc$. We have bijective correspondences between:
\begin{enumerate}
\renewcommand{\labelenumi}{{\rm (\roman{enumi})}}
\item Comonoidal distributive laws $\Omega\co TP \to PT$ of $T$ over $P$;
\item Products $p\co PTPT \to PT$ for which:
\begin{enumerate}
\item $(PT,p,u_T\eta)$ is a bimonad on $\cc$;
\item $u_T \co T \to PT$ and $P(\eta) \co P \to PT$ are morphisms of bimonads;
\item the middle unitary law $p_XP(\eta_{PT(X)} u_{T(X)})=\id_{PT(X)}$ holds.
\end{enumerate}
\item Lifts of the bimonad $P$ on $\cc$ to a bimonad $\Tilde{P}$ on $T\ti\cc$.
\end{enumerate}
Also, if $\Omega$ is a comonoidal distributive law of $T$ over $P$, the canonical isomorphism of categories
$\Tilde{P}^\Omega\ti(T\ti\cc)\simeq (P\circ_\Omega T)\ti\cc$ is strict monoidal.
\end{thm}

\begin{exa}\label{exa-monA}
Let $\bb$ be a braided category, $A$ and $B$ be two bialgebras in $\bb$, and $\Omega\co B \otimes A \to A \otimes B$ be a morphism in $\bb$. Then the following conditions are equivalent:
\begin{enumerate}
\renewcommand{\labelenumi}{{\rm (\roman{enumi})}}
\item $\Omega \otimes \id_{1_\bb}$ is a comonoidal distributive law of $B \otimes ?$ over $A \otimes ?$;
\item $\id_{1_\bb} \otimes \Omega$ is a comonoidal law of $? \otimes A$ over $? \otimes B$;
\item $\Omega$ is a comultiplicative distributive law of $B$ over $A$ (see Example~\ref{exa-cartierHA}).
\end{enumerate}
If such is the case,
we have the following equalities of bimonads:
\begin{align*}
 (A \otimes ?)\circ_{(\Omega\otimes \id_{1_\bb})} (B\otimes ?)&=(A\otimes_\Omega B)\otimes ?\\
(? \otimes B)\circ_{(\id_{1_\bb} \otimes \Omega)} (?\otimes A)&=\,? \otimes (A\otimes_\Omega B).
\end{align*}
\end{exa}

\begin{rem}\label{rem-crossprod-distlaw}
Let $\Omega \co TP \to PT$ be a distributive law between monads on a category $\cc$.
Then
$\Tilde{P}^\Omega \cp  T=P\circ_\Omega T$ as monads, where $\cp $ denotes the cross product (see
Section~\ref{sect-crossprod}). Moreover, if $\cc$ is monoidal, $T$ and $P$ are bimonads, and $\Omega$ is comonoidal,
then $\Tilde{P}^\Omega \cp  T=P\circ_\Omega T$ as bimonads.
\end{rem}

\subsection{Distributive laws and antipodes}
We show here that if $\Omega\co TP \to PT$ is a comonoidal distributive law between Hopf monads, then the composition
$P\circ_\Omega T$ and the lift $\Tilde{P}^\Omega$ are Hopf monads:

\begin{prop}\label{propdistribantip}
Let $T$ and $P$ be bimonads on a monoidal category~$\cc$ and let $\Omega\co TP \to PT$ be a comonoidal distributive law
of $T$ over $P$. Then:
\begin{enumerate}
\renewcommand{\labelenumi}{{\rm (\alph{enumi})}}
\item If $\cc$ is left autonomous, $T$ has a left antipode $s^l$,
and $P$ has a left antipode~$S^l$, then the bimonads $P\circ_\Omega T$ and $\Tilde{P}^\Omega$ have left antipodes,
denoted $a^l$ and $\Tilde{S}^l$ respectively, given by:
\begin{align*}
&a^l_X= S^l_X P(s^l_{P(X)}) PT(\ldual{\Omega_X})\co PT\bigl(\ldual{PT(X)}\bigr)
\to \ldual{X},\\
&\Tilde{S}^l_{(M,r)}=S^l_M \co \Tilde{P}^\Omega\bigl(\ldual{\Tilde{P}^\Omega(M,r)}\bigr) \to \ldual{(}M,r).
\end{align*}
\item If $\cc$ is right autonomous, $T$ has a right antipode $s^r$,
and $P$ has a right antipode~$S^r$, then the bimonads $P\circ_\Omega T$ and $\Tilde{P}^\Omega$ have right antipodes,
denoted $a^r$ and $\Tilde{S}^r$ respectively, given by:
\begin{align*}
&a^r_X= S^r_X P(s^r_{P(X)}) PT(\rdual{\Omega}_X)\co PT\bigl(\rdual{PT(X)}\bigr)
\to \rdual{X},\\
&\Tilde{S}^r_{(M,r)}=S^r_M \co \Tilde{P}^\Omega\bigl(\rdual{\Tilde{P}^\Omega(M,r)}\bigr) \to (M,r\rdual{)}.
\end{align*}
\end{enumerate}
\end{prop}

\begin{proof}
Let us prove Part (a). One first checks that $a^l_X$ satisfies the axioms of a left antipode, that is,
\begin{align*}
&(PT)_0 PT(\lev_X) PT(\ldual{(u_{T(X)}\eta_X)}\otimes \id_X)\\
&\qquad\qquad =\lev_{PT(X)} (a^l_{PT(X)} PT(\ldual{p_X}) \otimes \id_{PT(X)}) (PT)_2(\ldual{P}T(X),X);\\
&(u_{T(X)}\eta_X  \otimes \id_{\ldual{X}})\lcoev_X (PT)_0\\
&\qquad \qquad=(p_X \otimes a^l_X) (PT)_2(PT(X),\ldual{P}T(X))PT(\lcoev_{PT(X)}).
\end{align*}
This can be done applying the axioms for the left antipodes $s^l$ and $S^l$ of $T$ and $P$ and the axioms for the
distributive law $\Omega$. By  Theorem~\ref{thm-biHopfmon}(b), this implies that $(P\circ_\Omega T)\ti\cc$ is left
autonomous. Now recall that:
\begin{equation*} K\co\left\{\begin{array}{ccc}
\Tilde{P}^\Omega\ti(T\ti\cc) & \iso & (P\circ_\Omega T)\ti\cc \\
\bigl ((M,r),s \bigr) & \longmapsto & \bigl(M,sP(r)\bigr)
\end{array}\right.
\end{equation*}
is a strict monoidal isomorphism of categories (see Section~\ref{sect-dist-law-monads}). Therefore
$\Tilde{P}^\Omega\ti(T\ti\cc)$ is left autonomous and so, by Theorem~\ref{thm-biHopfmon}(b), $\Tilde{P}^\Omega$ has a
left antipode $\Tilde{S}^l$. Furthermore, given a  $\Tilde{P}^\Omega$\ti module $((M,r),s)$, we have:
\begin{equation*}
K^{-1}\bigl(\ldual{K}((M,r),s) \bigr )=\bigl ( (\ldual{M},s^l_MT(r)),\rho \bigr ).
\end{equation*}
where $U_T(\rho)=a^l_M PT(\ldual{P}(r)\ldual{s}) P(\eta_M) =S^l_M P(\ldual{s})$. Hence $\Tilde{S}^l_{(M,r)}=S^l_M$.

Part (b) results from Part (a) applied to the opposite Hopf monads
$$(P\circ_\Omega T)^\mop=P^\mop \circ_\Omega T^\mop \quad \text{and}
\quad (\Tilde{P}^\Omega)^\mop=(\widetilde{P^\mop})^\Omega,$$ see Remark~\ref{remoppositemonad}.
\end{proof}

From Proposition~\ref{propdistribantip}, Theorem~\ref{propdistribbimon}, and Remark~\ref{rem-crossprod-distlaw}, we deduce:
\begin{cor}\label{cordistribHopf}
If $T$ and $P$ are Hopf monads on an autonomous category $\cc$ and $\Omega\co TP \to PT$ is a comonoidal distributive
law, then $P\circ_\Omega T$ is a Hopf monad on $\cc$, $\Tilde{P}^\Omega$ is a Hopf monad on $T\ti\cc$, and $\Tilde{P}^\Omega \cp  T=P\circ_\Omega T$ as Hopf monads.
\end{cor}

\subsection{Invertible distributive laws}

Let $T$, $P$ be two monads on a category $\cc$ and $\Omega \co TP \to PT$ be an invertible distributive law of $T$ over
$P$. Then $\Omega^{-1}\co PT \to TP$ is a distributive law of $P$ over $T$, and $\Omega$ is a isomorphism of monads
from $T \circ_{\Omega^{-1}} P$ to $P \circ_\Omega T$.

If $\cc$ is monoidal, $P$, $T$ are bimonads, and $\Omega \co TP \to PT$ is a comonoidal distributive law of $T$ on $P$,
then $\Omega^{-1}\co PT \to TP$ is a comonoidal distributive law of $P$ over $T$, and $\Omega$ is a isomorphism of
bimonads from $T \circ_{\Omega^{-1}} P$ to $P \circ_\Omega T$.

\begin{prop}\label{propdistinv}
Let $P$, $T$ be Hopf monads on an autonomous category $\cc$. Then any comonoidal distributive law $\Omega \co TP \to
PT$ of $T$ over $P$ is invertible. Furthermore, for any object $X$ of $\cc$, we have:
\begin{equation*}
\Omega^{-1}_X=S^r_{\ldual{T}P(X)}P\bigl(s^r_{P(\ldual{T}P(X))}\bigr)PT\bigl(\Omega^\vee_{\ldual{T}P(X)}\bigr)
PT\bigl( P(s^l_{P(X)})^\vee\bigr)PT\bigl(S^{l\vee}_X\bigr),
\end{equation*}
where $s^l, s^r, S^l, S^r$ denote left and right antipodes of $T$ and $P$ respectively.
\end{prop}

\begin{proof}
The functors $T$, $P$ and $PT$ are Hopf monads by assumption and Corollary~\ref{cordistribHopf}. Therefore, by
Proposition~\ref{prop-right-adj}, the functors $\rexcla{T}$, $\rexcla{P}$, and $\rexcla{(PT)}$ are right adjoints for
$T$, $P$, and $PT$ respectively. On the other hand, by composition of adjunctions,
$\rexcla{P}\circ\rexcla{T}=\rexcla{(PT)}$ is a right adjoint for $TP$. As a left adjoint is unique up to unique natural
isomorphism, we obtain a canonical isomorphism $\alpha\co PT \to TP$. Denoting $e \co T \rexcla{T} \to 1_\cc$, $h \co 1_\cc
\to \rexcla{T} T$, $e' \co P \rexcla{P} \to 1_\cc$, $h' \co 1_\cc \to \rexcla{P} P$, $E \co PT \rexcla{(PT)} \to
1_\cc$, and $H \co 1_\cc \to\rexcla{(PT)} PT$ the adjunction morphisms, we have $\alpha=E_{TP}PT\rexcla{P}(h_P)PT(h')$.
Now the adjunction morphisms can be expressed in terms of the antipodes, see Remark~\ref{radj-T}. Therefore, using
Proposition~\ref{propdistribantip}, we get that, for any object $X$ of $\cc$,
\begin{equation*}
\alpha_X=S^r_{\ldual{TP(X)}}P(s^r_{P(\ldual{TP(X)})})PT(\Omega^\vee_{\ldual{TP(X)}})PT(
P(s^l_{P(X)})^\vee)PT(S^{l\vee}_X)
\end{equation*}
Furthermore: $E_X=S^r_{\ldual{X}}P(s^r_{P(\ldual{X})})PT(\Omega^\vee_{\ldual{X}})$ and $ H_X=\rexcla{(PT)}(\Omega_X)
P(s^l_{P(X)})^\vee S^{l\vee}_X$. Hence:
\begin{align*}
\id_{PT(X)}&=E_{PT(X)}PT(H_X) \\
&=E_{PT(X)}PT\rexcla{(PT)}(\Omega_X)PT(P(s^l_{P(X)})^\vee)PT(S^{l\vee}_X)\\
&=\Omega_X E_{TP(X)}PT( P(s^l_{P(X)})^\vee)PT(S^{l\vee}_X) \quad \text{by functoriality of $E$}\\
&=\Omega_X S^r_{\ldual{TP(X)}}P(s^r_{P(\ldual{TP(X)})})PT(\Omega^\vee_{\ldual{TP(X)}})PT( P(s^l_{P(X)})^\vee)PT(S^{l\vee}_X)\\
&=\Omega_X \alpha_X.
\end{align*}
This shows that $\Omega$, as inverse of the isomorphism $\alpha$, is an isomorphism.
\end{proof}

\begin{rem}\label{rem-invdistriHAlg}
Let $\Omega\co B \otimes A \to A \otimes B$ be a distributive law between two Hopf algebras $A$ and $B$ in a braided autonomous category $\bb$ with braiding $\tau$. Then, applying Proposition~\ref{propdistinv} to the distributive law of Example~\ref{exa-monA}, we find that $\Omega$ is invertible, and its inverse is given by:
$$\Omega^{-1}=(S^{-1}_B \otimes S^{-1}_A) \tau^{-1}_{B,A} \,\Omega \, \tau_{A,B} (S_A \otimes S_B),$$
where $S_A$ and $S_B$ are the antipodes of $A$ and $B$.
\end{rem}

\section{The centralizer of a Hopf monad}\label{sect-centHM}
In this section, we introduce the notion of centralizer of a Hopf monad, and interpret its category of modules as the categorical center relative to the Hopf monad.

\subsection{Centralizers of endofunctors}
Let $\cc$ be a monoidal category and $T$ be an endofunctor of $\cc$.

A \emph{centralizer of $T$ at an object $X$ of $\cc$} is a pair $(Z,\delta)$, where $Z\in \Ob(\cc)$ and
\begin{equation*}
\delta=\{\delta_Y \co X \otimes Y \to T(Y) \otimes Z \}_{Y \in \Ob(\cc)}\co X \otimes 1_\cc \to T \otimes Z
\end{equation*}
is a natural transformation, verifying the following universal property: for every object $W$ of $\cc$ and every
natural transformation $\xi\co X \otimes 1_\cc \to T \otimes W$, there exists a unique morphism $r\co Z \to W$ in
$\cc$ such that $\xi=(\id_T \otimes r)\delta$. Note that a centralizer of $T$ at $X$, if it exists, is unique up to
unique isomorphism.

\begin{rem}\label{rem-centleftright} The notion of  centralizer is not invariant under left/right symmetry. We should properly call it `left-handed' centralizer. We can as well define a `right-handed' centralizer of $T$ at $X$ to be a pair $(Z',\delta')$, with
\begin{equation*}
\delta'=\{\delta'_Y \co Y \otimes X \to Z' \otimes T(Y)\}_{Y \in \Ob(\cc)}\co 1_\cc \otimes X \to Z' \otimes T
\end{equation*}
satisfying the relevant universal property. Note that this is equivalent to saying that
$(Z',\delta')$ is a  `left-handed' centralizer of $T$ at $X$ in the monoidal category $\cc^{\otimes\opp}$.
By left/right symmetry, all notions and results concerning `left-handed' centralizers can be adapted to the `right-handed' version.
\end{rem}

The endofunctor $T$ is said to be \emph{centralizable at an object $X$ of $\cc$} if it admits a centralizer at~$X$.

A \emph{centralizer} of $T$ is a pair $(Z_T,\partial)$, where $Z_T$ is an endofunctor of $\cc$ and
\begin{equation*}
\partial=\{\partial_{X,Y} \co X \otimes Y \to T(Y) \otimes Z_T(X) \}_{X,Y \in \Ob(\cc)} \co \otimes \to
(T \otimes Z_T)\sigma_{\cc,\cc}
\end{equation*}
is a natural transformation, such that $(Z_T(X),
\partial_{X,1_\cc})$ is a centralizer of $T$ at $X$ for every object $X$ of $\cc$.

The endofunctor $T$ is said to be \emph{centralizable} if it admits a centralizer. An endofunctor of $\cc$ is
centralizable if and only if it is centralizable at every object of~$\cc$. In that case, its centralizer is essentially
unique. More precisely:

\begin{lem}\label{lem-centralisable}
Let $T$ be an endofunctor of a monoidal category $\cc$. We have:
\begin{enumerate}
\renewcommand{\labelenumi}{{\rm (\alph{enumi})}}
\item Given a centralizer
$(Z_T(X),\partial_X)$ of $T$ at every object $X$ of $\cc$, the assignment $Z_T \co X \mapsto Z_T(X)$ admits a unique
structure of functor such that:
\begin{equation*}
\partial=\{\partial_{X,Y}=(\partial_X)_Y\co X
\otimes Y \to T(Y) \otimes Z_T(X)\}_{X,Y \in \Ob(\cc)}
\end{equation*}
is a natural transformation. The pair $(Z_T,\partial)$ is then a centralizer of $T$.
\item If $(Z,\partial)$ and $(Z',\partial')$ are centralizers
of $T$, then there exists a unique natural isomorphism $\alpha\co Z \to Z'$ such that $\partial'=(\id_T \otimes
\alpha)\partial$.
\end{enumerate}
\end{lem}
\begin{proof} For each morphism $f\co X  \to X'$ in $\cc$, by the universal property of centralizers, there exists a unique morphism
$Z_T(f) \co Z_T(X) \to Z_T(X')$ such that:
\begin{equation*}
(\id_T \otimes Z_T(f))\,\partial_{X,1_\cc}=\partial_{X',1_\cc}\,(f \otimes 1_\cc),
\end{equation*}
and this assignment defines the only structure of functor on $Z_T$ such that $\partial$ is a natural transformation.
\end{proof}

\subsection{Centralizers and coends} In this section, we give a characterization
of centralizable endofunctors in a left autonomous category in terms of coends.

\begin{prop}\label{propcentral1}
Let $\cc$ be a left autonomous category, $T$ be an endofunctor of $\cc$, and $X$ be an object of $\cc$. Then $T$ is
centralizable at $X$ if and only if the coend
\begin{equation*}
Z_T(X)=\int^{Y \in \cc} \ldual{T(Y)} \otimes X \otimes Y
\end{equation*}
exists. If such is the case, denoting $i$ the universal dinatural transformation of the coend and setting:
\begin{equation*}
(\partial_X)_Y=\bigl(\id_{T(Y)} \otimes i_Y\bigr)(\lcoev_{T(Y)} \otimes \id_{X \otimes Y}),
\end{equation*}
the pair $(Z_T(X),\partial_X)$ is a centralizer of $T$ at $X$.
\end{prop}

\begin{proof}
Let $F\co \cc^\opp \times \cc \to \cc$ be the functor defined by $F(Y,Z)= \ldual{T(Y)} \otimes X \otimes Z$ and
$F(f,g)= \ldual{T(f)} \otimes X \otimes g$. By duality, we have a bijection:
\begin{equation*}
\psi\co\Dinat(F,Z) \to \Nat(X\otimes 1_\cc, T \otimes Z)
\end{equation*}
which is natural in $Z\in \Ob(\cc)$. It is defined by:
\begin{equation*}
\psi(j)_Y=\bigl(\id_{T(Y)} \otimes j_Y\bigr)(\lcoev_{T(Y)} \otimes \id_{X \otimes Y})\co X \otimes Y \to T(Y) \otimes Z
\end{equation*}
and its inverse by:
\begin{equation*}
\psi^{-1}(\delta)=(\lev_{T(Y)} \otimes \id_{Z})(\id_{\ldual{T(Y)}} \otimes \delta_Y) \co \ldual{T(Y)} \otimes X
\otimes Y \to Z.
\end{equation*}
Therefore $T$ is centralizable at~$X$ if
and only if $F$ admits a coend and, if so, the centralizer of $T$ at~$X$ is canonically isomorphic to the coend of $F$.
\end{proof}

\subsection{Extended factorization property of the centralizer}\label{sect-ext-prop}
Let $T$ be a centralizable endofunctor of a monoidal category $\cc$ and  $(Z_T,\partial)$ be a centralizer of $T$.  For
any non-negative integer~$n$, let
\begin{equation*}
\partial^n\co \otimes_{n+1} \rightarrow (T^{\otimes n} \otimes Z_T^n) \sigma_{\cc, \cc^n}
\end{equation*}
be the natural transformation defined by the following diagram:
\begin{center}
\psfrag{Y}[Bc][Bc]{\scalebox{.7}{$Y$}}
\psfrag{A}[Bc][Bc]{\scalebox{.7}{$T(Y_1)$}}
\psfrag{B}[Bc][Bc]{\scalebox{.7}{$T(Y_2)$}}
\psfrag{F}[Bc][Bc]{\scalebox{.7}{$T(Y_n)$}}
\psfrag{C}[Bc][Bc]{\scalebox{.7}{$Y_1$}}
\psfrag{R}[Bc][Bc]{\scalebox{.7}{$Y_2$}}
\psfrag{G}[Bc][Bc]{\scalebox{.7}{$Y_n$}}
\psfrag{X}[Bc][Bc]{\scalebox{.7}{$X$}}
\psfrag{T}[Bc][Bc]{\scalebox{.7}{$T(Y)$}}
\psfrag{Z}[Bc][Bc]{\scalebox{.7}{$Z_T(X)$}}
\psfrag{K}[Bc][Bc]{\scalebox{.7}{$Z^n_T(X)$}}
$\partial^n_{X,Y_1, \dots, Y_n}=$ \rsdraw{.45}{.9}{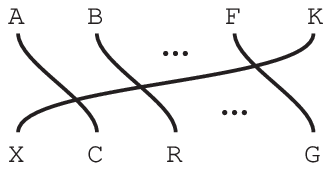} \quad where  \quad $\partial_{X,Y}=$
\rsdraw{.45}{.9}{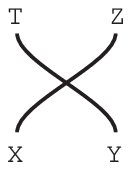}\;.
\end{center}
In other words, the morphisms:
\begin{equation*}
\partial^n_{X,Y_1, \dots, Y_n}\co X\otimes Y_1
\otimes \cdots \otimes Y_n \to T(Y_1) \otimes \cdots \otimes T(Y_n) \otimes Z_T^n(X)
\end{equation*}
are defined inductively by $\partial^0_X=\id_X$ and
\begin{equation*}
\partial^{n+1}_{X,Y_1, \dots, Y_{n+1}}=(\id_{T(Y_1) \otimes \cdots
\otimes T(Y_n)} \otimes
\partial_{Z_T^n(X),Y_{n+1}})(\partial^n_{X,Y_1, \dots, Y_n} \otimes
\id_{Y_{n+1}}).
\end{equation*}
Notice $\partial^1=\partial$ and $\partial^{p+q}=(\id_{T^{\otimes p}} \otimes
\partial^q)(\partial^p \otimes \id_{\otimes_q})$ for all non-negative integers~$p,q$.

\begin{lem}\label{lemcentralizer}
Assume $\cc$ is left autonomous. Let $\dd$ be category and $K,L\co \dd \to \cc$ be two functors. For every non-negative
integer $n$ and every natural transformation $\xi\co K\otimes \otimes_n \rightarrow (T^{\otimes n} \otimes L) \sigma_{\dd,
\cc^n}$, there exists a unique natural transformation $r\co Z_T^nK \rightarrow L$ such that:
\begin{equation*}
(\id_{T(Y_1) \otimes \cdots \otimes T(Y_n)} \otimes r_X)\partial^n_{X,Y_1, \dots, Y_n} =\xi_{X,Y_1, \dots, Y_n},
\end{equation*}
that is,
\begin{center}
\psfrag{A}[Bc][Bc]{\scalebox{.7}{$T(Y_1)$}} \psfrag{B}[Bc][Bc]{\scalebox{.7}{$T(Y_2)$}}
\psfrag{F}[Bc][Bc]{\scalebox{.7}{$T(Y_n)$}} \psfrag{C}[Bc][Bc]{\scalebox{.7}{$Y_1$}}
\psfrag{R}[Bc][Bc]{\scalebox{.7}{$Y_2$}} \psfrag{G}[Bc][Bc]{\scalebox{.7}{$Y_n$}}
\psfrag{X}[Bc][Bc]{\scalebox{.7}{$K(X)$}} \psfrag{L}[Bc][Bc]{\scalebox{.7}{$L(X)$}}
\psfrag{r}[Bc][Bc]{\scalebox{.9}{$r_X$}}
\psfrag{M}[Bc][Bc]{\scalebox{.9}{$\xi_{X,Y_1, \dots, Y_n}$}}
\rsdraw{.45}{.9}{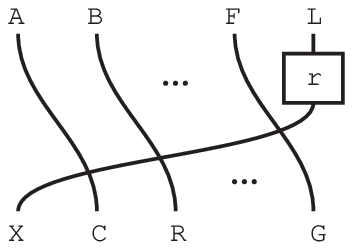} \; $=$ \; \rsdraw{.45}{.9}{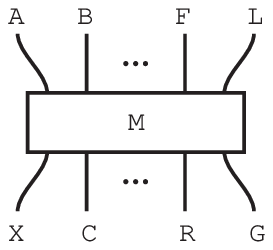}\,, \\[1em]
\end{center}
for all $X\in\Ob(\dd)$ and $Y_1,\dots,Y_n  \in \Ob(\cc)$.
\end{lem}
\begin{rem}
We will often write the equality defining $r$ in
Lemma~\ref{lemcentralizer} as:
\begin{equation*} (\id_{T^{\otimes n}}
\otimes r)\partial^n_{K, 1_{\cc^n}}=\xi.
\end{equation*}
Strictly speaking, it should be:
$(\id_{T^{\otimes n}} \otimes r)_{\sigma_{\dd,\cc^n}}\partial^n_{K,
1_{\cc^n}}=\xi$.
However, in this kind of formulae, we will usually  omit the
permutation $\sigma$  as it can easily be recovered from the
context.
\end{rem}
\begin{proof}[Proof of Lemma~\ref{lemcentralizer}]
The lemma can be verified by induction on $n$ using the Parameter Theorem and Fubini Theorem for coends (see
\cite{ML1}) and the fact that, by Proposition~\ref{propcentral1}, we have $Z_TK(X)=\int^{Y \in \cc} \ldual{T(Y)}
\otimes K(X) \otimes Y$ for all $X\in\Ob(\cc)$.
\end{proof}

\subsection{Structure of centralizers}
\begin{figure}[h]
\begin{center}
\psfrag{A}[Bc][Bc]{\scalebox{.7}{$T(Y_1)$}}
\psfrag{B}[Bc][Bc]{\scalebox{.7}{$T(Y_2)$}}
\psfrag{C}[Bc][Bc]{\scalebox{.7}{$Y_1$}}
\psfrag{R}[Bc][Bc]{\scalebox{.7}{$Y_2$}}
\psfrag{X}[Bc][Bc]{\scalebox{.7}{$X$}}
\psfrag{U}[Bc][Bc]{\scalebox{.7}{$Y_1 \otimes Y_2$}}
\psfrag{L}[Bc][Bc]{\scalebox{.7}{$Z_T(X)$}}
\psfrag{r}[Bc][Bc]{\scalebox{.9}{$m_X$}}
\psfrag{T}[Bc][Bc]{\scalebox{.8}{$T_2(Y_1,Y_2)$}}
\rsdraw{.45}{.9}{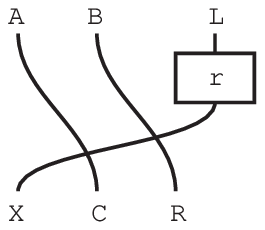} \, $=$ \, \rsdraw{.45}{.9}{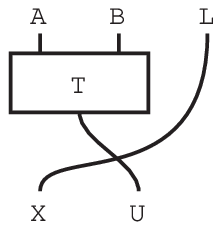},
\psfrag{X}[Bc][Bc]{\scalebox{.7}{$X$}}
\psfrag{L}[Bc][Bc]{\scalebox{.7}{$Z_T(X)$}}
\psfrag{r}[Bc][Bc]{\scalebox{.9}{$T_0$}}
\qquad \quad $u=$ \, \rsdraw{.45}{.9}{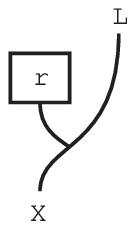},\\[1em]
\psfrag{A}[Bc][Bc]{\scalebox{.7}{$Z_T(X_1)$}}
\psfrag{B}[Bc][Bc]{\scalebox{.7}{$Z_T(X_2)$}}
\psfrag{C}[Bc][Bc]{\scalebox{.7}{$Y$}}
\psfrag{U}[Bc][Bc]{\scalebox{.7}{$T(Y)$}}
\psfrag{X}[Bc][Bc]{\scalebox{.7}{$X_1 \otimes X_2$}}
\psfrag{E}[Bc][Bc]{\scalebox{.7}{$X_1$}}
\psfrag{H}[Bc][Bc]{\scalebox{.7}{$X_2$}}
\psfrag{r}[Bc][Bc]{\scalebox{.9}{$\mu_X$}}
\psfrag{T}[Bc][Bc]{\scalebox{.8}{$(Z_T)_2(X_1,X_2)$}}
\rsdraw{.45}{.9}{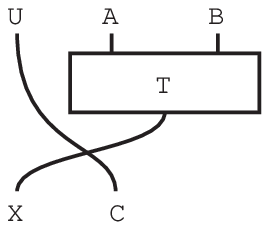} \, $=$ \, \rsdraw{.45}{.9}{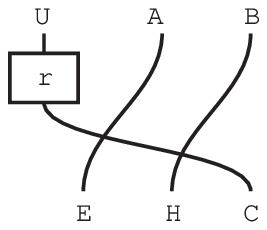},
\psfrag{U}[Bc][Bc]{\scalebox{.7}{$T(Y)$}}
\psfrag{X}[Bc][Bc]{\scalebox{.7}{$Y$}}
\psfrag{C}[Bc][Bc]{\scalebox{.7}{$Z_T(\un)$}}
\psfrag{r}[Bc][Bc]{\scalebox{.9}{$(Z_T)_0$}}
\psfrag{e}[Bc][Bc]{\scalebox{.9}{$\eta_X$}}
\qquad \quad \rsdraw{.45}{.9}{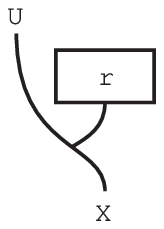} \, $=$ \, \rsdraw{.45}{.9}{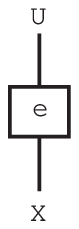}\;,\\[1em]
\psfrag{B}[Bc][Bc]{\scalebox{.7}{$T(Y)$}}
\psfrag{X}[Bc][Bc]{\scalebox{.7}{$\ldual{Z}_T(X)$}}
\psfrag{R}[Bc][Bc]{\scalebox{.7}{$Y$}}
\psfrag{L}[Bc][Bc]{\scalebox{.7}{$\ldual{X}$}}
\psfrag{r}[Bc][Bc]{\scalebox{.9}{$S^l_X$}}
\rsdraw{.45}{.9}{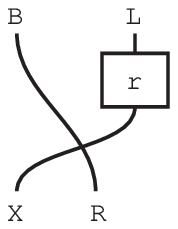}  $=$   \psfrag{r}[Bc][Bc]{\scalebox{.9}{$s^r_Y$}} \rsdraw{.45}{.9}{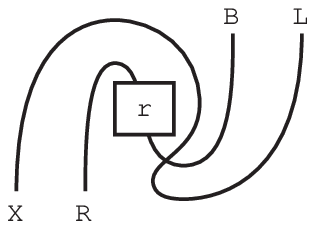}\;,
\qquad \psfrag{B}[Bc][Bc]{\scalebox{.7}{$T(Y)$}}
\psfrag{X}[B][Bc]{\scalebox{.7}{$Z_T(X)^\vee$}}
\psfrag{R}[Bc][Bc]{\scalebox{.7}{$Y$}}
\psfrag{L}[Bl][Bl]{\scalebox{.7}{$\rdual{X}$}}
\psfrag{r}[Bc][Bc]{\scalebox{.9}{$S^r_X$}}
\rsdraw{.45}{.9}{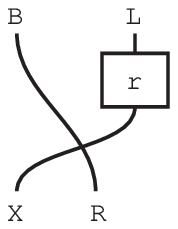}  $=$   \psfrag{r}[Bc][Bc]{\scalebox{.9}{$s^l_Y$}} \rsdraw{.45}{.9}{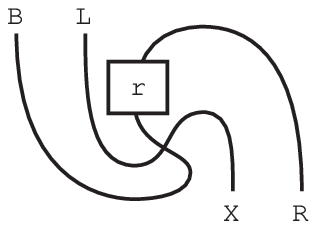}\;.
\end{center}
\caption{Structural morphisms of $Z_T$}
\label{fig-def-ZT}
\end{figure}

In this section, we show that the centralizer $Z_T$ of a Hopf monad $T$ is a Hopf monad. The structural morphisms of
$Z_T$ are defined as in Figure~\ref{fig-def-ZT} using the extended factorization property of $Z_T$ given in
Lemma~\ref{lemcentralizer}.
More precisely:

\begin{thm}\label{doub1}
Let $T$ be a centralizable endofunctor of a left autonomous category~$\cc$ and let $(Z_T,\partial)$ be its centralizer.
We have:
\begin{enumerate}
\renewcommand{\labelenumi}{{\rm (\alph{enumi})}}
\item If $T$ is comonoidal, then $Z_T$ is a monad on $\cc$, with product $m\co Z_T^2 \rightarrow Z_T$ and unit $u\co 1_\cc \rightarrow
Z_T$ defined by:
\begin{equation*}
(\id_{T^{\otimes 2}} \otimes m)\partial^2=(T_2 \otimes \id_{Z_T})\partial_{1_\cc,\otimes} \quad \text{and} \quad u=(T_0
\otimes \id_{Z_T})\partial_{1_\cc,\un}.
\end{equation*}
\item If $(T,\mu,\eta)$ is a monad, then $Z_T$ is comonoidal, with comonoidal structure
defined by:
\begin{align*}
&\bigl(\id_T \otimes (Z_T)_2\bigr)\partial_{\otimes , 1_\cc}=(\mu \otimes \id_{Z_T^{\otimes 2}})(\partial_{1_\cc,T}
\otimes
\id_{Z_T})(\id_{1_\cc} \otimes \partial);\\
& \bigl(\id_T \otimes (Z_T)_0\bigr)\partial_{\un,1_\cc}=\eta.
\end{align*}
\item If $T$ is a bimonad, then $Z_T$ is a bimonad on $\cc$, with the monad structure of Part (a) and the comonoidal structure of Part (b).
\item If $\cc$ is autonomous, $T$ is a bimonad, and $T$ has a right antipode $s^r$, then the bimonad $Z_T$ has a
left antipode $S^l$ defined by:
\begin{equation*}
(\id_T \otimes S^l)\partial_{\,\ldual{Z_T}, 1_\cc}=\ldual{\bigl((s^r \otimes \id_{Z_T})\partial_{1_\cc,
\rdual{T}}\bigr)}.
\end{equation*}
\item If $\cc$ is autonomous, $T$ is a bimonad, and $T$ has a left antipode $s^l$, then the bimonad $Z_T$
has a right antipode $S^r$ defined by:
\begin{equation*}
(\id_T \otimes S^r)\partial_{Z^\vee_T, 1_\cc}=\rdual{\bigl((s^l \otimes \id_{Z_T})\partial_{1_\cc, \ldual{T}}\bigr)}.
\end{equation*}
\end{enumerate}
In particular if $\cc$ is autonomous and $T$ is a Hopf monad, then $Z_T$ is a Hopf monad.
\end{thm}

\begin{rem}\label{rem-Zfunctorial}
The centralizer construction $T \mapsto Z_T$ is functorial, contravariant in $T$. More precisely, let $\cc$ be a left autonomous category and $T$, $T'$ be two centralizable endofunctors of $\cc$, with centralizers $(Z_T,\partial)$ and $(Z_{T'},\partial')$ respectively. Then, for each natural transformation $f\co T \to T'$, there exists a unique
natural transformation $Z_f\co  Z_{T'} \to Z_{T}$ such that:
$$(\id_{T'} \otimes Z_f)\partial'=(f \otimes \id_{Z_T})\partial.$$
We have: $Z_{fg}=Z_g Z_f$ and $Z_{\id_T}=\id_{Z_T}$. Moreover, if $f$ is comonoidal, then $Z_f$ is a morphism of monads. If $f$ if a morphism of monads, then $Z_f$ is comonoidal. Thus, if $f$ is a morphism of bimonads or Hopf monads, so is $Z_f$.
\end{rem}

\begin{rem}\label{rem-leftrightcentHopf}
Let $T$ be a centralizable Hopf monad on an autonomous category~$\cc$, with centralizer $(Z_T,\partial)$.
Set:
\begin{equation*}
\psfrag{D}[cc][cc]{\scalebox{.9}{$s^r_Y$}}
\psfrag{Y}[Bc][Bc]{\scalebox{.7}{$Y$}}
\psfrag{X}[Bc][Bc]{\scalebox{.7}{$X$}}
\psfrag{T}[Bc][Bc]{\scalebox{.7}{$T(Y)$}}
\psfrag{Z}[Bc][Bc]{\scalebox{.7}{$Z_T(X)$}}
\partial'_{X,Y}=
\rsdraw{.45}{.9}{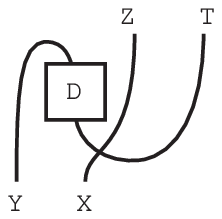} \co Y \otimes X \to Z_T(X) \otimes T(Y),
\end{equation*}
and $Z_{T^\mop}=(Z_T)^\mop$.
Then  $({Z_{T^\mop}},\partial')$ is a centralizer of $T^\mop$ in $\cc^{\otimes\opp}$.
Moreover, $Z_{T^\mop}=(Z_T)^\mop$ as Hopf monads when $Z_T$ and $Z_{T^\mop}$ are equipped with the Hopf monad structure of Theorem~\ref{doub1}.
In the language of Remark~\ref{rem-centleftright}, `left centralizability' and `right centralizability' are equivalent for a Hopf monad $T$, and a `left-handed' centralizer $Z'_T=(Z_{T\mop})^\mop$ can be identified with a `right-handed' centralizer $Z_T$ in a manner preserving the Hopf monad structures.
\end{rem}

\begin{proof}[Proof of Theorem~\ref{doub1}]
To simplify notations, set $Z=Z_T$.  Let us prove Part
(a). By definition of the product $m$ and unit $u$ of $Z$,
we have:
\begin{align*}
\bigl(\id_{T^{ \otimes 3}} \otimes m Z(m)\bigr)\partial^{3}
& =\bigl(T_2 \otimes \id_{T} \otimes m )\partial^2_{1_\cc,\otimes, 1_\cc}\\
& =\bigl((T_2 \otimes \id_{T})T_2 \otimes \id_Z )\partial_{1_\cc,\otimes_2}\\
& =\bigl((\id_{T} \otimes T_2)T_2 \otimes \id_Z )\partial_{1_\cc,\otimes_2}\\
& =\bigl(\id_{T} \otimes T_2 \otimes m )\partial^2_{1_\cc, 1_\cc,\otimes}\\
& =\bigl(\id_{T^{ \otimes 3}} \otimes m m_Z\bigr)\partial^{3}.
\end{align*}
Therefore $m Z(m)=mm_Z$ by the uniqueness assertion of
Lemma~\ref{lemcentralizer}.
 Likewise, since:
\begin{align*}
\bigl(\id_T \otimes mZ(u)\bigr)\partial &= (\id_T \otimes m)\partial_{Z,1_\cc}(u \otimes \id_{1_\cc}) \\
&= (T_0 \otimes \id_T \otimes m)\partial^2_{1_\cc,\un,1_\cc} \\
& = \bigl((T_0 \otimes \id_T)T_2(\un,-) \otimes \id_Z\bigr)\partial \\
& =(\id_T \otimes \id_Z)\partial
\end{align*}
and
\begin{align*}
(\id_T \otimes mu_Z)\partial
&= (\id_T \otimes T_0 \otimes m)\partial^2_{1_\cc,1_\cc,\un} \\
& = \bigl((\id_T \otimes T_0)T_2(-,\un) \otimes \id_Z\bigr)\partial\\
& =(\id_T \otimes \id_Z)\partial ,
\end{align*}
we get $mZ(u)=\id_Z=mu_Z$. Hence $(Z,m,u)$ is a monad on $\cc$.

Let us prove Part (b). By definition of the natural transformation $Z_2$, we have:
\begin{align*}
\bigl(\id_T & \otimes (\id_Z \otimes Z_2)Z_2\bigr)\partial_{\otimes_3,1_\cc}\\
&=(\mu T(\mu) \otimes \id_{Z^{\otimes 3}})(\partial_{T^2,1_\cc} \otimes \id_{Z^{\otimes 2}})
(\id_{1_\cc} \otimes \partial_{T,1_\cc} \otimes \id_Z)(\id_{\otimes} \otimes \partial)\\
&=(\mu \mu_T \otimes \id_{Z^{\otimes 3}})(\partial_{T^2,1_\cc} \otimes \id_{Z^{\otimes 2}})
(\id_{1_\cc} \otimes \partial_{T,1_\cc} \otimes \id_Z)(\id_{\otimes} \otimes \partial)\\
&=\bigl(\id_T \otimes (Z_2 \otimes \id_Z)Z_2\bigr)\partial_{\otimes_3,1_\cc},
\end{align*}
and so $(\id_Z \otimes Z_2)Z_2=(Z_2 \otimes \id_Z)Z_2$ by
Lemma~\ref{lemcentralizer}. Likewise, since:
\begin{align*}
\bigl(\id_T \otimes (\id_Z &\otimes Z_0)Z_2(-,\un)\bigr)\partial = (\mu \otimes \id_Z \otimes Z_0)(\partial_{1_\cc,T} \otimes \id_{Z(\un)})(\id_{1_\cc}\otimes\partial_{\un,1_\cc}) \\
&= (\mu \otimes \id_Z)\partial_{1_\cc, T}(\id_{1_\cc} \otimes \eta) =(\mu T(\eta) \otimes \id_Z)\partial=\partial
\end{align*}
and
\begin{align*}
\bigl(\id_T \otimes (Z_0 \otimes \id_Z)Z_2(\un,-)\bigr)\partial
&= (\mu \otimes Z_0 \otimes \id_Z)(\partial_{\un,T} \otimes \id_Z) \partial\\
&=(\mu \eta_T \otimes \id_Z)\partial=\partial,
\end{align*}
we get:  $(\id_Z \otimes Z_0)Z_2(1_\cc,\un)=\id_Z=(Z_0 \otimes\id_Z)Z_2(\un,1_\cc)$. Hence $Z$ is a co\-mo\-noi\-dal
functor.

Let us prove Part (c). We have to show that $m$ and $u$ are comonoidal morphisms. Since $\mu$ and $\eta$ are
comonoidal, we have:
\begin{align*}
(\id_{T^{\otimes 2}} \otimes Z_2 m)\partial^2_{\otimes,1_\cc,1_\cc}
&=(T_2\mu \otimes \id_{Z^{\otimes 2}})(\partial_{1_\cc,T} \otimes \id_Z)(\id_{1_\cc} \otimes \partial)\\
  &=\bigl((\mu \otimes \mu)T_2T(T_2) \otimes \id_{Z^{\otimes 2}} \bigr)(\partial_{1_\cc,T} \otimes \id_Z)(\id_{1_\cc} \otimes \partial)\\
&=\bigl(\id_{T^{\otimes 2}} \otimes (m \otimes m)Z_2 Z(Z_2)\bigr)\partial^2_{\otimes,1_\cc,1_\cc}
\end{align*}
and $ (\id_{T^{\otimes 2}} \otimes Z_0 m_\un)\partial^2_{\un,1_\cc,1_\cc}=T_2 \eta=\eta \otimes \eta=
\bigl(\id_{T^{\otimes 2}} \otimes Z_0 Z(Z_0)\bigr)\partial^2_{\un,1_\cc,1_\cc}$. Therefore $Z_2 m=(m \otimes m)Z_2
Z(Z_2)$ and $Z_0 m_\un=Z_0 Z(Z_0)$ by Lemma~\ref{lemcentralizer}, that is, $m$ is comonoidal. Moreover,
\begin{align*}
Z_2 u&=(T_0 \otimes Z_2)\partial_{\otimes, \un}=(T_0 \mu_\un \otimes \id_{Z^{\otimes 2}}) \partial^2_{1_\cc,1_\cc,\un} \\
&=(T_0 T(T_0) \otimes \id_{Z^{\otimes 2}}) \partial^2_{1_\cc,1_\cc,\un}
=u \otimes u
\end{align*}
and $Z_0 u_\un=(T_0 \otimes Z_0)\partial_{\un, \un}=T_0 \eta_\un =\id_\un$. Hence $u$ is comonoidal.

Parts (d) and (e) can be proved in a similar way, but we will
rather deduce them in Section~\ref{proofpartsde} from the next Theorem~\ref{mocentermon1}.
\end{proof}

\subsection{Categorical center relative to a Hopf monad} Let $T$ be a comonoidal endofunctor of a monoidal category~$\cc$.
The \emph{center of $\cc$ relative to $T$}, or shortly the \emph{$T$\trait center of $\cc$}, is the category
$\zz_T(\cc)$ defined as follows: objects are pairs $(M,\sigma)$, where $M$ is an object of $\cc$ and $\sigma\co M
\otimes 1_\cc \rightarrow T \otimes M$ is a natural transformation, such that:
\begin{center}
\psfrag{A}[Bc][Bc]{\scalebox{.7}{$T(Y)$}}
\psfrag{B}[Bc][Bc]{\scalebox{.7}{$T(Z)$}}
\psfrag{M}[Bc][Bc]{\scalebox{.7}{$M$}}
\psfrag{C}[Bc][Bc]{\scalebox{.7}{$Y\otimes Z$}}
\psfrag{Y}[Bc][Bc]{\scalebox{.7}{$Y$}}
\psfrag{Z}[Bc][Bc]{\scalebox{.7}{$Z$}}
\psfrag{s}[Bc][Bc]{\scalebox{.9}{$\sigma_{Y\otimes Z}$}}
\psfrag{r}[Bc][Bc]{\scalebox{.9}{$\sigma_Y$}}
\psfrag{u}[Bc][Bc]{\scalebox{.9}{$\sigma_Z$}}
\psfrag{o}[Bc][Bc]{\scalebox{.9}{$\sigma_\un$}}
\psfrag{T}[Bc][Bc]{\scalebox{.8}{$T_2(Y,Z)$}}
\psfrag{G}[Bc][Bc]{\scalebox{.8}{$T_0$}}
\rsdraw{.45}{.9}{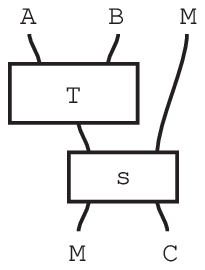} \, $=$ \, \rsdraw{.45}{.9}{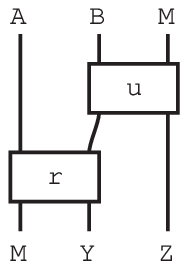} \quad and \quad
\rsdraw{.45}{.9}{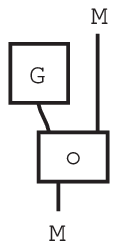} \, $=$ \, \rsdraw{.45}{.9}{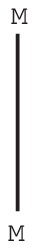} ,
\end{center}
that is,
\begin{align*}
& (T_2(Y,Z) \otimes \id_M) \sigma_{Y \otimes Z}
=(\id_{T(Y)} \otimes \sigma_Z)(\sigma_Y \otimes \id_{Z}) \quad \text{for all $Y,Z \in \Ob(\cc)$;}\\
& (T_0 \otimes \id_M) \sigma_\un =\id_M.
\end{align*}
A morphism $f \co (M,\sigma) \to (M',\sigma')$ is a morphism $f\co M \to M'$ in $\cc$ such that:
\begin{equation*}
(\id_{T(Y)} \otimes f) \sigma_Y=\sigma'_Y(f \otimes \id_Y)
\end{equation*}
for every object $Y$ of $\cc$. The composition and identities are
inherited from $\cc$.

Let $\uu_T\co \zz_T(\cc) \to \cc$ be the
forgetful functor defined by:
\begin{equation*}
\uu_T(M,\sigma)=M \quad \text{and} \quad \uu_T(f)=f.
\end{equation*}

If $\cc$ is autonomous and $T$ is a Hopf monad, then $\zz_T(\cc)$ is
autonomous. More precisely:
\begin{prop}\label{prop-monoidalZTC}
Let $(T,\mu,\eta)$ be a bimonad on a monoidal category $\cc$. Then $\zz_T(\cc)$ is monoidal, with unit object
$(\un,\eta )$ and monoidal product:
\begin{equation*}
(M,\sigma) \otimes (N,\gamma)=(M \otimes N,\rho) \quad \text{where} \quad \rho=(\mu \otimes \id_{M \otimes N})(\sigma_T \otimes \id_N)(\id_M \otimes
\gamma),
\end{equation*}
and the forgetful functor $\uu_T\co \zz_T(\cc) \to \cc$ is strict monoidal. Now assume $\cc$ is
autonomous. If $T$ has a right antipode $s^r$, then $\zz_T(\cc)$ is left autonomous with left duals given by
$\ldual{(M,\sigma)}=(\ldual{M},\sigma^l)$, where:
\begin{equation*}
\sigma^l_Y=\ldual{\bigl(}(s^r_Y \otimes \id_M)\sigma_{\rdual{T(Y)}} \bigr).
\end{equation*}
If $T$ has a left antipode $s^l$, then the category $\zz_T(\cc)$ is right autonomous with right duals given by
$\rdual{(M,\sigma)}=(\rdual{M},\sigma^r)$, where:
\begin{equation*}
\sigma^r_Y=\bigl ((s^l_Y \otimes \id_M) \sigma_{\ldual{T(Y)}} \bigr)^\vee.
\end{equation*}
In particular, if $T$ is a Hopf monad, then the category $\zz_T(\cc)$ is autonomous.
\end{prop}
We leave the proof to the reader. Pictorially, the morphisms $\rho$, $\sigma^l$, $\sigma^r$ of Proposition~\ref{prop-monoidalZTC} are:
\begin{center}
\psfrag{A}[Bc][Bc]{\scalebox{.7}{$T(Y)$}}
\psfrag{N}[Bc][Bc]{\scalebox{.7}{$N$}}
\psfrag{M}[Bc][Bc]{\scalebox{.7}{$M$}}
\psfrag{B}[Bc][Bc]{\scalebox{.7}{$\rdual{M}$}}
\psfrag{C}[Bc][Bc]{\scalebox{.7}{$\ldual{M}$}}
\psfrag{Y}[Bc][Bc]{\scalebox{.7}{$Y$}}
\psfrag{s}[Bc][Bc]{\scalebox{.9}{$\sigma_{T(Y)}$}}
\psfrag{r}[Bc][Bc]{\scalebox{.9}{$\gamma_Y$}}
\psfrag{u}[Bc][Bc]{\scalebox{.9}{$\sigma_{T(Y\rdual{)}}$}}
\psfrag{o}[Bc][Bc]{\scalebox{.9}{$\sigma_{\ldual{T}(Y)}$}}
\psfrag{m}[Bc][Bc]{\scalebox{.9}{$\mu_Y$}}
\psfrag{x}[Bc][Bc]{\scalebox{.9}{$s^r_Y$}}
\psfrag{v}[Bc][Bc]{\scalebox{.9}{$s^l_Y$}} $\rho=$
\rsdraw{.45}{.9}{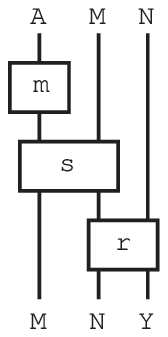} , \quad $\sigma^l=$
\rsdraw{.45}{.9}{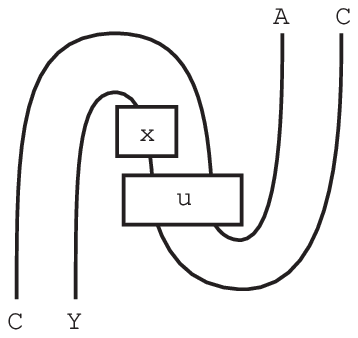}, \quad and \quad $\sigma^r=$
\rsdraw{.45}{.9}{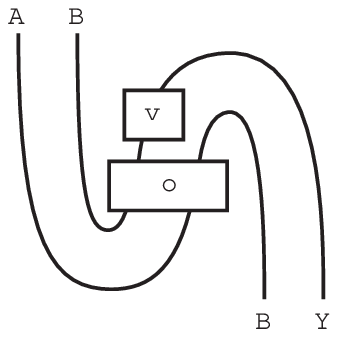} .
\end{center}

\begin{rem}\label{remcenterusual}
If $\cc$ is an autonomous category, then $\zz_{1_\cc}(\cc)$ coincides with the usual center $\zz(\cc)$ of $\cc$ (see Section~\ref{sect-centerusual}).
\end{rem}

\begin{rem}\label{remcenterleft}
The definition of the category $\zz_T(\cc)$ is not left/right symmetric. One may also consider the category $\zz'_T(\cc) = \zz_{T^\mop}(\cc^{\otimes \opp})^{\otimes\opp}$,
whose objects are pairs $(M,\sigma)$, where $M$ is an object of $\cc$ and $\sigma\co 1_\cc \otimes M \to M \otimes T$ is a natural transformation satisfying the obvious conditions.
If $\cc$ is autonomous and $T$ is a Hopf monad, then the category $\zz'_T(\cc)$ is autonomous and isomorphic to $\zz_T(\cc)$ via the strict monoidal functor $\zz_T(\cc) \to \zz'_T(\cc)$ defined
by $(M,\sigma) \mapsto (M,\sigma')$, where: $$\sigma'_Y=(\rev_Y (\id_Y \otimes s^r_Y)\otimes \id_{M \otimes T(Y)})
(\id_Y \otimes \sigma_{\rdual{T(Y)}} \otimes \id_{T(Y)})(\id_{Y \otimes M} \otimes \rcoev_{T(Y)}).$$
In particular $\zz'_{1_\cc}(\cc)=\zz'(\cc)$, see Remark~\ref{rem-alter-center}.
\end{rem}

\subsection{Monadicity of centers}
In this section, we show that the center relative to a centralizable Hopf monad is monoidally equivalent to the category of modules of a its centralizer.

\begin{thm}\label{mocentermon1}
Let $T$ be a centralizable comonoidal endofunctor of a left autonomous category~$\cc$, with centralizer
$(Z_T,\partial)$. The functor $E\co Z_T \ti \cc \to \zz_T(\cc)$, defined by:
\begin{equation*}
E(M,r)=\bigl(M,(\id_T \otimes r) \partial_{M,1_\cc}\bigr) \quad \text{and} \quad E(f)=f,
\end{equation*}
is an isomorphism of categories such that the following triangle commutes:
\begin{equation*}
\xymatrix @!0 @C=3pc @R=2.5pc {Z_T \ti \cc \ar[rr]^(.47)E \ar[dr]_{U_{Z_T}} & \ar@{}[d]|(.4){\circlearrowright} & \zz_T(\cc) \ar[dl]^{\uu_T}\\
& \cc & }
\end{equation*}
Furthermore, if $T$ is a bimonad, so that $Z_T$ is a bimonad and
$\zz_T(\cc)$ is monoidal, then $E$ is strict monoidal (and so $\uu_T
E=U_{Z_T}$ as monoidal functors).
\end{thm}
We prove Theorem~\ref{mocentermon1} in Section~\ref{sect-proof-center}.

\begin{rem}
The functor $\ff_T=EF_{Z_T}\co
\cc \to \zz_T(\cc)$ is left adjoint to $\uu_T$ and the adjunction
$(\ff_T,\uu_T)$ is monadic with monad $Z_T$ (see Remark~\ref{rem-def-monadic}). If $T$ is a bimonad, this adjunction is monoidal and $Z_T$ is its associated bimonad (see Theorem~\ref{thm-gen-ex}).
\end{rem}

A monoidal category~$\cc$ is
said to be \emph{centralizable} if its identity endofunctor $1_\cc$ is centralizable. In such case, the centralizer of
$1_\cc$ is called the \emph{centralizer} of $\cc$. In view of Remark~\ref{remcenterusual}, we have:

\begin{cor}\label{cormoncenter}
Let $\cc$ be a centralizable autonomous category, with centralizer $(Z,\partial)$.
Then the forgetful functor $\uu\co
\zz(\cc) \to \cc$ is monadic with monad $Z$.
In fact $Z$ is a Hopf monad and the functor $Z\ti \cc \to \zz(\cc)$, defined by:
\begin{equation*}
(M,r)\mapsto \bigl(M,(\id_{1_\cc} \otimes r)
\partial_{M,1_\cc}\bigr) \quad \text{and} \quad f \mapsto f,
\end{equation*}
is a strict monoidal isomorphism of categories.
\end{cor}

\begin{rem}
The monadicity assertion of Corollary~\ref{cormoncenter}
is a consequence of~\cite[Theorem~4.3]{DayStreet}.
\end{rem}

\begin{rem}
We will see in Section~\ref{sect-doubleHM} that $R=\bigl ( u \otimes \id \bigr)\partial$ is an \Rt matrix for~$Z$
(where $u$ denotes the unit of $Z$), making the isomorphism of Corollary~\ref{cormoncenter} an isomorphism of braided
categories.
\end{rem}

\subsection{Proof of Theorem~\ref{mocentermon1}} \label{sect-proof-center}
Throughout this section, let $T$ be a centralizable endofunctor of a left autonomous category~$\cc$, with centralizer
$(Z_T,\partial)$. Recall $Z_T$ is a monad by Theorem~\ref{doub1}(a). Denote $m$ and $u$ its product and unit.

Remark first that, by Lemma~\ref{lemcentralizer}, for any object $M$ of $\cc$, we have a bijection:
\begin{equation*}
\left \{
\begin{array}{ccl}
\Hom_\cc(Z_T(M),M) & \to & \Nat(M\otimes 1_\cc, T \otimes M) \\
r & \mapsto & \sigma_{(M,r)}=\{(\id_{T(Y)}\otimes r)\partial_{M,Y}\}_{Y \in \Ob(\cc)}
\end{array}\right..
\end{equation*}
\begin{lem} \label{lemproofmon1}
Let $M$ be an object of $\cc$ and $r\co Z_T(M) \to M$ be a morphism in $\cc$. Then $(M,r)$ is a $Z_T$\ti module if and
only if $(M,\sigma_{(M,r)})$ is an object of $\zz_T(\cc)$.
\end{lem}
\begin{proof}
By definition of the multiplication $m$ of $Z_T$, we have:
\begin{equation*}
(T_2 \otimes \id_M)(\sigma_{(M,r)})_\otimes=(T_2 \otimes r)\partial_{M,\otimes}=(\id_{T^{\otimes 2}} \otimes
rm_M)\partial^2_{M,1_\cc,1_\cc}.
\end{equation*}
Moreover:
\begin{equation*}
(\id_T \otimes \sigma_{(M,r)})(\sigma_{(M,r)} \otimes \id_{1_\cc})= \bigl(\id_{T^{\otimes 2}} \otimes
rZ_T(r)\bigr)\partial^2_{M,1_\cc,1_\cc}.
\end{equation*}
Therefore, by Lemma~\ref{lemcentralizer}, $(T_2 \otimes \id_M)(\sigma_{(M,r)})_\otimes=(\id_T \otimes
\sigma_{(M,r)})(\sigma_{(M,r)} \otimes \id_{1_\cc})$ if and only if $rm_M=rZ_T(r)$. Also, since $(T_0 \otimes \id_M)
(\sigma_{(M,r)})_\un=(T_0 \otimes r) \partial_{M,\un}=ru_M$, we have $(T_0 \otimes \id_M) (\sigma_{(M,r)})_\un=\id_M$
if and only if $ru_M=\id_M$.
\end{proof}
\begin{lem} \label{lemproofmon2}
Let $(M,r)$ and $(N,s)$ be two $Z_T$\ti modules. Let $f\co M \to N$ be a morphism in $\cc$. Then $f$ is $Z_T$\ti linear
if and only if it is a morphism  from $(M,\sigma_{(M,r)})$ to $(N,\sigma_{(N,s)})$ in $\zz_T(\cc)$.
\end{lem}
\begin{proof}
We have: $(\id_T \otimes f)\sigma_{(M,r)}=(\id_T \otimes fr)\partial_{M,1_\cc}$ and $$\sigma_{(N,s)}(
f\otimes\id_T)=(\id_T \otimes s)\partial_{N,1_\cc}(f\otimes \id_{1_\cc})=(\id_T \otimes sZ_T(f))\partial_{M,1_\cc}.$$
Therefore, by Lemma~\ref{lemcentralizer}, we obtain: $(\id_T \otimes f)\sigma_{(M,r)}=\sigma_{(N,s)}( f\otimes\id_T)$
if and only if $fr=sZ_T(f)$.
\end{proof}
Using Lemmas~\ref{lemproofmon1} and~\ref{lemproofmon2}, one sees
that the functor $E\co Z_T \ti \cc \to \zz_T(\cc)$, given by
$E(M,r)=(M,\sigma_{(M,r)})$ and $E(f)=f$, is a well-defined
isomorphism of categories. Furthermore it clearly satisfies $\uu_T
E=U_{Z_T}$.

Assume now that $(T,\mu,\eta)$ is a bimonad. Then $Z_T$ is a bimonad by Theorem~\ref{doub1}(c) and the category
$\zz_T(\cc)$ is monoidal by Proposition~\ref{prop-monoidalZTC}. Since, for all $Z_T$\ti modules $(M,r)$ and $(N,s)$, we
have:
\begin{align*}
E(M,r) \otimes E(N,s)
& =(M,\sigma_{(M,r)}) \otimes (N,\sigma_{(N,s)})\\
& =\bigl (M \otimes N,(\mu \otimes r \otimes s)(\partial_{M,T} \otimes \id_{Z_T(N)})(\id_M \otimes
\partial_{N,1_\cc})\bigr)\\
& =\bigl (M \otimes N,(\id_T\otimes (r \otimes s)(Z_T)_2(M,N))\partial_{M\otimes N,1_\cc}\bigr)\\
&=E\bigl ((M,r) \otimes (N,s)\bigr)
\end{align*}
and $E\bigl (\un,(Z_T)_0\bigr)=\bigl (\un, (\id_T \otimes (Z_T)_0)\partial_{\un,1_\cc}\bigr)=(\un,\eta)$, the functor
$E$ is strict monoidal.
Finally, we have: $\uu_T E=U_{Z_T}$ as monoidal functors because  the forgetful functors $U_{Z_T}\co Z_T\ti\cc \to \cc$ and $\uu_T\co \zz_T(\cc)
\to \cc$ are strict monoidal.

\subsection{End of proof of Theorem~\ref{doub1}}\label{proofpartsde}
Let us prove Part (d) of Theorem~\ref{doub1}. Let $(T,\mu,\eta)$ be
a centralizable bimonad on an autonomous category~$\cc$, with
centralizer $(Z_T,\partial)$.
By Theorem~\ref{doub1}(c), $Z_T$ is a bimonad. By Theorem~\ref{mocentermon1}, the functor $E\co Z_T \ti \cc \to \zz_T(\cc)$, defined by:
\begin{equation*}
E(M,r)=\bigl(M,(\id_T \otimes r) \partial_{M,1_\cc}\bigr) \quad \text{and} \quad E(f)=f,
\end{equation*}
is a strict monoidal isomorphism. Assume~$T$ admits a right antipode $s^r$. Then the category $\zz_T(\cc)$ is
left autonomous by Proposition~\ref{prop-monoidalZTC}. Hence the category $Z_T\ti \cc$ is
left autonomous, and so $Z_T$ admits a left antipode by
Theorem~\ref{thm-biHopfmon}(b).
Denote $m$ the product  of $Z_T$, $u$ its unit, and $S^l$ its right
antipode. Let $X$ be an object of $\cc$. In the category $Z_T\ti\cc$, we have a duality:
\begin{equation*}
\bigl(\ldual{(}Z_T(X),m_X),(Z_T(X),m_X), \lev_{Z_T(X)},\lcoev_{Z_T(X)}\bigr),
\end{equation*}
where $\ldual{(}Z_T(X),m_X)=\bigl(\ldual{Z}_T(X),S^l_{Z_T(X)}Z_T(\ldual{m}_T)
\bigr)$. Hence, $E$ being strict monoidal, a duality in
the category $\zz_T(\cc)$:
\begin{equation*}
\bigl(E\bigl(\ldual{(Z_T(X),m_X)}\bigr),E(Z_T(X),m_X), \lev_{Z_T(X)},\lcoev_{Z_T(X)}\bigr),
\end{equation*}
where
$E\bigl(\ldual{(Z_T(X),m_X)}\bigr)=\bigl(\ldual{Z}_T(X),(\id_T \otimes
S^l_{Z_T(X)}Z_T(\ldual{m}_T))\partial_{\ldual{Z}_T(X),1_\cc}\bigr)$.
Now, by Proposition~\ref{prop-monoidalZTC}, we also have the following duality in $\zz_T(\cc)$:
\begin{equation*}
\bigl(\ldual{E}(Z_T(X),m_X),E(Z_T(X),m_X), \lev_{Z_T(X)},\lcoev_{Z_T(X)}\bigr)
\end{equation*}
where
$
\ldual{E}(Z_T(X),m_X)=\bigr(\ldual{Z_T(X)},\ldual{(}(s^r \otimes
m_X)\partial_{Z_T(X),T^\vee})\bigr)
$.
Hence, by uniqueness of duals up to unique isomorphism:
\begin{equation*}
\bigl(\id_T \otimes S^l_{Z_T(X)}Z_T(\ldual{m}_T)\bigr)
\partial_{\ldual{Z}_T(X),1_\cc}= \ldual{\bigl(}(s^r \otimes
m_X)\partial_{Z_T(X),T^\vee}\bigr).
\end{equation*}
Composing on the left with $(\id_T \otimes \ldual{u}_X)=\ldual{(}
u_X \otimes \id_{T^\vee})$, we get:
\begin{equation*}
(\id_T \otimes S^l_X)\partial_{\,\ldual{Z_T(X)}, 1_\cc}=\ldual{\bigl((s^r \otimes \id_{Z_T(X)})\partial_{X,
\rdual{T}}\bigr)},
\end{equation*}
which is the defining relation of
Theorem~\ref{doub1}(d).   Hence  Part (d) of Theorem~\ref{doub1}. Part
(e) can be shown similarly.

\section{The double of a Hopf monad}\label{sect-double}
Given a centralizable Hopf monad $T$ on an autonomous category $\cc$, we construct the canonical distributive
law $\Omega$ of $T$ over its centralizer $Z_T$, which serves two purposes.
Firstly $\Omega$ gives rise to a new Hopf monad $D_T=Z_T \circ_\Omega T$, called the double of $T$. The double $D_T$ is
actually quasitriangular  and $Z(T\ti\cc)\simeq D_T\ti\cc$ as braided categories,
see Section~\ref{sect-doubleHM}.
Secondly $\Omega$ defines a lift of the Hopf monad $Z_T$ to a Hopf monad $\Tilde{Z}_T^\Omega$ on $T\ti\cc$, which turns
out to be the centralizer of the category~$T\ti\cc$, and so $\Tilde{Z}_T^\Omega(\un,T_0)$ is the coend of $T\ti\cc$, see Section~\ref{sect-centralizerTC}.

Most of the results of this section are special cases of results of Section~\ref{sect-cent-mod}. We state them here for
convenience.

\subsection{The canonical distributive law}\label{sect-conodistlaw}
Let $T$ be a centralizable Hopf monad on an autonomous category $\cc$ and $(Z_T,\partial)$ be its centralizer.

Recall (see Proposition~\ref{propcentral1}) that $ Z_T(X)=\int^{Y \in \cc} \ldual{T(Y)} \otimes X \otimes Y$, with
universal dinatural transformation:
$$
i_{X,Y}=(\lcoev_{T(Y)} \otimes \id_{Z_T(X)})(\id_{\ldual{T}(X)} \otimes \partial_{X,Y}),
$$
which is natural in $X$ and dinatural in $Y$. Since $T(i)$ is a universal dinatural trans\-formation (see
Section~\ref{sect-coendHopfmon}), we can define a natural transformation $\Omega\co TZ_T \to Z_TT$ by:
\begin{equation*}
\Omega_X T(i_{X,Y})= i_{T(X),T(Y)}\bigl(\ldual{\mu_Y}s^l_{T(Y)}T(\ldual{\mu_Y}) \otimes \id_{T(X)\otimes
T(Y)}\bigr)T_3\bigl(\ldual{T(Y)},X, Y\bigr),
\end{equation*}
where $\mu$ and $s^l$ are the product and left antipode of $T$ and
$T_3\co T\otimes_3 \to T^{\otimes 3}$ is defined as in
Section~\ref{sect-comonofunctor}.

\begin{thm}\label{thm-candistlaw}
The natural transformation $\Omega\co TZ_T \to Z_TT$ is an invertible co\-mon\-oidal distributive law.
\end{thm}
We call $\Omega$ the \emph{canonical distributive law of $T$}. We
prove Theorem~\ref{thm-candistlaw} in
Section~\ref{sect-proof-thm-candistlaw}.

The inverse $\Omega^{-1}\co Z_TT \to TZ_T$ of the distributive law $\Omega$ is the natural transformation  defined by:
\begin{align*}
\Omega&^{-1}_X i_{T(X),Y}
=\bigl(\lev_{T(Y)}(\id_{\ldual{T}(Y)} \otimes \mu_Y T(\mu_Y) ) \otimes  T(i_{X,T(Y)}) \otimes \lev_Y
 (s^l_Y \otimes \id_Y)\bigr)\\
&\circ T_3(T^2(Y),\ldual{T}^2(Y) \otimes X \otimes T(Y),\ldual{T}(Y)) T(\lcoev_{T^2(Y)}\otimes \id_X \otimes \lcoev_{T(Y)}).
\end{align*}

\begin{rem}
The canonical distributive law of $T$ is the only natural transformation $\Omega\co TZ_T \to Z_TT$ satisfying:
\begin{equation*}
(\mu \otimes \Omega)T_2T(\partial)= (\mu \otimes
\id_{Z_TT})\partial_{T,T} T_2.
\end{equation*}
\end{rem}

\begin{rem} \label{rem-rmatf}
One can show that \Rt matrices for $T$ correspond bijectively with morphisms of Hopf  monads $f\co Z_T \to T$ satisfying $\mu T(f)=\mu f_T \Omega$.
The \Rt matrix associated with such a morphism $f$ is $R= (\id_T \otimes f)\partial$.
\end{rem}

\subsection{The double of a Hopf monad}\label{sect-doubleHM}
Let $T$ be a centralizable Hopf monad on an autonomous category~$\cc$, with centralizer $(Z_T,\partial)$. Let
$\Omega\co TZ_T \to Z_TT$ be the canonical distributive law of $T$. By Corollary~\ref{cordistribHopf},
\begin{equation*}
D_T=Z_T \circ_{\Omega} T,
\end{equation*}
is a Hopf monad on $\cc$. Denote $\eta$ and $u$ the units of $T$ and
$Z_T$ respectively.

\begin{thm}\label{thm-Rmat-double}
The natural transformation $R=\{R_{X,Y}\}_{X,Y \in \Ob(\cc)}$,
defined by:
\begin{equation*}
R_{X,Y}=\bigl ( u_{T(Y)} \otimes Z_T(\eta_X)\bigr)\partial_{X,Y} \co
X \otimes Y \to D_T(Y) \otimes D_T(X),
\end{equation*}
is \Rt matrix for the Hopf monad $D_T$.
\end{thm}
The quasitriangular Hopf monad $D_T$ is called the \emph{double of
$T$}. This terminology is justified by the fact that the braided
categories $Z(T\ti\cc)$ and $D_T\ti\cc$ coincide. More precisely,
let $\uu\co \zz(T\ti\cc) \to \cc$ be the strict monoidal forgetful
functor defined by:
\begin{equation*}
\uu\bigl((M,r),\sigma\bigr)=M \quad \text{and} \quad \uu(f)=f.
\end{equation*}
Let $I\co D_T \ti \cc \to \zz(T\ti\cc)$ be the functor defined by
$I(f)=f$ and:
\begin{equation*}
I(M,r)=\bigl((M,ru_{T(M)}),\sigma\bigr) \quad \text{with} \quad
\sigma_{(N,s)}=(s \otimes rZ_T(\eta_M))\partial_{M,N}.
\end{equation*}
\begin{thm}\label{thm-doublable}
The functor $I$ is a strict monoidal isomorphism of braided
categories such that the following triangle of monoidal functors
commutes:
\begin{equation*}
\xymatrix @!0 @C=3pc @R=2.5pc {D_T \ti \cc \ar[rr]^(.47)I \ar[dr]_{U_{D_T}}& \ar@{}[d]|(.4){\circlearrowright} & \zz(T\ti\cc) \ar[dl]^{\uu}\\
& \cc & }
\end{equation*}
\end{thm}

We prove Theorems~\ref{thm-Rmat-double} and \ref{thm-doublable} in
Section~\ref{sect-proof-thm-doublable}.

\begin{rem}
The functor $\ff=IF_{D_{T}}\co
\cc \to \zz(T\ti\cc)$ is left adjoint to $\uu$ and the adjunction
$(\ff,\uu)$ is monadic with monad $D_{T}$ (see Remark~\ref{rem-def-monadic}). Moreover $D_T$ is the Hopf monad associated with this monoidal adjunction (see Theorem~\ref{thm-gen-ex}).
\end{rem}

\begin{rem}
According to Remark~\ref{rem-centleftright}, the construction of the double of a Hopf monad $T$ admits a `right-handed' version: if $Z'_T$ is a `right-handed' centralizer of~$T$, there exists a `right-handed' canonical law $\Omega'$ of $T$ over $Z'_T$, and hence a Hopf monad $D'_T=Z'_T \circ_{\Omega'} T$ endowed with an \Rt matrix $R'$ such that $D'_T\ti \cc\simeq \zz'(T\ti\cc)$ as braided category. If we identify $Z'_T$ to $Z_T$ as in Remark~\ref{rem-leftrightcentHopf}, then $\Omega'=\Omega$, $D'_T=D_T$ as Hopf monads, and $R'=R^{* -1}$.
\end{rem}

\begin{rem}\label{rem-dri-caract}
Let $T$ be a centralizable Hopf monad on an autonomous category $\cc$ and $(Z_T,\partial)$ be its centralizer. Denote
$\eta$ and $u$ the units of $T$ and $Z_T$ respectively. Assuming $u_T\co T \to Z_TT$ is a monomorphism, one can show
that the canonical distributive law of~$T$ is the only comonoidal distributive law $\Omega\co TZ_T \to Z_TT$ such that:
\begin{equation*}
R=\bigl ( u_{T} \otimes Z_T(\eta)\bigr)\partial
\end{equation*}
is an \Rt matrix for the Hopf monad $Z_T \circ_\Omega T$. This
generalizes Drinfeld's original characterization of the double of
a finite-dimensional Hopf algebra.
\end{rem}

\subsection{The centralizer and the coend of a category of modules}\label{sect-centralizerTC}
Let~$T$ be a centralizable Hopf monad on an autonomous category~$\cc$. Let $(Z_T,\partial)$ be the centralizer of $T$
and
$\Omega\co TZ_T \to Z_TT$ be the canonical distributive law of $T$.
By Corollary~\ref{cordistribHopf}, $\Tilde{Z}_T^\Omega$ is Hopf monad which is a lift of the Hopf monad $Z_T$ to
$T\ti\cc$. Recall:
\begin{equation*}
\Tilde{Z}_T^\Omega(M,r)=(Z_T(M), Z_T(r)\Omega_M) \quad \text{and} \quad \Tilde{Z}_T^\Omega(f)=Z_T(f).
\end{equation*}
For any $T$-modules $(M,r)$ and $(N,s)$, set:
\begin{equation*}
\Tilde{\partial}_{(M,r),(N,s)}=(s\otimes\id_{Z_T(M)})\partial_{M,N}\co (M,r) \otimes (N,s) \to (N,s) \otimes
\Tilde{Z}_T^\Omega(M,r).
\end{equation*}
\begin{thm}\label{thm-centraZ}
The pair $(\Tilde{Z}_T^\Omega,\Tilde{\partial})$ is a centralizer of the category $T\ti\cc$.
\end{thm}
We prove Theorem~\ref{thm-centraZ} in Section~\ref{sect-proof-centraZ}.

Recall that:
\begin{equation*}
Z_T(\un)=\int^{Y \in \cc} \ldual{T}(Y) \otimes Y,
\end{equation*}
with universal dinatural transformation
$i_Y=(\lev_{T(Y)} \otimes \id_{Z_T(\un)})\partial_{\un,Y}$. Denote $\alpha=Z_T(T_0)\Omega_\un$  the $T$\ti action of
$\Tilde{Z}_T^\Omega(\un,T_0)$. It is characterized by:
$$
\alpha T(i_Y)=i_{T(Y)}\bigl(\ldual{\mu_Y}s^l_{T(Y)}T(\ldual{\mu_Y}) \otimes \id_{T(Y)}\bigr)T_2\bigl(\ldual{T(Y)}, Y\bigr).
$$
By Theorem~\ref{thm-centraZ} and Proposition~\ref{propcentral1}, $\Tilde{Z}_T^\Omega(\un,T_0)=(Z_T(\un),\alpha)$ is the coend of $T\ti\cc$, that is:
$$(Z_T(\un),\alpha)=\int^{(M,r) \in T\ti\cc}\hspace*{-3em} \ldual{(}M,r) \otimes (M,r),$$ with universal dinatural transformation $\Tilde{\imath}_{(M,r)}=i_M(\ldual{r} \otimes M)$.

\begin{figure}[t]
\begin{center}
\psfrag{X}[Bc][Bc]{\scalebox{.8}{$X$}} \psfrag{Y}[Bc][Bc]{\scalebox{.8}{$Y$}}
\psfrag{M}[Bc][Bc]{\scalebox{.8}{$\ldual{T}(X)$}} \psfrag{W}[Bc][Bc]{\scalebox{.8}{$\ldual{T}(Y)$}}
 \psfrag{t}[Bc][Bc]{\scalebox{.8}{$T_2(T(X),Y)$}}
 \psfrag{d}[Bc][Bc]{\scalebox{.8}{$\partial_{\un,T(X) \otimes Y}$}}
 \psfrag{a}[Bc][Bc]{\scalebox{.8}{$\ldual{\mu_X}s^l_{T(X)}T(\ldual{\mu_X})$}}
  \psfrag{T}[Bc][Bc]{\scalebox{.8}{$T_2(\ldual{T}(X),X)$}}
   \psfrag{Z}[Bc][Bc]{\scalebox{.8}{$Z_T(\un)$}}
 \psfrag{R}[Bc][Bc]{\scalebox{.8}{$R_{\ldual{T}(X) \otimes X,\ldual{T}(Y)}$}}
\psfrag{D}[cc][cc]{\scalebox{.9}{$s^l_{T(Y)}T(\ldual{\mu_Y})$}}
$m(i_X \otimes i_Y)=$ \rsdraw{.45}{1}{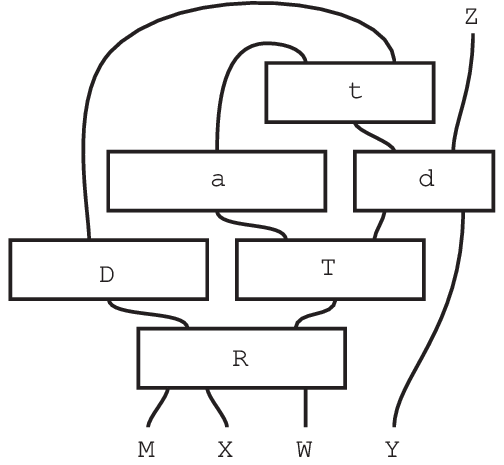}\,, \\
 \psfrag{Y}[Bc][Bc]{\scalebox{.8}{$Y$}}
\psfrag{W}[Bc][Bc]{\scalebox{.8}{$\ldual{T}(Y)$}}
 \psfrag{d}[Bc][Bc]{\scalebox{.8}{$\partial_{\un,\ldual{T}(Y)}$}}
 \psfrag{n}[Bc][Bc]{\scalebox{.8}{$\mu_Y$}}
  \psfrag{a}[Bc][Bc]{\scalebox{.8}{$\alpha$}}
   \psfrag{Z}[Bc][Bc]{\scalebox{.8}{$Z_T(\un)$}}
 \psfrag{R}[Bc][Bc]{\scalebox{.8}{$R_{Z_T(\un),T(Y)}$}}
\psfrag{D}[cc][cc]{\scalebox{.9}{$s^l_{T(Y)}T(\ldual{\mu_Y})$}}
$S i_Y=$ \rsdraw{.45}{1}{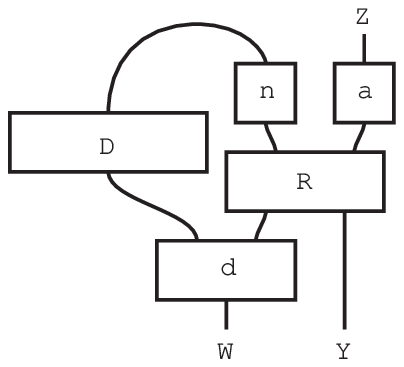}\, , \qquad
\psfrag{X}[Bc][Bc]{\scalebox{.8}{$X$}} \psfrag{Y}[Bc][Bc]{\scalebox{.8}{$Y$}}
\psfrag{M}[Bc][Bc]{\scalebox{.8}{$\ldual{T}(X)$}} \psfrag{W}[Bc][Bc]{\scalebox{.8}{$\ldual{T}(Y)$}} \psfrag{R}[Bc][Bc]{\scalebox{.8}{$R_{X,\ldual{T}(Y)}$}}
\psfrag{U}[Bc][Bc]{\scalebox{.8}{$R_{T(\ldual{T}(Y)),T(X)}$}}
\psfrag{u}[Bc][Bc]{\scalebox{.8}{$\mu_X$}}
\psfrag{D}[cc][cc]{\scalebox{.9}{$s^l_Y\mu_{\ldual{T}(Y)}$}}
$\omega(i_X \otimes i_Y)=$ \rsdraw{.45}{1}{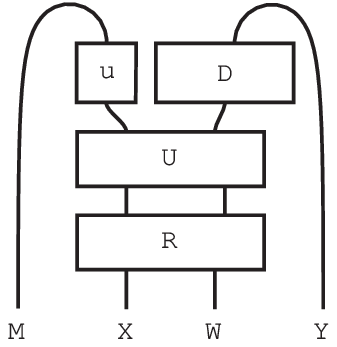}\,.
\end{center}
\caption{Hopf algebra structure of the coend of $T\ti\cc$}
\label{fig-coendTC}
\end{figure}

The coend $(Z_T(\un),\alpha)$ of $T\ti\cc$ is a coalgebra in $T\ti\cc$, with coproduct  and counit given by:
\begin{equation*}
\Delta=(Z_T)_2(\un,\un)\co Z_T(\un)\to Z_T(\un) \otimes Z_T(\un) \quad \text{and} \quad \varepsilon=(Z_T)_0 \co Z_T(\un) \to\un.
\end{equation*}
Assume now that $T$ is furthermore quasitriangular, with \Rt matrix $R$, so that the autonomous category $T\ti\cc$ is braided. Then the coalgebra $\bigl((Z_T(\un),\alpha),\Delta,\varepsilon\bigr)$ becomes a Hopf algebra in~$T\ti\cc$ endowed with a self-dual Hopf pairing (see Section~\ref{sect-coendbraided}). Its unit is:
\begin{equation*}
u=(T_0 \otimes \id_{Z_T(\un)})\partial_{\un,\un} \co \un \to Z_T(\un).
\end{equation*}
Its product $m$, antipode $S$, and Hopf pairing $\omega$ are given in Figure~\ref{fig-coendTC}.

\begin{rem}
In Section~\ref{sect-coendcenterfusion}, we treat the case of the centralizer
of a fusion category~$\ff$ (which is a quasitriangular Hopf monad by Theorem~\ref{thm-doublable}) to get a convenient description of the coend
of $\zz(\ff)$.
\end{rem}

\section{The centralizer of a Hopf monad on a category of modules}\label{sect-cent-mod}

In this section, we study the centralizer of a Hopf monad $Q$ on the category $T\ti\cc$ of modules over a Hopf monad $T$ on an autonomous category $\cc$. We show that it is centralizable whenever the cross product $Q\cp T$ is centralizable.
In that case, the centralizer of $Q\cp T$ lifts naturally to a centralizer of $Q$, which turns out to be also a lift of Hopf monads.  Hence a canonical distributive law $\Omega$ of $T$ over $Z_{Q\cp T}$ and a Hopf monad  $D_{Q,T}=Z_{Q\cp T} \circ_\Omega T$ on $\cc$. We interpret the category of $D_{Q,T}$ modules as the center of $T\ti\cc$ relative to $Q$.

\subsection{Centralizability on categories of modules}\label{sect-lift1}
In this section, given a Hopf monad $T$ on an autonomous category $\cc$, we give a criterion for an
endofunctor~$Q$ of~$T\ti\cc$ to be centralizable in terms of the centralizability of the cross product $Q\cp  T$ on $\cc$
(see Section~\ref{sect-crossprod} for the definition of cross-products).

\begin{prop} \label{lem-lift-cent1}
Let $T$ be a Hopf monad on an autonomous category $\cc$ and let~$Q$ be a endofunctor of $T\ti\cc$. Let $(M,r)$ be a
$T$\ti module. Then:
\begin{enumerate}
\renewcommand{\labelenumi}{{\rm (\alph{enumi})}}
\item  The endofunctor $Q$ is centralizable at $(M,r)$ if and only if $Q\cp  T$ is centralizable at $U_T(M,r)=M$.
\item   Assume $Q\cp  T$ is centralizable at $M$, with centralizer $(Z,\delta)$.
    Then $Q$ admits a unique centralizer
$(\Tilde{Z},\Tilde{\delta})$ at $(M,r)$ such that:
\begin{equation*}
U_T(\Tilde{Z})=Z \quad\mbox{and}\quad  \Tilde{\delta}_{(N,s)}  =(Q(s) \otimes \id_Z)\delta_{N}
\end{equation*}
for any $T$\ti module $(N,s)$.
\end{enumerate}
\end{prop}
\begin{rem}
In the second formula of Proposition~\ref{lem-lift-cent1}(b), $Q(s)$ makes
sense because $s\co(T(N),\mu_N) \to (N,s)$ is a morphism in $T\ti\cc$.  This formula can be written:
\begin{equation*}
 \Tilde{\delta}=(Q(\varepsilon) \otimes \id_Z)\delta
\end{equation*}
where $\varepsilon$ denotes the counit of the adjunction $(U_T,F_T)$.
\end{rem}
\begin{proof}[Proof of Proposition~\ref{lem-lift-cent1}]
Let us prove Part (a). Fix a $T$\ti module $(M,r)$. By Proposition~\ref{propcentral1}, $Q$ is centralizable at $(M,r)$
if and only if the coend:
\begin{equation*}
\int^{(N,s) \in T\ti\cc} \hspace*{-2.8em}\ldual{Q}(N,s) \otimes (M,r) \otimes (N,s)
\end{equation*}
exists. Since the functor $U_T$ creates and preserves coends (see Section~\ref{sect-coendHopfmon}) and is strict
monoidal, this is equivalent to the existence of the coend:
\begin{equation*}
\int^{(N,s) \in T\ti\cc} \hspace*{-2.8em}U_T\bigl(\ldual{Q}(N,s) \otimes (M,r) \otimes (N,s)\bigr)=\int^{(N,s) \in
T\ti\cc} \hspace*{-2.8em} \ldual{U_T}Q(N,s) \otimes M \otimes U_T(N,s).
\end{equation*}
By Lemma~\ref{lem-coend-adj}, this is equivalent to the existence of the coend:
\begin{equation*}
\int^{Y \in \cc} \hspace*{-2em} \ldual{U_T}QF_T(Y) \otimes M \otimes Y=\int^{Y \in \cc} \hspace*{-2em} \ldual{Q}\cp  T(Y)
\otimes M \otimes Y,
\end{equation*}
and so, by Proposition~\ref{propcentral1}, to the fact that $Q \cp  T$ is centralizable at $M$.

Let us prove Part (b). By Proposition~\ref{propcentral1}, we have:
\begin{equation*}
Z=\int^{Y \in \cc} \hspace*{-2em} \ldual{Q}\cp  T(Y) \otimes M \otimes Y,
\end{equation*}
with universal dinatural transformation $i_Y=(\lev_{Q\cp  T(Y)} \otimes \id_Z)(\id_{\ldual{Q}\cp  T(Y)} \otimes \delta_Y)$.
By Lemma~\ref{lem-coend-adj}, we have also:
\begin{equation*}
Z=\int^{(N,s) \in T\ti\cc} \hspace*{-2.8em} \ldual{U_T}Q(N,s) \otimes M \otimes U_T(N,s).
\end{equation*}
with universal dinatural transformation $j_{(N,s)}=i_N(\ldual{U_T}Q(s) \otimes \id_{M \otimes N})$. Set:
\begin{equation*}
\Tilde{\delta}_{(N,s)}= (\id_{Q(N,s)} \otimes j_{(N,s)})(\lcoev_{Q(N,s)} \otimes \id_M ).
\end{equation*}
By Proposition~\ref{HM-creat-colim}, there exists a unique $T$\ti action $\alpha\co T(Z) \to Z$ such that $j_{(N,s)}$,
or equivalently $\Tilde{\delta}_{(N,s)}$, is $T$\ti linear for all $T$\ti modules $(N,s)$. Furthermore we have:
\begin{equation*}
(Z,\alpha)=\int^{(N,s) \in T\ti\cc} \hspace*{-2.8em}\ldual{Q}(N,s) \otimes (M,r) \otimes (N,s)
\end{equation*}
with universal dinatural transformation $j$. Set $\Tilde{Z}=(Z,\alpha)$. By Proposition~\ref{propcentral1},
$(\Tilde{Z}, \Tilde{\delta})$ is a  centralizer of $Q$ at $(M,r)$.  By construction, we have $U_T(\Tilde{Z})=Z$ and
$\Tilde{\delta}_{(N,s)}=(U_TQ(s) \otimes \id_Z)\delta_{N}$ for every $T$\ti module $(N,s)$. Furthermore, since
$\alpha$ is the only action of $T$ on $Z=U_T(\Tilde{Z})$ such that every $\Tilde{\delta}_{(N,s)}$ is $T$\ti
linear, $(\Tilde{Z}, \Tilde{\delta})$ is the only centralizer of $Q$ at $(M,r)$ satisfying the conditions of Part~(b).
\end{proof}

Applying Lemma~\ref{lem-centralisable} and Proposition~\ref{lem-lift-cent1}(a), we deduce immediately:
\begin{cor} \label{prop-lift-cent1}
Let $T$ be a Hopf monad on an autonomous category $\cc$ and let~$Q$ be an endofunctor of $T\ti\cc$. Then $Q$ is
centralizable if and only if, for any $T$-module $(M,r)$, the  endofunctor $Q\cp T$ of $\cc$ is centralizable at $M$.

\end{cor}

\subsection{Lifting centralizers}\label{sect-lift1bis}  In this section, given a centralizable
Hopf monad $T$ on an autonomous category $\cc$ and an endofunctor $Q$ of $T\ti\cc$, we show that a centralizer of $Q\cp
T$ lifts uniquely to a centralizer of $Q$. Furthermore, if $Q$ is comonoidal (resp.\@ a bimonad), then it is also a
lift as a monad (resp.\@ a bimonad).

\begin{thm}\label{thm-lift-cent} Let $T$ be a Hopf monad on an autonomous category $\cc$ and let~$Q$ be an endofunctor of $T\ti\cc$.
Assume $Q\cp  T$ is centralizable, with centralizer $(Z_{Q\cp  T},\partial)$. Then:
\begin{enumerate}
\renewcommand{\labelenumi}{{\rm (\alph{enumi})}}
\item The centralizer of $Q\cp  T$ lifts uniquely to a centralizer of $Q$. More precisely,
$Q$ admits  a unique centralizer $(Z_Q,\Tilde{\partial})$ such that  $U_TZ_Q=Z_{Q\cp  T}U_T$ and:
\begin{equation*}
\Tilde{\partial}_{(M,r),(N,s)}=(Q(s)\otimes \id_{Z_{Q\cp  T}(M)})
\partial_{M,N}
\end{equation*}
for all $T$\ti modules $(M,r)$ and $(N,s)$.
\item If $Q$ is comonoidal, the monad $Z_Q$ is a lift of the monad $Z_{Q\cp  T}$ to $T\ti\cc$.
\item If $Q$ is a bimonad, the bimonad $Z_Q$ is a lift of the bimonad $Z_{Q\cp  T}$ to $T\ti\cc$.
\end{enumerate}
\end{thm}

\begin{proof}
Part (a) is a direct consequence of Lemma~\ref{lem-centralisable} and Proposition~\ref{lem-lift-cent1}(b).
Let $(Z_Q,\Tilde{\partial})$ be the centralizer of $Q$ given by Part (a).

Assume $Q$ is comonoidal. Then $Q\cp  T$ is comonoidal  by
Section~\ref{sect-crossprod}. Therefore both $Z_Q$ and
$Z_{Q\cp  T}$ are monads by Theorem~\ref{doub1}(a). Denote $\eta$ and
$\varepsilon$  the unit and counit of the adjunction $(U_T,F_T)$. By
Part (a), we have:
\begin{equation*}
U_T(\Tilde{\partial})=(U_TQ(\varepsilon)\otimes \id_{Z_{Q\cp  T}U_T})
\partial_{U_T,U_T}.
\end{equation*}
By definition of the product $\Tilde{m}$ of $Z_Q$, we have:
\begin{equation*}
(\id_{QF_T \otimes QF_T} \otimes \Tilde{m})\Tilde{\partial}^2_{1_{T\ti\cc}, F_T,F_T} =\bigl(Q_2(F_T,F_T) \otimes
\id_{Z_Q}\bigr)\Tilde{\partial}_{1_{T\ti\cc}, F_T\otimes F_T}.
\end{equation*}
Hence, we get:
\begin{align*}
\bigl(U_TQ(\varepsilon_{F_T}) &\otimes U_TQ(\varepsilon_{F_T}) \otimes U_T(\Tilde{m}) \bigr)\partial^2_{U_T, T,T}\\
&=\bigl(U_T\bigl(Q_2(F_T,F_T) Q(\varepsilon_{F_T\otimes F_T})\bigr ) \otimes \id_{U_TZ_Q}\bigr)\partial_{U_T, T\otimes
T}.
\end{align*}
Composing on the right with $(\id_{U_T} \otimes \eta \otimes\eta)$
and then using the expression of the comonoidal structure of $Q\cp  T$ (see Section~\ref{sect-crossprod}) and the identity $\varepsilon_{F_T} F_T(\eta) =\id_{F_T}$, we obtain:
\begin{equation*}
\bigl(\id_{(Q\cp T)^{\otimes 2}} \otimes U_T(\Tilde{m}) \bigr)
\partial^2_{U_T, 1_\cc,1_\cc} =\bigl((Q\cp T)_2 \otimes \id_{Z_{Q\cp  T}U_T}\bigr)\partial_{U_T, \otimes}
\end{equation*}
and so, by definition of the product $m$ of $Z_{Q\cp T}$, we have: $U_T(\Tilde{m})=m_{U_T}$.
Moreover, denoting $\Tilde{u}$ and $u$ the units of
$Z_Q$ and $Z_{Q\cp T}$ respectively, we have:
\begin{align*}
U_T(\Tilde{u})&= U_T\bigl( (Q_0 \otimes \id_{Z_Q}) \Tilde{\partial}_{1_{T\ti\cc},(\un,T_0)} \bigr ) \\
&= \bigl( U_T(Q_0)U_TQ(T_0) \otimes \id_{U_TZ_Q}\bigr) \partial_{U_T,\un}   \\
&=\bigl((Q\cp T)_0 \otimes \id_{Z_{Q\cp  T}U_T}\bigr) \partial_{U_T,\un}=u_{U_T}.
\end{align*}
Hence Part (b).

Suppose now $Q$ is a bimonad. Then $Q\cp  T$ is a bimonad (see Section~\ref{sect-crossprod}).
Therefore both $Z_Q$ and $Z_{Q\cp  T}$ are bimonads by
Theorem~\ref{doub1}(c). By definition of the morphism $(Z_Q)_2$, we
have:
\begin{equation*}
\bigl(\id_{QF_T} \otimes (Z_Q)_2\bigr)\Tilde{\partial}_{\otimes , F_T}=(q_{F_T} \otimes \id_{Z_Q^{\otimes
2}})(\Tilde{\partial}_{1_{T\ti\cc},QF_T} \otimes \id_{Z_Q})(\id_{1_{T\ti\cc}} \otimes
\Tilde{\partial}_{1_{T\ti\cc},F_T}),
\end{equation*}
where $q$ is the product of $Q$. Thus:
\begin{align*}
\bigl(U_TQ(\varepsilon_{F_T}) \otimes U_T&((Z_Q)_2)\bigr)\partial_{U_T\otimes U_T, T}=\bigl(U_T(q_{F_T}Q(
\varepsilon_{QF_T})) \otimes \id_{(U_TZ_Q)^{\otimes 2}}\bigr)\\ & \circ \bigl(\partial_{U_T,Q\cp T} \otimes
\id_{U_TZ_Q}\bigr)\bigl(\id_{U_T} \otimes (U_TQ(\varepsilon_{F_T}) \otimes \id_{U_TZ_Q})
\partial_{U_T,T}\bigr).
\end{align*}
Composing on the right with $(\id_{U_T} \otimes \id_{U_T} \otimes\eta)$, since the product of $Q\cp T$ is given by
$p=U_T\bigl(q_{F_T}Q( \varepsilon_{QF_T})\bigr)$, we obtain:
\begin{align*}
\bigl(\id_{Q\cp T}& \otimes U_T((Z_Q)_2)\bigr)\partial_{U_T\otimes U_T, 1_\cc}\\ &=\bigl(p \otimes \id_{(Z_{Q\cp
T}U_T)^{\otimes 2}}\bigr)(\partial_{U_T,Q\cp T} \otimes \id_{Z_{Q\cp  T}U_T})(\id_{U_T} \otimes
\partial_{U_T,1_\cc}),
\end{align*}
and so, by definition of the morphism $(Z_{Q\cp T})_2$, we obtain: $(U_TZ_Q)_2=(Z_{Q\cp T}U_T)_2$. Now, by definition of the morphism $(Z_{Q})_0$, we have:
\begin{equation*}
\bigl(\id_{QF_T} \otimes (Z_Q)_0\bigr)\Tilde{\partial}_{(\un,T_0),F_T}=v_{F_T}.
\end{equation*}
where $v$ is the unit of $Q$. Applying $U_T$ and composing with $\eta$, we get:
\begin{equation*}
\bigl(\id_{Q\cp T} \otimes U_T((Z_Q)_0)\bigr)\partial_{\un,1_\cc}=U_T(v_{F_T})\eta.
\end{equation*}
Since $U_T(v_{F_T})\eta$ is the unit of $Q\cp T$ and by definition of the morphism $(Z_{Q\cp T})_0$, we have: $(U_TZ_Q)_0=(Z_{Q\cp T}U_T)_0$.
Hence $U_TZ_Q=Z_{Q\cp T}U_T$ as comonoidal functors, and Part (c).
\end{proof}

\subsection{The canonical distributive law and the double}\label{sect-candistlawQ}
Throughout this section,  let $T$ be a Hopf monad on an autonomous category $\cc$ and $Q$ be a comonoidal endofunctor of
$T\ti\cc$, such that $Q\cp  T$ is centralizable with centralizer $(Z_{Q\cp  T},\partial)$.

By Theorem~\ref{thm-lift-cent},
the centralizer $(Z_{Q\cp  T},\partial)$ lifts to a centralizer $(Z_Q,\Tilde{\partial})$ of $Q$ and the monad $Z_Q$ is
a lift of the monad $Z_{Q\cp  T}$ to $T\ti\cc$. The monad:
\begin{equation*}
D_{Q,T}=Z_Q \cp  T,
\end{equation*}
is called the \emph{double of the pair $(Q,T)$}. Since lifts correspond bijectively with distributive laws (see
Theorem~\ref{propdistrib}), there exists a unique distributive law $\Omega$ of $T$ over $Z_{Q\cp T}$ such that:
\begin{equation*}
Z_Q=\Tilde{Z}_{Q\cp T}^\Omega.
\end{equation*}
This distributive law is called the \emph{canonical distributive law of the pair $(Q,T)$}. It provides a description of structure of the monad $D_{Q,T}$:
\begin{equation*}
D_{Q,T}=Z_{Q\cp  T} \circ_{\Omega} T.
\end{equation*}

\begin{prop}\label{prop-canlawQT}
\begin{enumerate}
\renewcommand{\labelenumi}{{\rm (\alph{enumi})}}
\item If $Q$ is a bimonad,  then the canonical distributive law~$\Omega$ is comonoidal, $D_{Q,T}$ is a bimonad,
and $Z_Q=\Tilde{Z}_{Q\cp T}^\Omega$ as bimonads.
\item If $Q$ is a Hopf monad, then $D_{Q,T}$ is a Hopf monad.
\end{enumerate}
\end{prop}
\begin{proof}
Let us prove Part (a). By Theorem~\ref{thm-lift-cent}, $Z_Q$ is a lift of $Z_{Q\cp  T}$ as a bimonad. Therefore, by Theorem~\ref{propdistribbimon}, $\Omega$ is comonoidal and $D_{Q,T}$ is a bimonad.

Let us prove Part (b). Since $Q\cp T$ is a Hopf monad (see
Section~\ref{sect-crossprod}),  so is $Z_{Q\cp  T}$ (by
Theorem~\ref{doub1}). Therefore $D_{Q,T}$ is a Hopf monad (by
Corollary~\ref{cordistribHopf}).
\end{proof}

Let $\uu\co \zz_Q(T\ti\cc) \to \cc$ be the functor defined as the composition of the forgetful functors $\uu_Q \co \zz_Q(T\ti\cc) \to T\ti\cc$ and $U_T\co T\ti\cc \to \cc$, that is:
\begin{equation*}
\uu\bigl((M,r),\sigma\bigr)=M \quad \text{and} \quad \uu(f)=f.
\end{equation*}
Denoting $\eta$ and $u$ the units of $T$ and $Z_{Q\cp  T}$,  let $I\co
D_{Q,T} \ti \cc \to \zz_Q(T\ti\cc) $ be the functor defined by:
\begin{equation*}
I(M,r)=\bigl((M,ru_{T(M)}),\sigma\bigr) \quad \text{and} \quad I(f)=f,
\end{equation*}
where $\sigma_{(N,s)}=(U_TQ(s) \otimes rZ_{Q \cp
T}(\eta_M))\partial_{M,N}$.
\begin{thm}\label{thm-doublepairQT}
The functor $I$ is an isomorphism of categories such that the
following triangle commutes:
\begin{equation*}
\xymatrix @!0 @C=3pc @R=2.5pc {D_{Q,T} \ti \cc \ar[rr]^(.47)I \ar[dr]_{U_{D_{Q,T}}}& \ar@{}[d]|(.4){\circlearrowright} & \zz_Q(T\ti\cc) \ar[dl]^{\uu}\\
& \cc & }
\end{equation*}
Furthermore, if $Q$ is a bimonad (so that $D_{Q,T}$ is a bimonad and $\zz_Q(T\ti\cc)$ is monoidal), then the functor
$I$ is strict monoidal (and so $\uu I=U_{D_{Q,T}}$ as monoidal functors).
\end{thm}
\begin{rem}
The functor $\ff=IF_{D_{Q,T}}\co
\cc \to \zz_Q(T\ti\cc)$ is left adjoint to $\uu$ and the adjunction
$(\ff,\uu)$ is monadic with monad $D_{Q,T}$ (see Remark~\ref{rem-def-monadic}). If $Q$ is a bimonad, this adjunction is monoidal and $D_{Q,T}$ is its associated bimonad (see Theorem~\ref{thm-gen-ex}).
\end{rem}

\begin{proof}
By Section~\ref{sect-dist-law-monads}, since $D_{Q,T}=Z_{Q\cp  T} \circ_{\Omega} T$ and
$Z_Q=\Tilde{Z}_{Q\cp T}^\Omega$, the functor:
\begin{equation*}
L\co\left\{\begin{array}{ccl}
D_{Q,T}\ti\cc & \longrightarrow & Z_Q\ti(T\ti\cc) \\
(M,r)& \longmapsto & \bigl ((M,ru_{T(M)}), r Z_{Q\cp  T}(\eta_M) \bigr)
\end{array}\right.
\end{equation*}
is an isomorphism of categories. By Theorem~\ref{mocentermon1}, the functor:
\begin{equation*}
E\co\left\{\begin{array}{ccl}
Z_Q\ti(T\ti\cc) & \longrightarrow & Z_Q(T\ti\cc) \\
\bigl((M,r),s\bigr)& \longmapsto & \bigl((M,r),(\id_Q \otimes s) \Tilde{\partial}_{(M,r),1_{T\ti\cc}}\bigr)
\end{array}\right.
\end{equation*}
is an isomorphism of categories. Using Theorem~\ref{thm-lift-cent}(a), one verifies that $I=EL$. Thus $I$ is an
isomorphism of categories, and it clearly satisfies $\uu I=U_{D_{Q,T}}$.

Assume $Q$ is a bimonad. Then $L$ is strict monoidal (by Theorem~\ref{propdistribbimon}) and $E$ is strict monoidal (by
Theorem~\ref{mocentermon1}). Hence $I=EL$ is strict monoidal, and so $\uu I=U_{D_{Q,T}}$ as monoidal functors (since
$\uu$ and $U_{D_{Q,T}}$ are strict monoidal).
\end{proof}

The canonical distributive law $\Omega$ can be described explicitly as follows. By Proposition~\ref{propcentral1}, we
have:
\begin{equation*}
Z_{Q\cp  T}(X)=\int^{Y \in \cc} \hspace*{-1em} \ldual{Q}\cp  T(Y) \otimes X \otimes Y,
\end{equation*}
with universal dinatural transformation:
\begin{equation*}
i_{X,Y}=(\lev_{Q\cp  T(Y)} \otimes \id_{Z_{Q\cp T}(X)})(\id_{\ldual{Q}\cp  T(Y)} \otimes
\partial_{X,Y}).
\end{equation*}
Recall that $T(i)$ is a universal dinatural trans\-formation (see
Proposition~\ref{HM-creat-colim}). Denote $s^l$ the left antipode of $T$ and $\varepsilon$ the counit of the adjunction $(F_T,U_T)$.

\begin{prop}\label{prop-canlawQT-inv}
The canonical distributive law $\Omega$ of the pair $(Q,T)$ is invertible, and $\Omega$ and $\Omega^{-1}$
are characterized as natural transformations by:
\begin{align*}
&\Omega_X T(i_{X,Y})
= i_{T(X),T(Y)} \bigl(\ldual{b}_Y s^l_Y T(\ldual{a}_Y)  \otimes \id_{T(X) \otimes T(Y)}\bigr)  T_3(\ldual{Q \cp  T}(Y),X,Y),\\
&\Omega^{-1}_X i_{T(X),Y}
=\bigl(\lev_{Q\cp T(Y)} \otimes  T(i_{X,T(Y)}) \otimes \lev_Y \bigr)(\id_{\ldual{Q}\cp T (Y)} \otimes E_{X,Y}\otimes \id_{Y}),
\end{align*}
where $a_Y=U_T(\varepsilon_{QF_T(Y)})$, $b_Y=U_TQ(\varepsilon_{F_T(Y)})$, and:
\begin{align*}
E_{X,Y}=&(a_Y T(b_Y) \otimes \id_{T(\ldual{Q}\cp T\,T(Y) \otimes X \otimes T(Y))} \otimes s^l_Y) \\
&\circ T_3(Q\cp  T\,T(Y),\ldual{Q}\cp T\,T(Y) \otimes X \otimes T(Y),\ldual{T}(Y))\\
&\circ T(\lcoev_{Q\cp T\,T(Y)}\otimes \id_X \otimes \lcoev_{T(Y)}).
\end{align*}
\end{prop}

\begin{rem}\label{rem-samelaw}
In the special case $Q=1_{T\ti\cc}$, we have: $1_{T\ti\cc}\cp  T=T$
and so, by Proposition~\ref{prop-canlawQT-inv}, the canonical
distributive law of the pair $(1_{T\ti\cc},T)$ is nothing but the
canonical law of $T$ defined in Section~\ref{sect-conodistlaw},
and the double $D_{1_{T\ti\cc},T}$ of the pair
$(1_{T\ti\cc},T)$ coincides with the double $D_T$ of $T$ defined in
Section~\ref{sect-doubleHM}.
\end{rem}

\begin{proof}
Note that $a_Y$ and $a'_Y=s^l_{Q\cp T(Y)} T(\ldual{a_Y})$ are the $T$\ti actions of $QF_T(Y)$ and $\ldual{Q}F_T(Y)$ respectively. By adjunction we have: $b_Y Q\cp T(\eta_Y)=\id_{Q\cp T(Y)}$.

Recall that $\Tilde{Z}_{Q\cp T}^\Omega$ is the centralizer of $Q$,
with universal dinatural transformation:
\begin{equation*}
j_{(M,r),(N,s)}=i_{M,N}\bigl(\ldual{U}_TQ(\varepsilon_{(N,s)}) \otimes \id_M \otimes \id_N\bigr).
\end{equation*}
In particular, given two objects $X,Y$ of $\cc$, the morphism $j_{F_T(X),F_T(Y)}$ is $T$\ti linear, that is,
\begin{equation*}
Z_{Q\cp T}(\mu_X)\Omega_{T(X)} T(j_{F_T(X),F_T(Y)})=j_{F_T(X),F_T(Y)}\gamma_{X,Y},
\end{equation*}
where:
\begin{equation*}
\gamma_{X,Y}=\bigl(a'_Y \otimes \mu_X \otimes \mu_Y \bigr)T_3(\ldual{Q \cp  T}(Y),T(X),T(Y))
\end{equation*}
is the $T$\ti action of $\ldual{Q}F_T(Y) \otimes F_T(X) \otimes F_T(Y)$. Composing on the right with $T(\id_{\ldual{Q \cp  T}(Y)} \otimes \eta_X \otimes \eta_Y)$, the left-hand side becomes:
\begin{align*}
Z_{Q\cp T}&(\mu_X)\Omega_{T(X)} T\bigl(i_{T(X),T(Y)} (\ldual{b}_Y \otimes \eta_X \otimes \eta_Y)\bigr)\\
&=Z_{Q\cp T}(\mu_X)\Omega_{T(X)}TZ_{Q\cp T}(\eta_X)T(i_{X,Y})=\Omega_X T(i_{X,Y}),
\end{align*}
and the right-hand side becomes:
\begin{align*}
i_{T(X),T(Y)} (\ldual{b}_Y a'_Y\otimes \mu_X T(\eta_X)\otimes \mu_Y T(\eta_Y))T_3(\ldual{Q}\cp T(Y),X,Y).
\end{align*}
Hence the formula for $\Omega$.

Let $\Omega'\co Z_{Q\cp  T} T \to TZ_{Q\cp  T} $ be the natural
transformation defined by:
\begin{equation*}
\Omega'_X  i_{T(X),Y} =\bigl(\lev_{Q\cp T(Y)} \otimes  T(i_{X,T(Y)}) \otimes \lev_Y \bigr)(\id_{\ldual{Q}\cp T (Y)} \otimes E_{X,Y}\otimes \id_{Y})
\end{equation*}
Using the axioms of a left antipode, one shows that $\Omega'\Omega=\id_{TZ_{Q\cp  T}}$ and
$\Omega'\Omega=\id_{Z_{Q\cp  T}T}$ by verifying that
$\Omega'_X\Omega_XT(i_{X,Y})=T(i_{X,Y})$ and
$\Omega_X\Omega'_Xi_{T(X),Y}=i_{T(X),Y}$. This is left to the reader.
Note that when $Q$ is a Hopf monad, the invertibility of $\Omega$ follows from Proposition~\ref{propdistinv}, since in this case both $T$ and  $Z_{Q\cp  T}$ are Hopf monads.
\end{proof}

\begin{rem}
Let $T$ be a Hopf monad on an autonomous  category $\cc$ and $Q$ be
a comonoidal endofunctor of $T\ti\cc$ such that $Q\cp  T$ is
centralizable. Consider the following diagram:
\begin{equation*}
\xymatrix @!0 @C=4.5pc @R=3.5pc {& \ar[dl]_{\mathcal{W}} \ar@/^.3pc/[dr] Z_Q\ti(T\ti\cc) & \\
Z_{Q\cp  T}\ti\cc \ar@/^.3pc/[dr] & & \ar@/^.3pc/[ul] \ar@/^.3pc/[dl] T\ti\cc \\
& \ar@/^.3pc/[ul] \ar@/^.3pc/[ur] \cc  &}
\end{equation*}
where a double arrow represents the adjunction of the corresponding
monad and the functor $\mathcal{W}$ is defined by
$\mathcal{W}\bigl((M,r),s\bigr)=(M,s)$ and $\mathcal{W}(f)=f$. This
diagram is a distributive adjoint square in the sense of
Beck~\cite{Beck1} whose distributive law is precisely the canonical
distributive law $\Omega$ of the pair $(Q,T)$. Furthermore,
since~$\Omega$ is invertible by Proposition~\ref{prop-canlawQT-inv},
the monad $T$ lifts to a monad $\Tilde{T}^{\Omega^{-1}}$ on
$Z_{Q\cp T}\ti\cc$.  Therefore, since $Z_Q=\Tilde{Z}_{Q\cp T}^\Omega$, we have an isomorphism of categories:
\begin{equation*}
\Tilde{T}^{\Omega^{-1}}\ti(Z_{Q\cp T}\ti\cc) \simeq (T\circ_{\Omega^{-1}}Z_{Q\cp T})\ti\cc \simeq (Z_{Q\cp T}\circ_\Omega T)\ti\cc \simeq
Z_Q\ti(T\ti\cc).
\end{equation*}
Via this isomorphism, $\mathcal{W}$ is the forgetful functor
$U_{\Tilde{T}^{\Omega^{-1}}}$. Hence $\mathcal{W}$ is monadic.
Note that when $Q$ is a bimonad, the four monadic adjunctions are
monoidal.
\end{rem}

\subsection{Proof of
Theorem~\ref{thm-candistlaw}}\label{sect-proof-thm-candistlaw} This
is a direct consequence of Remark~\ref{rem-samelaw} and
Propositions~\ref{prop-canlawQT} and~\ref{prop-canlawQT-inv} applied
to the Hopf monad $Q=1_{T\ti\cc}$.

\subsection{Proof of Theorems~\ref{thm-Rmat-double} and \ref{thm-doublable}}\label{sect-proof-thm-doublable}
By Theorem~\ref{thm-doublepairQT} applied to the Hopf monad
$Q=1_{T\ti\cc}$ and Remark~\ref{rem-samelaw}, the functor $I\co D_T
\ti \cc \to \zz(T\ti\cc)$ of Theorem~\ref{thm-doublable} is a strict
monoidal isomorphism of monoidal categories such that $\uu
I=U_{D_T}$.

Now, by
Remark~\ref{remcenterusual}, the category $\zz(T\ti\cc)$ is a braided category with braiding:
\begin{equation*}
\tau_{\bigl((M,r),\gamma\bigr),\bigl((N,s),\delta\bigr)}=\gamma_{(N,s)}.
\end{equation*}
Therefore, since $I$ is a strict monoidal isomorphism, there exists
a unique braiding~$c$ on $D_T\ti\cc$ such that $I$ is braided. By
Theorem~\ref{thm-quasitrigHopfmon}, $c$ is encoded by an \Rt matrix
$R$ for $D_T$.  Let $p$ and $e=u_T\eta$ be the product and unit of
$D_T$. Then $R$ is given by:
\begin{align*}
R_{X,Y}&=c_{F_{D_T(X)},F_{D_T(Y)}}(e_X \otimes
e_Y)\\
&=\tau_{IF_{D_T(X)},IF_{D_T(Y)}}(e_X \otimes
e_Y)\\
&= \bigl(p_Y u_{D_T(Y)}\otimes p_X
Z_T(\eta_{D_T(X)})\bigr)\partial_{D_T(X),D_T(Y)}(
e_X \otimes e_Y )\\
&= \bigl(p_Y u_{D_T(Y)} T(e_Y) \otimes
p_X Z_T(\eta_{D_T(X)}e_X)\bigr)\partial_{X,Y}\\
&= \bigl(p_Y D_T(e_Y) u_{T(Y)} \otimes
p_X D_T(e_X)Z_T(\eta_X)\bigr)\partial_{X,Y}\\
&= \bigl( u_{T(Y)} \otimes Z_T(\eta_X)\bigr)\partial_{X,Y}.
\end{align*}
This concludes the proof of Theorems~\ref{thm-Rmat-double} and
\ref{thm-doublable}.

\subsection{Proof of
Theorem~\ref{thm-centraZ}}\label{sect-proof-centraZ} This is a
direct consequence of Remark~\ref{rem-samelaw} and
Theorem~\ref{thm-lift-cent} applied to the Hopf monad
$Q=1_{T\ti\cc}$.

\section{The double of a Hopf algebra in a braided category}\label{sect-cas-hopf-ds-braided}

 In this section, we extend several classical notions concerning
a Hopf algebra over a field to a Hopf algebra $A$ in a braided
autonomous category $\bb$, namely: quasitriangularity and \Rt
matrices and the double $D(A)$ of $A$. Our approach consists
in applying the results of previous sections  to the Hopf monad $?
\otimes A$.

We need to assume that $\bb$ admits a coend $C$.
Then the Hopf monad $? \otimes A$ is centralizable and its
centralizer is of the form $? \otimes Z(A)$,  where $Z(A)$ is  a
certain Hopf algebra in $\bb$ called the centralizer of $A$. As an
object of $\bb$, $Z(A)=\ldual{A} \otimes C$.

We then define the \emph{double of $A$} as $D(A)=A \otimes_\Omega
Z(A)=A \otimes \ldual{A} \otimes C$, where $\Omega$  is an explicit
distributive law.  The double $D(A)$ is a quasitriangular Hopf
algebra in $\bb$  such that $D_{?\otimes A}=? \otimes D(A)$ (as
quasitriangular Hopf monads).  It satisfies:
$\zz(\rmod{\bb}{A})\simeq \rmod{\bb}{{D(A)}}$ (as braided
categories).   When $\bb=\vect_\kk$,
we have:  $C=\kk$, $A$ is a
finite-dimensional Hopf algebra over $\kk$, $Z(A)=(A^*)^\cop$, and $D(A)$ is the usual Drinfeld double of $A$.

\subsection{Hopf monads represented by Hopf algebras}\label{sect-repres}
Let $\bb$ be a braided autonomous category and let $A$ be a Hopf algebra in $\bb$.
A Hopf monad $T$ on a $\bb$ is said to
be \emph{represented on the left} (resp.\@ \emph{on the right}) by $A$ if it is isomorphic to
the Hopf monad $A\otimes ?$ (resp.\@ $? \otimes A$) defined in Example~\ref{exa-HMfromHA}.

More generally, let $T$ be a Hopf monad on an autonomous category $\cc$. If $(A,\sigma)$ is a Hopf algebra in the center $\zz(\cc)$ of $\cc$, then the Hopf monad $T$ is said to be \emph{represented on the left} by $(A,\sigma)$ if it is isomorphic to the Hopf monad $A \otimes_\sigma ?$ on~$\cc$ defined in Example~\ref{exa-HMfromHAcent}. Likewise, if $(A,\sigma)$ is a Hopf algebra in $\Bar{\zz}(\cc)$, then the Hopf monad $T$ is said to be \emph{represented on the right} by $(A,\sigma)$ if it is isomorphic to the Hopf monad $? \otimes_\sigma A$ on $\cc$.

Not all Hopf monads can be so represented by Hopf algebras (see Remark~\ref{rem-ZTnotrepresent} for an example).

\subsection{Coends as Hopf algebras}\label{sect-CTfromHM}
Let $T$ be an endofunctor of an autonomous category $\cc$. If $\cc$ admits a braiding $\tau$,
then, by Proposition~\ref{propcentral1}, $T$ is centralizable if and only if the coend:
\begin{equation*}
C_T=\int^{Y \in \bb} \ldual{T(Y)} \otimes Y
\end{equation*}
exists.
Assume this is the case. By Lemma~\ref{lem-coend-adj}, if $T$ is a monad, then $C_T$ coincides with the coend $\int^{(M,r) \in T\ti\cc} \ldual{U}_T(M,r) \otimes U_T(M,r)$ of $U_T$. According to Majid~\cite{Maj2}, the (co)end of a strong monoidal functor from an autonomous category to a braided category is a Hopf algebra. In particular, if $T$ is a Hopf monad and $\tau$ a braiding on $\cc$, then $C_T$ is a Hopf algebra in $\cc$ braided by $\tau$. In this section we recover this structure explicitly in terms of
the braiding $\tau$ and the Hopf monad structure of $T$.

Let $T$ be an endofunctor of an autonomous category $\cc$ such that the coend $C_T=\int^{Y \in \bb} \ldual{T(Y)} \otimes Y$ exists.
Denote $i_Y\co \ldual{T(Y)} \otimes Y \to C_T$ the universal dinatural transformation of $C_T$,
and set:
\begin{equation*}
\delta_Y = \psfrag{Y}[Bc][Bc]{\scalebox{.8}{$Y$}} \psfrag{X}[Bc][Bc]{\scalebox{.8}{$X$}}
\psfrag{T}[Bc][Bc]{\scalebox{.8}{$T(Y)$}} \psfrag{K}[cc][cc]{\scalebox{.9}{$i_Y$}}
\psfrag{C}[Bl][Bl]{\scalebox{.8}{$C_T$}} \rsdraw{.45}{.95}{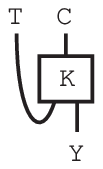} \; =(\id_{T(Y)} \otimes i_Y)(\lcoev_{T(Y)} \otimes
\id_Y)\co X \to T(Y) \otimes C_T, \; \text{depicted} \quad \rsdraw{.45}{.95}{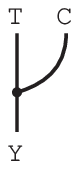}.
\end{equation*}

If $T$ is a monad on $\cc$, then $C_T$ is a coalgebra in $\cc$, with coproduct $\Delta$ and counit~$\varepsilon$ defined by:
\begin{equation*}
\psfrag{Y}[cc][cc]{\scalebox{.8}{$Y$}}  \psfrag{X}[cc][cc]{\scalebox{.8}{$X$}}
\psfrag{T}[Bc][Bc]{\scalebox{.8}{$T(X)$}} \psfrag{G}[Bc][Bc]{\scalebox{.8}{$T(Y)$}}
\psfrag{m}[cc][cc]{\scalebox{.9}{$m$}} \psfrag{C}[Bl][Bl]{\scalebox{.8}{$C_T$}}
\psfrag{a}[cc][cc]{\scalebox{.9}{$T_2(X,Y)$}}
 \psfrag{Z}[Bc][Bc]{\scalebox{.8}{$X\otimes Y$}} \psfrag{D}[cc][cc]{\scalebox{.9}{$\Delta$}}
\rsdraw{.45}{.9}{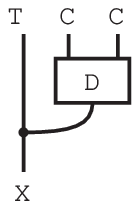} \, = \, \psfrag{u}[cc][cc]{\scalebox{.9}{$\mu_X$}}\rsdraw{.45}{.9}{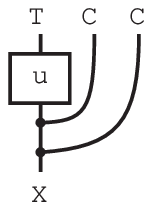} \quad \text{and} \quad \;
\psfrag{X}[cc][cc]{\scalebox{.8}{$X$}}
\psfrag{T}[Bc][Bc]{\scalebox{.8}{$T(X)$}} \psfrag{C}[Bl][Bl]{\scalebox{.8}{$C_T$}}
\psfrag{D}[cc][cc]{$\varepsilon$} \rsdraw{.45}{.9}{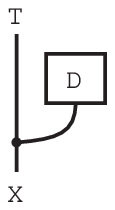} \, = \, \psfrag{D}[cc][cc]{\scalebox{.9}{$\eta_X$}}
\rsdraw{.45}{.9}{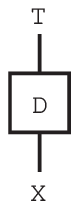}\;,
\end{equation*}
where $\mu$ and $\eta$ are the product and unit of $T$.

If $T$ is comonoidal and $\tau$ is a braiding on $\cc$, then $C_T$ becomes an algebra in $\cc$ with product $m_\tau$ and unit $u$ defined by:
\begin{equation*}
\psfrag{Y}[cc][cc]{\scalebox{.8}{$Y$}}  \psfrag{X}[cc][cc]{\scalebox{.8}{$X$}}
\psfrag{T}[Bc][Bc]{\scalebox{.8}{$T(X)$}} \psfrag{G}[Bc][Bc]{\scalebox{.8}{$T(Y)$}}
\psfrag{m}[cc][cc]{\scalebox{.9}{$m_\tau$}} \psfrag{C}[Bl][Bl]{\scalebox{.8}{$C_T$}} \rsdraw{.45}{.9}{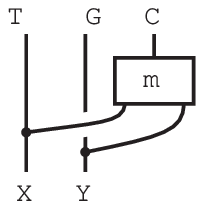} \, = \,
\psfrag{a}[cc][cc]{\scalebox{.9}{$T_2(X,Y)$}}
 \psfrag{Z}[Bc][Bc]{\scalebox{.8}{$X\otimes Y$}} \rsdraw{.45}{.9}{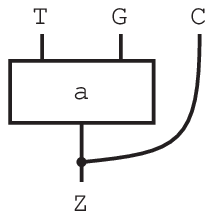}  \quad \text{and} \quad u=\, \psfrag{P}[Bc][Bc]{\scalebox{.9}{$T_0$}} \rsdraw{.35}{.9}{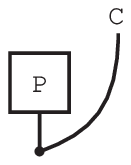}\,, \quad \text{where} \quad \psfrag{X}[Bc][Bc]{\scalebox{.8}{$X$}} \psfrag{Y}[Bc][Bc]{\scalebox{.8}{$Y$}} \tau_{X,Y}=\,\rsdraw{.45}{.9}{braiding}.
\end{equation*}

If $T$ is a bimonad and $\tau$ a braiding on $\cc$, then $(C_T,m^\tau,u,\Delta,\varepsilon)$ is a bialgebra in $\cc$ braided by $\tau$. Furthermore, if $T$ is a Hopf monad, then $C_T$ is a Hopf algebra, whose antipode $S_\tau$ and its inverse $S^{-1}_\tau$ are defined by:
\begin{equation*}
  \psfrag{X}[cc][cc]{\scalebox{.8}{$X$}}
\psfrag{T}[Bc][Bc]{\scalebox{.8}{$T(X)$}}
 \psfrag{C}[Bl][Bl]{\scalebox{.8}{$C_T$}}
\psfrag{D}[cc][cc]{\scalebox{.9}{$S_\tau$}} \rsdraw{.45}{.9}{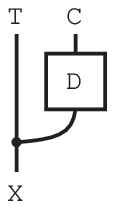} \, = \, \psfrag{D}[cc][cc]{\scalebox{.9}{$s^l_X$}}
\rsdraw{.45}{.9}{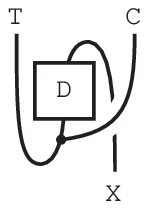} \quad \text{and} \quad
\psfrag{D}[cc][cc]{\scalebox{.9}{$S^{-1}_\tau$}} \rsdraw{.45}{.9}{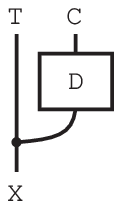}
\, = \, \psfrag{D}[cc][cc]{\scalebox{.9}{$s^r_X$}} \rsdraw{.45}{.9}{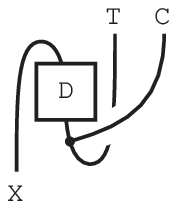} .
\end{equation*}
We denote this Hopf algebra by $C_T^\tau$.

\subsection{The coend of a braided autonomous category}\label{sect-coendbraided} Let $\bb$ be an autonomous category.
The coend:
\begin{equation*}
C=\int^{Y \in \bb} \ldual{Y} \otimes Y,
\end{equation*}
if it exists, is called the \emph{coend of  $\bb$}.

Assume that $\bb$ admits a coend $C$ and denote by
$i_Y\co \ldual{Y} \otimes Y \to C$ its universal dinatural transformation. The
\emph{universal coaction} of $C$ on the objects of $\bb$ is the natural transformation $\delta$ defined by:
\begin{equation*}
\delta_Y=(\id_Y \otimes i_Y)(\lcoev_Y \otimes \id_Y)\co Y \to Y \otimes C, \quad \text{depicted} \quad \psfrag{C}[Bl][Bl]{\scalebox{.8}{$C$}} \psfrag{Y}[Bc][Bc]{\scalebox{.8}{$Y$}}\delta_Y=\rsdraw{.45}{.95}{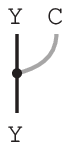}.
\end{equation*}

If $\bb$ is braided, then $C$ is a Hopf algebra in $\bb$. This well-known fact may be viewed as a special case of the construction of Section~\ref{sect-centralrepres}, as  $1_\bb$ is a Hopf monad on $\bb$ and $C=C_{1_\bb}^\tau$ where $\tau$ is the braiding of $\bb$.  Furthermore, the morphism $\omega\co C \otimes C \to \un$, defined by:
\begin{equation*}
\psfrag{Y}[Bc][Bc]{\scalebox{.8}{$Y$}}  \psfrag{X}[Bc][Bc]{\scalebox{.8}{$X$}}
\psfrag{w}[Bc][Bc]{$\omega$}\rsdraw{.45}{.9}{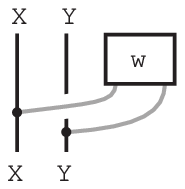} \; = \,
\rsdraw{.45}{.9}{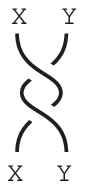}\,,
\end{equation*}
is a Hopf pairing for $C$, that is, it satisfies:
\begin{align*}
&\omega(m \otimes \id_C)=\omega (\id_C \otimes \omega \otimes \id_C)(\id_{C^{\otimes 2}} \otimes \Delta), && \omega(u
\otimes
\id_C)=\varepsilon,\\
&\omega(\id_C \otimes m)=\omega (\id_C \otimes \omega \otimes \id_C)(\Delta \otimes \id_{C^{\otimes 2}}), &&
\omega(\id_C \otimes u)=\varepsilon.
\end{align*}
These axioms imply: $\omega(S \otimes \id_C)=\omega(\id_C \otimes S)$. Moreover the pairing $\omega$ satisfies the self-duality condition:    $\omega \tau_{C,C} (S \otimes
S)=\omega$.

In this section, the structural morphisms of $C$ are drawn in grey and the
Hopf pairing $w\co C \otimes C \to \un$ is depicted as:
\begin{center}
\psfrag{C}[Bc][Bc]{\scalebox{.8}{$C$}} $\omega=$\rsdraw{.45}{.9}{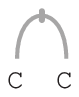}.
\end{center}
\begin{rem}
The category $\bb$ is symmetric if and only if $\omega=\epsilon \otimes \epsilon$. In particular, this is the case when $C=\un$.
\end{rem}

\begin{rem}
The universal coaction of the coend on itself can be expressed in terms of its Hopf algebra structure as follows:
\begin{equation*}
\delta_C=\, \psfrag{C}[Bc][Bc]{\scalebox{.8}{$C$}}\rsdraw{.45}{.95}{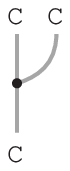}=\;\rsdraw{.45}{.95}{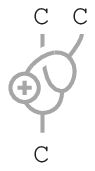}.
\end{equation*}
\end{rem}
\begin{rem}\label{rem-coendmirror}
The coend of the mirror $\bbb$ of $\bb$ is the Hopf algebra $C^\opp$, with self-dual pairing
$\omega(S \otimes \id_C)$.
\end{rem}

\subsection{Centralizers in braided categories}\label{sect-centralrepres}
Let $T$ be an endofunctor of a braided autonomous category $\bb$, with braiding~$\tau$.
Assume that the coend:
\begin{equation*}
C_T=\int^{Y \in \bb} \ldual{T(Y)} \otimes Y
\end{equation*}
exists. Set:
\begin{equation*}
\partial_{X,Y} = \psfrag{Y}[Bc][Bc]{\scalebox{.8}{$Y$}}  \psfrag{X}[Bc][Bc]{\scalebox{.8}{$X$}}
\psfrag{T}[Bc][Bc]{\scalebox{.8}{$T(Y)$}} \psfrag{K}[cc][cc]{$i_Y$} \psfrag{C}[Bl][Bl]{\scalebox{.8}{$C_T$}}
\rsdraw{.45}{.95}{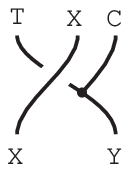}\, = (\tau_{X,T(Y)} \otimes \id_{C_T})(\id_X \otimes \delta_Y) \co X \otimes Y \to T(Y)
\otimes X \otimes C_T.
\end{equation*}
Then $(? \otimes C_T,\partial)$ is a centralizer of $T$.
Likewise, set:
\begin{equation*}
\partial'_{X,Y} = \psfrag{Y}[Bc][Bc]{\scalebox{.8}{$Y$}}  \psfrag{X}[Bc][Bc]{\scalebox{.8}{$X$}}
\psfrag{T}[Bc][Bc]{\scalebox{.8}{$T(Y)$}} \psfrag{K}[cc][cc]{$i_Y$} \psfrag{C}[Bl][Bl]{\scalebox{.8}{$C_T$}}
\rsdraw{.45}{.95}{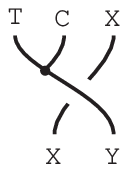}\, = (\delta_Y\otimes \id_X)\tau_{Y,X}^{-1} \co X \otimes Y \to T(Y)
\otimes C_T \otimes X,
\end{equation*}
Then $(C_T\otimes ?,\partial')$ is also a centralizer of $T$.

Assume furthermore that $T$ is a Hopf monad. By Section~\ref{sect-CTfromHM}, the object $C_T$ is endowed with two Hopf algebra structures in $\bb$, namely: $C_T^\tau$ and $(C_T^{\overline{\tau}})^\opp$, where~$\overline{\tau}$ is the mirror of $\tau$.  One verifies that the Hopf monad structure on $? \otimes C_T$ (resp.\@ $C_T \otimes ?$) given by Theorem~\ref{doub1} coincides with that induced by the Hopf algebra $C_T^\tau$  (resp.\@ $(C_T^{\overline{\tau}})^\opp$). Thus:
\begin{thm}\label{thm-cent-br}
Let $T$ be a Hopf monad on a braided autonomous category $\bb$, with braiding $\tau$. Then
$T$ is centralizable if and only if the coend:
\begin{equation*}
C_T=\int^{Y \in \bb} \ldual{T(Y)} \otimes Y
\end{equation*}
exists. If such is the case,
the centralizer of $T$ is
represented on the right by the Hopf algebra $C_T^\tau$ and on the left by the Hopf algebra $(C_T^{\overline{\tau}})^\opp$.
\end{thm}

\begin{rem}\label{rem-ZTnotrepresent}
In general, the centralizer $Z_T$ of a Hopf monad $T$ on an autonomous category $\cc$
is isomorphic neither to $Z_T(\un) \otimes ?$ nor to $? \otimes Z_T(\un)$ as an endofunctor of~$\cc$, see Remark~\ref{rem-Znotrepres} for a counter-example. In particular, it cannot be represented on the left by a Hopf algebra of $\zz(\cc)$, nor on the right by a Hopf algebra of $\Bar{\zz}(\cc)$, in the sense of Section~\ref{sect-repres}.
\end{rem}

\subsection{Centralizers of Hopf algebras}\label{sect-centralizersHA}
Let $\bb$ be a braided autonomous category, with braiding $\tau$, admitting a coend $C=\int^{Y \in \bb} \ldual{Y} \otimes Y$.
 Recall that $C$ is a Hopf algebra in
$\bb$ endowed with a Hopf pairing, and denote by $\delta$ the universal coaction of $C$ (see Section~\ref{sect-coendbraided}).

Let $A$ be a Hopf algebra in $\bb$. Set $Z(A)=\ldual{A}
\otimes C$. Then:
\begin{equation*}
Z(A)=\int^{Y \in \bb} \ldual{(Y \otimes A)} \otimes
Y,
\end{equation*}
with universal dinatural transformation given by:
\begin{equation*}
i_Y=  \psfrag{Y}[Bc][Bc]{\scalebox{.8}{$Y$}}  \psfrag{X}[Bc][Bc]{\scalebox{.8}{$\ldual{Y}$}}
\psfrag{A}[Bc][Bc]{\scalebox{.8}{$\ldual{A}$}} \psfrag{C}[Bc][Bc]{\scalebox{.8}{$C$}} \rsdraw{.45}{.9}{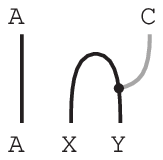}
\; \co \ldual{(Y \otimes A)} \otimes
Y=\ldual{A} \otimes  \ldual{Y} \otimes Y \to\ldual{A}
\otimes C.
\end{equation*}
We endow the object $Z(A)$ with the Hopf algebra structure $C_{? \otimes A}^\tau$ defined in Section~\ref{sect-CTfromHM}. Explicitly, the structural morphisms of $Z(A)$  are:
\begin{center}
\psfrag{B}[Br][Br]{\scalebox{.8}{$\ldual{A}$}} \psfrag{C}[Bc][Bc]{\scalebox{.8}{$C$}}
$m_{Z(A)}=$\;\rsdraw{.4}{1}{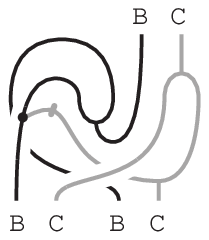}\,,\quad
$\Delta_{Z(A)}=$\;\rsdraw{.4}{1}{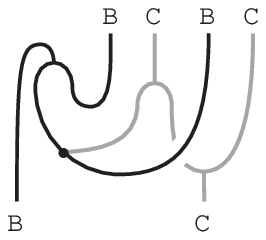}, \quad $S_{Z(A)}=$\;\rsdraw{.4}{1}{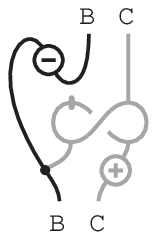}\,, \\[.8em]
$u_{Z(A)}=$\;\rsdraw{.4}{1}{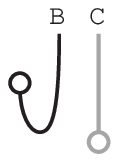},\qquad
$\varepsilon_{Z(A)}=$\;\rsdraw{.4}{1}{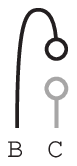}\,.
\end{center}
\begin{exa}
Let $H$ be a finite-dimensional Hopf algebra over a field $\kk$. Note that $\kk$ is the coend of the category $\vect_\kk$ of finite-dimensional vector spaces. Then the centralizer of $H$ is $Z(H)=(H^*)^\mathrm{cop}$.
\end{exa}

By Theorem~\ref{thm-cent-br}, $(? \otimes Z(A),\partial)$ is a centralizer of $? \otimes A$, where:
\begin{equation*}
\partial_{X,Y}=
   \psfrag{Y}[Bc][Bc]{\scalebox{.8}{$Y$}}  \psfrag{X}[Bc][Bc]{\scalebox{.8}{$X$}} \psfrag{A}[Bc][Bc]{\scalebox{.8}{$A$}}
\psfrag{B}[Bc][Bc]{\scalebox{.8}{$\ldual{A}$}} \psfrag{C}[Bc][Bc]{\scalebox{.8}{$C$}} \rsdraw{.45}{.9}{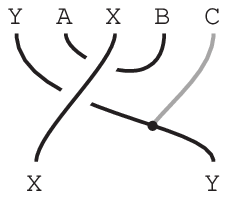},
\end{equation*}
and the Hopf monad structure on $? \otimes Z(A)$ given by Theorem~\ref{doub1} is that induced by the Hopf algebra $Z(A)$.
Hence, by Theorem~\ref{mocentermon1}, we
have:
\begin{equation*}
  \zz_{? \otimes A}(\bb)\simeq(? \otimes Z(A))\ti\bb=\bb_{Z(A)}
\end{equation*}
as monoidal categories.

\begin{rem}
One can show that the Hopf algebra $Z(A)$ represents on the left the centralizer $Z_{A\otimes ?}$ of $A\otimes ?$, and so:
\begin{equation*}
\zz_{A \otimes ?}(\bb)\simeq(Z(A) \otimes ?)\ti\bb=\lmod{\bb}{Z(A)}
\end{equation*}
as monoidal categories.
\end{rem}

\subsection{\Rt matrices for Hopf algebras in braided categories}\label{sect-RmatrixA} In~\cite{Drin}, Drinfeld introduced the notion of \Rt matrix for a Hopf algebra $H$ over a field $\kk$. When $H$ is finite-dimensional, \Rt matrices for $H$ are in bijection with braidings on the category of finite-dimensional $H$-modules. The aim of this section is to extend the notion of an \Rt matrix to a Hopf algebra $A$ in braided autonomous category so as to preserve this bijective correspondence. Note that the definition of an \Rt matrix for~$A$ as a morphism $\mathfrak{r}\co \un \to A \otimes A$ by straightforward extension of Drinfeld's axioms (sometimes found in the literature) does not fulfil this objective. Recall that braidings on the autonomous category $\rmod{\bb}{A}=(? \otimes A)\ti\bb$
are encoded by \Rt matrices for the Hopf monad $? \otimes A$. When $\bb$ admits a coend, we can encode \Rt matrices for $? \otimes A$ in terms of $A$ by introducing \Rt matrices for $A$.

Let $A$ be a Hopf algebra in a braided autonomous category $\bb$, with braiding~$\tau$. Assume that $\bb$ admits a coend $C$.  Any \Rt matrix $R_{X,Y} \co X \otimes Y \to Y \otimes A \otimes X \otimes A$ for the Hopf monad $? \otimes A$ gives rise to a unique morphism $\mathfrak{r} \co C \otimes C \to A \otimes A$ in $\bb$, defined by:
\begin{equation*}
\psfrag{X}[Bc][Bc]{\scalebox{.8}{$X$}} \psfrag{Y}[Bc][Bc]{\scalebox{.8}{$Y$}}
\psfrag{M}[Bc][Bc]{\scalebox{.8}{$\ldual{X}$}} \psfrag{W}[Bc][Bc]{\scalebox{.8}{$\ldual{Y}$}}
\psfrag{C}[Bc][Bc]{\scalebox{.8}{$C$}} \psfrag{A}[Bc][Bc]{\scalebox{.8}{$A$}} \psfrag{R}[Bc][Bc]{\scalebox{.8}{$R_{X,Y}$}}
\psfrag{r}[Bc][Bc]{\scalebox{.9}{$\mathfrak{r}$}} \rsdraw{.45}{1}{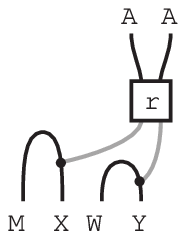}
=\;\rsdraw{.45}{1}{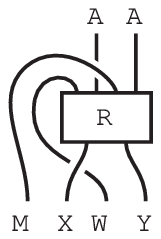}\;, \quad \text{so that:} \quad
R_{X,Y}= \rsdraw{.45}{1}{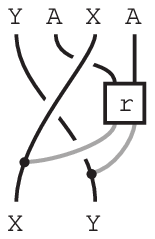}\,.
\end{equation*}
Re-writing the axioms for $R_{X,Y}$ (see Section~\ref{sect-QTHM}) in terms of $\mathfrak{r}$ leads to the following definition:
a \emph{\Rt matrix} for $A$ is a morphism
\begin{equation*}
\mathfrak{r}\co C \otimes C \to A \otimes A
\end{equation*}
in $\bb$, which satisfies:
\begin{center}
\psfrag{C}[Bc][Bc]{\scalebox{.8}{$C$}} \psfrag{A}[Bc][Bc]{\scalebox{.8}{$A$}}
\psfrag{r}[Bc][Bc]{\scalebox{.9}{$\mathfrak{r}$}}
\rsdraw{.45}{1}{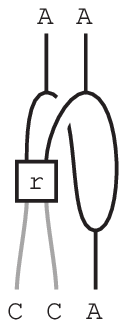} \, $=$ \, \rsdraw{.45}{1}{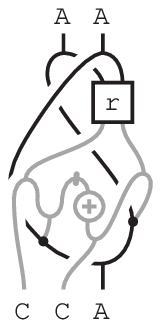}\;, \qquad
\rsdraw{.45}{1}{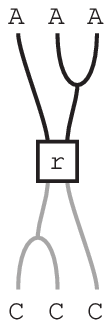} $=$ \rsdraw{.45}{1}{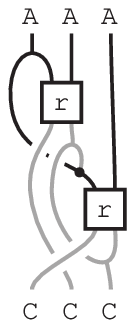}\;, \qquad
\rsdraw{.45}{1}{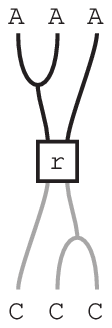} $=$ \rsdraw{.45}{1}{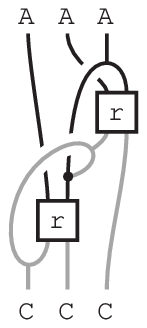}\,,\\[.8em]
\rsdraw{.45}{1}{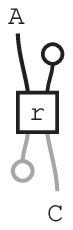} \, $=$ \, \rsdraw{.45}{1}{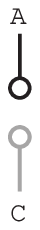} \, $=$ \, \rsdraw{.45}{1}{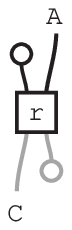} \,.
\end{center}

\begin{rem}
For finite-dimensional Hopf algebras over a field $\kk$, our definition of an \Rt matrix coincides with Drinfeld's definition (as the coend of $\vect_\kk$ is $\kk$).
\end{rem}

An \Rt matrix $\mathfrak{r}$ for $A$ defines an \Rt matrix for $?\otimes A$ (by definition) and so a braiding~$c$ on $\rmod{\bb}{A}=(? \otimes A)\ti \bb$  (by Theorem~\ref{thm-quasitrigHopfmon})  as:
\begin{equation*}
c_{(M,r),(N,s)}=(s\otimes r)R_{M,N}=\psfrag{X}[Bc][Bc]{\scalebox{.8}{$M$}} \psfrag{Y}[Bc][Bc]{\scalebox{.8}{$N$}}
\psfrag{n}[Bc][Bc]{\scalebox{.9}{$r$}}\psfrag{s}[Bc][Bc]{\scalebox{.9}{$s$}} \psfrag{r}[Bc][Bc]{\scalebox{.9}{$\mathfrak{r}$}}
\rsdraw{.45}{1}{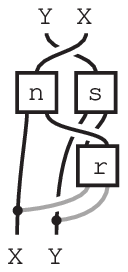}\;.
\end{equation*}

As braidings on $\lmod{\bb}{A}$ are in bijective correspondence with braidings on  $\rmod{\bb}{A}$ (see Remark~\ref{rem-braidleftright}), an \Rt matrix $\mathfrak{r}$ for $A$ defines also a braiding $c'$ on  $\lmod{\bb}{A}$ as:
\begin{equation*}
c'_{(M,r),(N,s)}=\psfrag{X}[Bc][Bc]{\scalebox{.8}{$M$}} \psfrag{Y}[Bc][Bc]{\scalebox{.8}{$N$}}
\psfrag{n}[Bc][Bc]{\scalebox{.9}{$r$}}\psfrag{s}[Bc][Bc]{\scalebox{.9}{$s$}} \psfrag{r}[Bc][Bc]{\scalebox{.9}{$\mathfrak{r}$}}
\rsdraw{.45}{1}{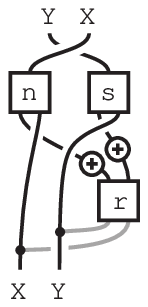}\;.
\end{equation*}
Furthermore, the map $\mathfrak{r} \mapsto c$ (resp.\@ $\mathfrak{r} \mapsto c'$) is a bijection between \Rt matrices for $A$ and braidings on $\bb_A$ (resp.\@ on $\lmod{\bb}{A}$).

A \emph{quasitriangular Hopf algebra} in $\bb$ is a Hopf algebra in $\bb$ endowed with an \Rt matrix.

\begin{rem}\label{rem-rmatmirror}
Let $A$ be a quasitriangular Hopf algebra in $\bb$. By construction, the monoidal isomorphism $F_A \co (\lmod{\bb}{A})^{\otimes\opp} \to \rmod{\bb}{A}$ of Remark~\ref{rem-braidleftright} is braided.
\end{rem}

\begin{rem}\label{rem-rmatmirrorbis}
Let $A$ be a quasitriangular Hopf algebra in $\bb$. Combining Remark~\ref{rem-rmatmirror} with Example~\ref{rem-isobraidedop}, we obtain that $\lmod{\bb}{A}$ and $\rmod{\bb}{A}$ are braided isomorphic.
\end{rem}

\subsection{The canonical distributive law of a Hopf algebra}\label{sect-canoniclawA}
Let $A$ be a Hopf algebra in a braided autonomous category $\bb$ which admits a coend $C$. By Section~\ref{sect-centralizersHA}, the centralizer of $? \otimes
A$ is $Z_{? \otimes A}=? \otimes Z(A)$, where $Z(A)=\ldual{A}
\otimes C$ is the centralizer of $A$.  It turns out that the
canonical distributive law of $? \otimes A$ over $Z_{?
\otimes A}$ is of the form
$\id_{1_\bb} \otimes \Omega$, where $\Omega \co Z(A) \otimes A \to A \otimes Z(A)$
is a comultiplicative distributive law of  $Z(A)$ over $A$ (see Example~\ref{exa-cartierHA}). We have:
\begin{equation*}
\psfrag{B}[Br][Br]{\scalebox{.8}{$\ldual{A}$}} \psfrag{C}[Bc][Bc]{\scalebox{.8}{$C$}}
\psfrag{A}[Bc][Bc]{\scalebox{.8}{$A$}} \Omega=\;\rsdraw{.4}{1}{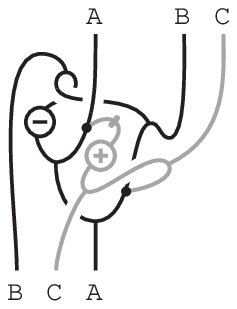} \quad \text{and} \quad \Omega^{-1}=\;\rsdraw{.4}{1}{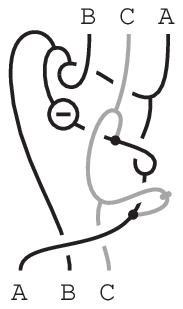} .
\end{equation*}
We call $\Omega$ the  \emph{canonical distributive law of $A$}.

\begin{rem}
By Theorem~\ref{thm-centraZ}, $Z_{\bb_A}(M,r)=\bigl(M \otimes Z(A), (r \otimes \id_{Z(A)})(\id_M \otimes\Omega)\bigr)$ is the centralizer of the category $\bb_A$.
In particular, the coend of~$\bb_A$ is $Z_{\bb_A}(\un,\varepsilon_A)=(Z(A),\alpha)$, where:
\begin{equation*}
\psfrag{B}[Br][Br]{\scalebox{.8}{$\ldual{A}$}} \psfrag{C}[Bc][Bc]{\scalebox{.8}{$C$}}
\psfrag{A}[Bc][Bc]{\scalebox{.8}{$A$}} \alpha=(\varepsilon_A \otimes \id_{Z(A)})\Omega=\rsdraw{.4}{1}{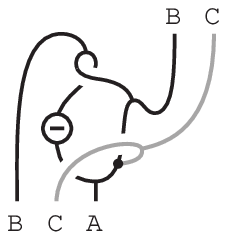} \co Z(A) \otimes A \to Z(A).
\end{equation*}
If $A$ is quasitriangular, so that $\bb_A$ is braided, then $(Z(A),\alpha)$ is a Hopf algebra in~$\bb_A$ which represents $Z_{\bb_A}$ on the right (see Theorem~\ref{thm-cent-br}). However, in this case, this Hopf algebra $(Z(A),\alpha)$ is not in general a lift to $\bb_A$ of the Hopf algebra $Z(A)$.
\end{rem}

\subsection{The Double of a Hopf algebra in a braided category}\label{sect-doubleA} Let $A$ be a Hopf algebra in a braided autonomous category $\bb$ which admits a coend $C$.

Let $Z(A)$ be the centralizer of $A$ (see Section~\ref{sect-centralizersHA}) and $\Omega$ be the canonical distributive law of $A$ (see Section~\ref{sect-canoniclawA}). By Example~\ref{exa-cartierHA},
\begin{equation*}
D(A)=A \otimes_\Omega Z(A)=A\otimes \ldual{A} \otimes C
\end{equation*}
is a Hopf algebra in $\bb$. Since $? \otimes Z(A)$ is the centralizer of $?\otimes A$ (see Section~\ref{sect-centralizersHA}), the Hopf monad $?\otimes D(A)$ is the double of $? \otimes A$, and so admits an \Rt matrix by Theorem~\ref{thm-doublable}, which turns out to be encoded by the following \Rt matrix for  $D(A)$:
\begin{equation*}
\psfrag{B}[Br][Br]{\scalebox{.8}{$\ldual{A}$}}
\psfrag{C}[Bc][Bc]{\scalebox{.8}{$C$}}
\psfrag{A}[Bc][Bc]{\scalebox{.8}{$A$}}
\mathfrak{r}=\,\rsdraw{.4}{1}{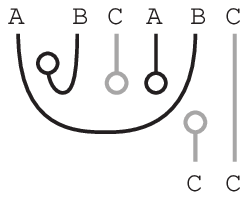} \;\co C \otimes C \to D(A)
\otimes D(A).
\end{equation*}
The quasitriangular Hopf algebra $D(A)$ is called the \emph{double of $A$}.

\begin{rem}
The canonical distributive law of $A$ is the unique comultiplicative distributive law
$\Omega$ of $Z(A)$ over  $A$ such
that the morphism $\mathfrak{r}$ above is an \Rt matrix for $A \otimes_\Omega Z(A)$, see
Remark~\ref{rem-dri-caract}.
\end{rem}

\begin{thm}
Let $A$ be a Hopf algebra in a braided autonomous category $\bb$ admitting a coend $C$ and
$D(A)=A \otimes_\Omega Z(A)$ be the double of $A$.
We have isomorphisms of braided categories:
\begin{equation*}
\zz(\rmod{\bb}{A}) \simeq \rmod{\bb}{{D(A)}}\simeq\lmod{\bb}{{D(A)}} \simeq \zz'(\lmod{\bb}{A}) \simeq \overline{\zz(\lmod{\bb}{A})}.
\end{equation*}
\end{thm}
\begin{proof}
By construction,  the quasitriangular Hopf monad $? \otimes D(A)$ is the double of $? \otimes A$. Hence the first braided isomorphism by Theorem~\ref{thm-doublable}. By Remark~\ref{rem-rmatmirrorbis}, $\lmod{\bb}{D(A)} \simeq
\rmod{\bb}{D(A)}$ as braided categories since $D(A)$ is quasitriangular. Finally, we  have the following isomorphisms of braided categories:
\begin{align*}
\zz(\rmod{\bb}{A})
& \simeq \zz(\rmod{\bb}{A})^{\otimes \opp} \quad \text{by Remark~\ref{rem-isobraidedop}}\\
  & \simeq \zz\bigl((\lmod{\bb}{A})^{\otimes \opp}\bigr)^{\otimes \opp} \quad \text{by Remark~\ref{rem-rmatmirror}}\\
& \simeq \zz'(\lmod{\bb}{A}) \simeq \overline{\zz(\lmod{\bb}{A})} \quad \text{by Remark~\ref{rem-alter-center}}.
\end{align*}
This completes the proof of the theorem.
\end{proof}

\begin{rem}
When $\bb=\vect_\kk$ is the category of finite-dimensional vector
spaces over a field $\Bbbk$, we recover the usual Drinfeld double
and the interpretation of its category of modules in terms of the
center. More precisely,  let $H$ be a finite-dimensional Hopf algebra
over $\Bbbk$ and $(e_i)$ be a basis of $H$ with dual basis $(e^i)$. Then $D(H)=H \otimes (H^*)^\mathrm{cop}$ is a quasitriangular Hopf algebra over $\Bbbk$, with \Rt matrix $\mathfrak{r}=\sum_i   e_i \otimes \varepsilon \otimes 1_H \otimes e_i$, such that:
\begin{equation*}
 \zz\bigl((\vect_\kk)_H\bigr)\simeq (\vect_\kk)_{D(H)} \simeq \lmod{(\vect_\kk)}{D(H)} \simeq \zz'\bigl(\lmod{(\vect_\kk)}{H}\bigr) \simeq \overline{\zz\bigl(\lmod{(\vect_\kk)}{H}\bigr)}
\end{equation*}
as braided categories.
\end{rem}

\begin{rem}
By Remark~\ref{rem-cartieriso}, $\Omega^{-1}$ is a distributive law of $Z(A)$ over $A$ and induces an isomorphism of Hopf algebras:
\begin{equation*}
D(A)=A \otimes_\Omega Z(A) \iso Z(A) \otimes_{\Omega^{-1}} A.
\end{equation*}
Via this isomorphism, the \Rt matrix $\mathfrak{r}$ of $D(A)$ is sent to the \Rt matrix:
\begin{equation*}
\psfrag{B}[Br][Br]{\scalebox{.8}{$\ldual{A}$}}
\psfrag{C}[Bc][Bc]{\scalebox{.8}{$C$}}
\psfrag{A}[Bc][Bc]{\scalebox{.8}{$A$}}
\mathfrak{r}'=\,\rsdraw{.4}{1}{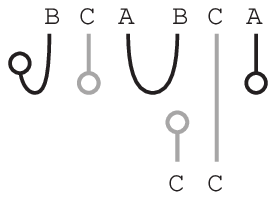} \; \co C \otimes C \to \bigl(Z(A) \otimes_{\Omega^{-1}} A \bigr) \otimes \bigl(Z(A) \otimes_{\Omega^{-1}} A\bigr)
\end{equation*}
of $Z(A) \otimes_{\Omega^{-1}} A$.
\end{rem}

\begin{rem}
Let $\bb$ a braided autonomous category which admits a coend $C$.  Then $D(\un)=C$ as a Hopf algebra. Hence $C$ is quasitriangular, with \Rt matrix $\mathfrak{r}=u_C\varepsilon_C
\otimes \id_C$, and
$\zz(\bb)\simeq\bb_C\simeq\lmod{\bb}{C}\simeq\overline{\zz(\bb)}$ as braided categories.  In other words, the
center of~$\bb$ is self-mirror and is the category of $C$\ti modules in $\bb$.
\end{rem}

\section{Hopf monads and fusion categories}\label{sect-fusion-cat}
In this section, given a \kt linear Hopf monad $T$ of a fusion category $\ff$, we describe explicitly the centralizer of $T$ and the canonical distributive law of $T$. Hence, in particular, a description of the coend of $\zz(\ff)$, which is used in \cite{BV4} to show that the center of a
spherical fusion category is modular and in \cite{BV5} to provide an algorithm for computing the Reshetikhin-Turaev invariant $\mathrm{RT}_{\zz(\ff)}$ in terms of~$\cc$ itself.

\subsection{Fusion categories}\label{sect-fusion}
A \emph{fusion category} over a commutative ring $\kk$ is a \kt
linear autonomous category $\ff$, whose monoidal product
$\otimes$ is  \kt linear in each variable, endowed with a finite
family $\{V_i\}_{i \in I}$ of objects of $\ff$ satisfying:
\begin{itemize}
\item $\Hom_\ff(V_i,V_j)=\delta_{i,j} \, \kk$ for all $i,j \in I$;
\item each object of $\ff$ is a finite direct sum of objects of $\{V_i\}_{i \in I}$;
\item $\un$ is isomorphic to $V_0$ for some $0 \in I$.
\end{itemize}

Let $\ff$ be a fusion category. The family $\{V_i\}_{i \in I}$ is a representative family of scalar objects of
$\ff$ (an object $X$ of \kt linear category is said to be \emph{scalar} if $\End(X)=\kk$).  The $\Hom$ spaces in
$\ff$ are free $\kk$-modules of finite rank. The \emph{multiplicity} of $i\in I$ in an object $X$ of $\ff$ is defined
as:
\begin{equation*}
N^i_X=\rank_\kk \,\Hom_\ff(V_i,X)=\rank_\kk \,\Hom_\ff(X,V_i).
\end{equation*}
For each object $X$ of $\ff$, we choose families of morphisms $(p_X^{i,\alpha} \co X \to V_i)_{1 \leq \alpha \leq N^i_X}$ and $(q_X^{i,\alpha}\co V_i \to
X)_{1 \leq \alpha \leq N^i_X}$ such that:
\begin{equation*}
\id_X=\!\!\!\!\!\sum_{\substack{i \in I \\ 1 \leq \alpha \leq N^i_X}} \!\!\!q_X^{i,\alpha}p_X^{i,\alpha} \quad
\text{and} \quad p_X^{i,\alpha}q_X^{j,\beta}=\delta_{i,j}\delta_{\alpha,\beta} \, \id_{V_i}.
\end{equation*}

\subsection{Centralizers in fusion categories}\label{sect-centerfusion} Let $\ff$ be a fusion
category over a commutative ring $\kk$ and $T$ be a \kt linear endofunctor $T$ of $\ff$. Then $T$ is centralizable,
with centralizer $(Z_T,\partial)$ given by:
\begin{equation*}
Z_T(X)=\bigoplus_{i \in I} \ldual{T}(V_i) \otimes X \otimes V_i
\end{equation*}
and
\begin{equation*}
\partial_{X,Y}=\sum_{\substack{i\in I \\ 1 \leq \alpha \leq N_Y^i}}
\bigl(T(q_Y^{i,\alpha})\otimes \id_{\ldual{T}(V_i) \otimes X} \otimes p_Y^{i,\alpha}\bigr)\bigl(\lcoev_{T(V_i)} \otimes
\id_{X\otimes Y}\bigr).
\end{equation*}
In particular, a fusion category is centralizable, with centralizer $Z=Z_{1_\ff}$ given by:
\begin{equation*}
Z(X)=\bigoplus_{i \in I} \ldual{V}_i \otimes X \otimes V_i.
\end{equation*}
\begin{rem} By Corollary~\ref{cormoncenter},
the centralizer $Z$ of $\ff$ provides in particular a left adjoint $F_Z$ to the forgetful functor $\uu\co \zz(\ff) \simeq Z\ti\ff\to \ff$, which is called the induction functor in \cite{ENO}.
\end{rem}
\begin{rem}\label{rem-Znotrepres}
In general, the centralizer  $Z$ of $\ff$ is not isomorphic (as an endofunctor of $\ff$) to $Z(\un) \otimes ?$ nor to $? \otimes Z(\un)$, as shown by the following counter-example. Let $G$ be a non-commutative finite group and let $\ff$ be the fusion category of finite-dimensional $G$\ti graded vector spaces over a field $\kk$. The elements of $G$ form a representative set of scalar objects of $\ff$. Then $Z(x)=\bigoplus_{g \in G} g^{-1} x g$ for $x \in G$. In particular $Z(\un)=\un^{\# G}$. Now, if $x\in G$ is not central, $Z(x)$ is not isomorphic to $Z(\un) \otimes x \simeq x^{\#  G} \simeq x \otimes Z(\un)$.
\end{rem}

\begin{figure}[t]
\begin{center}
         $\displaystyle (Z_T)_2(X,Y)=\hspace*{-.2cm}\sum_{\substack{i,k \in I \\ 1 \leq \alpha \leq N^k_{T(V_i)}}}$\hspace*{-.2cm}
\psfrag{o}[Bc][Bc]{\scalebox{.65}{$V_i$}}
\psfrag{a}[Bc][Bc]{\scalebox{.65}{$V_k$}}
\psfrag{n}[Bc][Bc]{\scalebox{.65}{$X$}}
\psfrag{u}[Bc][Bc]{\scalebox{.65}{$Y$}}
\psfrag{c}[Bc][Bc]{\scalebox{.65}{$\ldual{T}(V_i)$}}
\psfrag{s}[Bc][Bc]{\scalebox{.65}{$\ldual{T}(V_k)$}}
\psfrag{p}[c][c]{\scalebox{.7}{$\ldual{T}\bigl(q^{k,\alpha}_{T(V_i)}\bigr)$}}
\psfrag{q}[c][c]{\scalebox{.7}{$p^{k,\alpha}_{T(V_i)}$}}
\psfrag{m}[c][c]{\scalebox{.7}{$\ldual{\mu}_{V_i}$}}
\rsdraw{.5}{1}{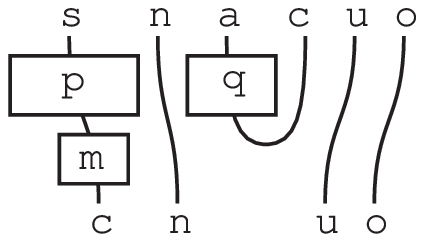} , \quad $\displaystyle (Z_T)_0=\sum_{i \in I}$\,
\psfrag{o}[Bc][Bc]{\scalebox{.65}{$V_i$}}
\psfrag{c}[Bc][Bc]{\scalebox{.65}{$\ldual{T}(V_i)$}}
\psfrag{q}[c][c]{\scalebox{.7}{$\eta_{V_i}$}}
\rsdraw{.5}{1}{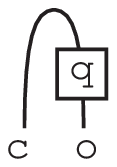} , \\[.5em]
$\displaystyle m_X=\hspace*{-.2cm}\sum_{\substack{i,j,k \in I \\ 1 \leq \alpha \leq N_{V_i\otimes V_j}^k }}$\,
\psfrag{u}[Bc][Bc]{\scalebox{.65}{$V_i$}}
\psfrag{r}[Bc][Bc]{\scalebox{.65}{$V_j$}}
\psfrag{a}[Bc][Bc]{\scalebox{.65}{$V_k$}}
\psfrag{o}[Bc][Bc]{\scalebox{.65}{$\ldual{T}(V_k)$}}
\psfrag{c}[Bc][Bc]{\scalebox{.65}{$\ldual{T}(V_j)$}}
\psfrag{s}[Bc][Bc]{\scalebox{.65}{$\ldual{T}(V_i)$}}
\psfrag{n}[Bc][Bc]{\scalebox{.65}{$X$}}
\psfrag{p}[c][c]{\scalebox{.7}{$\ldual{T}\bigl(q^{k,\alpha}_{V_i\otimes V_j}\bigr)$}}
\psfrag{q}[c][c]{\scalebox{.7}{$p^{k,\alpha}_{V_i\otimes V_j}$}}
\psfrag{m}[c][c]{\scalebox{.7}{$\ldual{T}_2(V_i,V_j)$}}
\rsdraw{.5}{1}{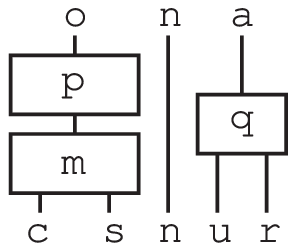} , \qquad $\displaystyle u_X=$
\psfrag{n}[Bc][Bc]{\scalebox{.65}{$V_0$}}
\psfrag{c}[Bc][Bc]{\scalebox{.65}{$X$}}
\psfrag{o}[Bc][Bc]{\scalebox{.65}{$\ldual{T}(V_0)$}}
\psfrag{R}[c][c]{\scalebox{.7}{$T_0$}}
\rsdraw{.5}{1}{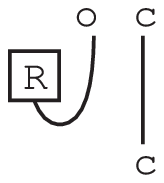} , \\[.5em]
$\displaystyle S^l_X=\hspace*{-.3cm}\sum_{\substack{i,j \in I \\ 1 \leq \alpha \leq N^i_{T(V_j)^\vee}}}$\hspace*{-.2cm}
\psfrag{u}[Bc][Bc]{\scalebox{.65}{$\lldual{T}(V_i)$}}
\psfrag{r}[Bc][Bc]{\scalebox{.65}{$V_j$}}
\psfrag{c}[Bc][Bc]{\scalebox{.65}{$\ldual{T}(V_j)$}}
\psfrag{s}[Bc][Bc]{\scalebox{.65}{$\ldual{V}_i$}}
\psfrag{n}[Bc][Bc]{\scalebox{.65}{$\ldual{X}$}}
\psfrag{p}[c][c]{\scalebox{.7}{$\lldual{p}^{i,\alpha}_{T(V_j\rdual{)}}$}}
\psfrag{m}[c][c]{\scalebox{.7}{$\lldual{T}\bigl(q^{i,\alpha}_{T(V_j\rdual{)}} \bigr)$}}
\psfrag{o}[c][c]{\scalebox{.7}{$\ldual{(}s^r_{V_j})$}}
\rsdraw{.35}{1}{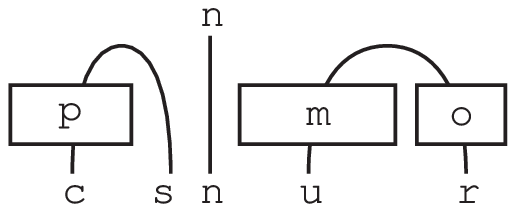} , \\[.5em]
$\displaystyle S^r_X=\hspace*{-.3cm}\sum_{\substack{i,j \in I \\ 1 \leq \alpha \leq
N^i_{\ldual{T}(V_j)}}}$\hspace*{-.2cm}
\psfrag{u}[Bc][Bc]{\scalebox{.65}{$T(V_i)$}}
\psfrag{r}[Bc][Bc]{\scalebox{.65}{$V_j$}}
\psfrag{c}[Bc][Bc]{\scalebox{.65}{$\ldual{T}(V_j)$}}
\psfrag{s}[Bc][Bc]{\scalebox{.65}{$\rdual{V}_i$}}
\psfrag{n}[Bc][Bc]{\scalebox{.65}{$\rdual{X}$}}
\psfrag{p}[c][c]{\scalebox{.7}{$p^{i,\alpha}_{\ldual{T}(V_j)}$}}
\psfrag{m}[c][c]{\scalebox{.7}{$T\bigl(q^{i,\alpha}_{\ldual{T}(V_j)} \bigr)$}}
\psfrag{o}[c][c]{\scalebox{.7}{$(s^l_{V_j}\rdual{)}$}}
\rsdraw{.35}{1}{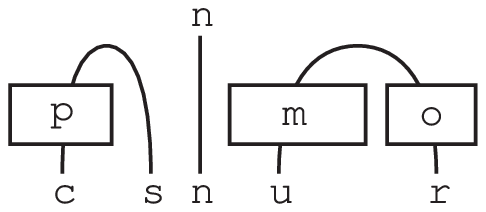}  .
\end{center}
\caption{Structural morphisms of $Z_T$}
\label{morphZ}
\end{figure}

Assume $T$ is a centralizable Hopf monad. By Theorem~\ref{doub1}, its centralizer $Z_T$ is a Hopf monad on $\ff$ and its structural morphisms can be
described purely in terms of those of $T$ and of the category~$\ff$ (that is, the $p,q$'s  and the
duality morphisms). They are depicted in Figure~\ref{morphZ}, where $\mu$, $\eta$, $s^l$, $s^r$ (resp.\@ $m$, $u$,
$S^l$, $S^r$) denote the product, unit, left antipode, and right antipode of $T$ (resp.\@ $Z_T$).
The canonical distributive law of $T$ is:
\begin{equation*}
\Omega_X=\hspace*{-1em}\sum_{\substack{i,j \in I \\ 1 \leq \alpha \leq N^j_{T(V_i)}}}  \hspace*{-1em} \Bigl(
\ldual{T}\bigl(q^{j,\alpha}_{T(V_i)}\bigr) \ldual{\mu}_{V_i} s^l_{T(V_i)} T\bigl(\ldual{\mu}_{V_i}\bigr)\otimes \id_{T(X)}
\otimes p^{j,\alpha}_{T(V_i)} \Bigr) T_3\bigl(\ldual{T}(V_i),X,V_i\bigr).
\end{equation*}
Hence an explicit description of the double $D_T=Z_T \circ_\Omega T$ of $T$ and of the lift $\Tilde{Z}_T^\Omega$ of~$Z_T$ to
$T\ti\ff$. Note that the \Rt matrix of $D_T$ is:
\begin{equation*}
R_{X,Y}=\sum_{\substack{i\in I \\ 1 \leq \alpha \leq N_Y^i}}
\bigl(\leftidx{^\vee}{T}{_0}\otimes T(q_Y^{i,\alpha})\otimes \id_{\ldual{T}(V_i)} \otimes \eta_X \otimes p_Y^{i,\alpha}\bigr)\bigl(\lcoev_{T(V_i)} \otimes
\id_{X\otimes Y}\bigr).
\end{equation*}

\subsection{The coend of the center of a fusion category}\label{sect-coendcenterfusion} 
\begin{figure}[t]
   \begin{center}
       $\displaystyle \Delta_C=\hspace*{-.4cm}\sum_{\substack{i,j,k,m,n \in I \\ 1 \leq \alpha \leq N_{k,m}^k \\ 1 \leq
\beta \leq N_{\ldual{k},j,k}^n}}$
 \psfrag{u}[Bc][Bc]{\scalebox{.55}{$\ldual{V}_i$}}
 \psfrag{r}[Bc][Bc]{\scalebox{.55}{$\ldual{V}_j$}}
 \psfrag{i}[Bc][Bc]{\scalebox{.55}{$\lldual{V}_i$}}
 \psfrag{b}[Bc][Bc]{\scalebox{.55}{$V_j$}}
 \psfrag{o}[Bc][Bc]{\scalebox{.55}{$\ldual{V}_m$}}
 \psfrag{c}[Bc][Bc]{\scalebox{.55}{$\ldual{V}_n$}}
 \psfrag{s}[Bc][Bc]{\scalebox{.55}{$\lldual{V}_m$}}
 \psfrag{n}[Bc][Bc]{\scalebox{.55}{$V_n$}}
 \psfrag{w}[Bc][Bc]{\scalebox{.55}{$\ldual{V}_k$}}
 \psfrag{e}[Bc][Bc]{\scalebox{.55}{$\ldual{V}_j$}}
 \psfrag{z}[Bc][Bc]{\scalebox{.55}{$\lldual{V}_k$}}
 \psfrag{p}[c][c]{\scalebox{.7}{$\ldual{p}^{n,\beta}_{\ldual{k},j,k}$}}
 \psfrag{q}[c][c]{\scalebox{.7}{$\ldual{q}^{n,\beta}_{\ldual{k},j,k}$}}
 \psfrag{a}[c][c]{\scalebox{.7}{$\ldual{p}^{i,\alpha}_{k,m}$}}
 \psfrag{x}[c][c]{\scalebox{.7}{$\lldual{q}^{i,\alpha}_{k,m}$}}
\rsdraw{.6}{.8}{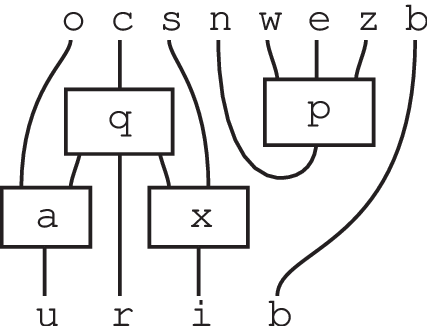} , \qquad $\displaystyle\varepsilon_C= \sum_{j \in I}$ \,
 \psfrag{o}[Bc][Bc]{\scalebox{.55}{$\ldual{V}_0$}}
 \psfrag{c}[Bc][Bc]{\scalebox{.55}{$\ldual{V}_j$}}
 \psfrag{n}[Bc][Bc]{\scalebox{.55}{$\lldual{V}_0$}}
 \psfrag{e}[Bc][Bc]{\scalebox{.55}{$V_j$}}
\rsdraw{.4}{.8}{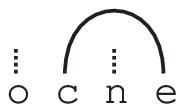} ,\\[.5em]
       $\displaystyle m_C=\hspace*{-.6cm}\sum_{\substack{i,j,k,l,m,n \in I \\ 1 \leq \alpha \leq N_{\ldual{k},l,k}^n \\ 1 \leq \beta \leq N_{\ldual{n},k,n}^i \\  1 \leq \gamma \leq N_{\lldual{n},j,\ldual{n},l}^m}}$
 \psfrag{u}[Bc][Bc]{\scalebox{.55}{$\ldual{V}_k$}}
 \psfrag{r}[Bc][Bc]{\scalebox{.55}{$\ldual{V}_l$}}
 \psfrag{i}[Bc][Bc]{\scalebox{.55}{$\lldual{V}_k$}}
 \psfrag{b}[Bc][Bc]{\scalebox{.55}{$V_l$}}
 \psfrag{o}[Bc][Bc]{\scalebox{.55}{$\ldual{V}_i$}}
 \psfrag{c}[Bc][Bc]{\scalebox{.55}{$\ldual{V}_j$}}
 \psfrag{s}[Bc][Bc]{\scalebox{.55}{$\lldual{V}_i$}}
 \psfrag{n}[Bc][Bc]{\scalebox{.55}{$V_j$}}
  \psfrag{k}[Bc][Bc]{\scalebox{.55}{$\ldual{V}_k$}}
 \psfrag{b}[Bc][Bc]{\scalebox{.55}{$\ldual{V}_m$}}
 \psfrag{l}[Bc][Bc]{\scalebox{.55}{$\lldual{V}_k$}}
 \psfrag{h}[Bc][Bc]{\scalebox{.55}{$V_m$}}
 \psfrag{e}[Bc][Bc]{\scalebox{.7}{$\ldual{q}^{m,\gamma}_{\lldual{n},j,\ldual{n},l}$}}
 \psfrag{v}[Bc][Bc]{\scalebox{.7}{$p^{m,\gamma}_{\lldual{n},j,\ldual{n},l}$}}
 \psfrag{q}[c][c]{\scalebox{.7}{$\ldual{q}^{n,\alpha}_{\ldual{k},l,k}$}}
 \psfrag{p}[c][c]{\scalebox{.7}{$\ldual{p}^{n,\alpha}_{\ldual{k},l,k}$}}
 \psfrag{a}[c][c]{\scalebox{.7}{$\ldual{p}^{i,\beta}_{\ldual{n},k,n}$}}
 \psfrag{x}[c][c]{\scalebox{.7}{$\lldual{q}^{i,\beta}_{\ldual{n},k,n}$}}
\rsdraw{.6}{.8}{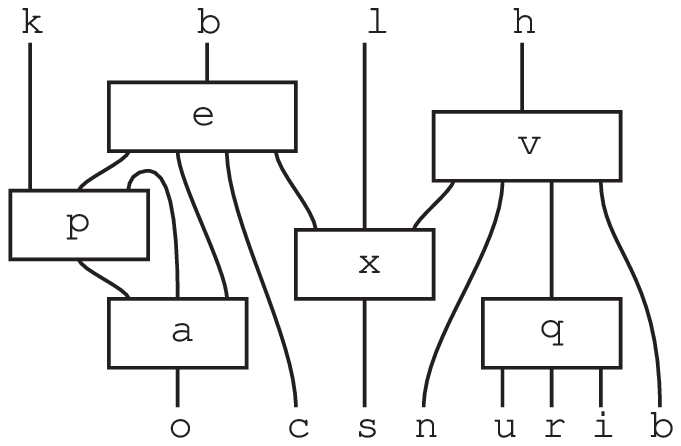}, \qquad $\displaystyle u_C= \sum_{i \in I}$ \,
 \psfrag{o}[Bc][Bc]{\scalebox{.55}{$\ldual{V}_i$}}
 \psfrag{c}[Bc][Bc]{\scalebox{.55}{$\ldual{V}_0$}}
 \psfrag{n}[Bc][Bc]{\scalebox{.55}{$V_0$}}
 \psfrag{e}[Bc][Bc]{\scalebox{.55}{$\lldual{V}_i$}}
\rsdraw{.4}{.8}{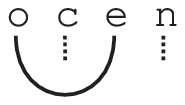} ,\\[.5em]
$\displaystyle S_C=\hspace*{-.8cm}\sum_{\substack{i,j,k,l \in I \\ 1 \leq \alpha \leq N^{\ldual{i}}_{j,k,\rdual{j}} \\ 1
\leq \beta \leq N^l_{\ldual{j},\ldual{i}, \ldual{j}, \lldual{i}, j}}}$
 \psfrag{x}[c][c]{\scalebox{.7}{$p^{l,\beta}_{\ldual{j},\ldual{i}, \ldual{j}, \lldual{i}, j}$}}
 \psfrag{a}[c][c]{\scalebox{.7}{$\ldual{q}^{l,\beta}_{\ldual{j},\ldual{i}, \ldual{j}, \lldual{i}, j}$}}
 \psfrag{p}[c][c]{\scalebox{.7}{$\ldual{p}^{\ldual{i},\alpha}_{j,k,\rdual{j}}$}}
 \psfrag{q}[c][c]{\scalebox{.7}{$\lldual{q}^{\ldual{i},\alpha}_{j,k,\rdual{j}}$}}
 \psfrag{o}[Bc][Bc]{\scalebox{.55}{$\ldual{V}_k$}}
 \psfrag{c}[Bc][Bc]{\scalebox{.55}{$\ldual{V}_l$}}
 \psfrag{s}[Bc][Bc]{\scalebox{.55}{$\lldual{V}_k$}}
 \psfrag{n}[Bc][Bc]{\scalebox{.55}{$V_l$}}
\rsdraw{.5}{.8}{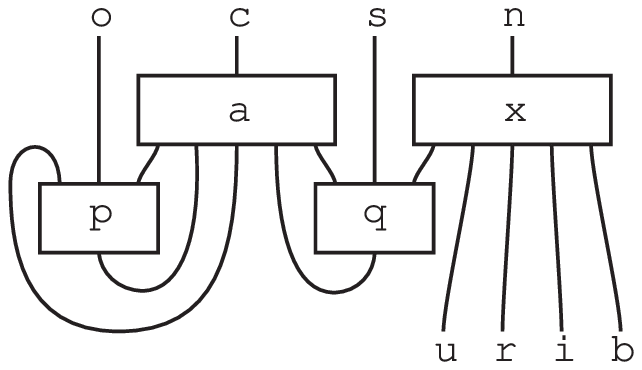} ,\\[.5em]
$\displaystyle \omega_C=\hspace*{-.6cm}\sum_{\substack{i,j,k,l \in I \\ 1 \leq \alpha \leq N^{\ldual{k}}_{\ldual{i},j,i} \\
1 \leq \beta \leq N^i_{\ldual{k},\ldual{l},\lldual{k}}}}$\;
 \psfrag{x}[c][c]{\scalebox{.7}{$p^{\ldual{k},\alpha}_{\ldual{i},j,i}$}}
 \psfrag{a}[c][c]{\scalebox{.7}{$q^{i,\beta}_{\ldual{k},\ldual{l},\lldual{k}}$}}
 \psfrag{p}[c][c]{\scalebox{.7}{$\ldual{q}^{\ldual{k},\alpha}_{\ldual{i},j,i}$}}
 \psfrag{q}[c][c]{\scalebox{.7}{$p^{i,\beta}_{\ldual{k},\ldual{l},\lldual{k}}$}}
 \psfrag{o}[Bc][Bc]{\scalebox{.55}{$\ldual{V}_i$}}
 \psfrag{c}[Bc][Bc]{\scalebox{.55}{$\ldual{V}_j$}}
 \psfrag{s}[Bc][Bc]{\scalebox{.55}{$\lldual{V}_i$}}
 \psfrag{n}[Bc][Bc]{\scalebox{.55}{$V_j$}}
 \psfrag{u}[Bc][Bc]{\scalebox{.55}{$\ldual{V}_k$}}
 \psfrag{r}[Bc][Bc]{\scalebox{.55}{$\ldual{V}_l$}}
 \psfrag{i}[Bc][Bc]{\scalebox{.55}{$\lldual{V}_k$}}
 \psfrag{b}[Bc][Bc]{\scalebox{.55}{$V_l$}}
\rsdraw{.55}{.8}{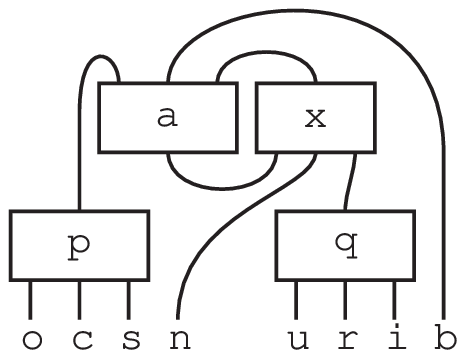} .
   \end{center}
     \caption{Structural morphisms of the coend of $\zz(\cc)$}
     \label{morphcoendcenter}
\end{figure}

Let $\ff$ be a fusion
category over a commutative ring $\kk$, and denote $Z$ the centralizer of $\ff$. Recall that $Z$ is a quasitriangular Hopf monad on $\ff$ such that $\zz(\ff) \simeq Z\ti\ff$, see Section~\ref{sect-centerfusion}.
Since $Z$ is \kt linear, it is centralizable. Denote $Z_Z$ its centralizer and $\Omega$
the canonical distributive law of $Z$ over $Z_Z$.
Then the coend of $Z(\ff)$ is:
\begin{equation*}
C=\Tilde{Z}_Z^\Omega(\un,Z_0)=\bigl(Z_Z(\un),Z_Z(Z_0)\Omega_\un\bigr).
\end{equation*}
Note that:
\begin{equation*}
Z_Z(\un)=\bigoplus_{j \in I} \ldual{Z(V_j)} \otimes V_j= \bigoplus_{i,j \in I} \ldual{V}_i \otimes \ldual{V}_j \otimes
\lldual{V}_i \otimes V_j.
\end{equation*}
Using the results of Section~\ref{sect-centralizerTC}, one computes the Hopf algebra structure of $C$ and its self-dual Hopf pairing. These are depicted in Figure~\ref{morphcoendcenter}, where we denote $A_{V_{i_1} \otimes \cdots \otimes V_{i_n}}$ by $A_{i_1, \dots, i_n}$ for $A=p^{i,\alpha}$,
$q^{i,\alpha}$, or~$N^i$, and the dotted lines represent the relevant isomorphism between $\un$ and $V_0$ or its duals.

In \cite{BV4}, we use this explicit description of the coend of $\zz(\ff)$ to show that the center $\zz(\ff)$ of a
spherical fusion category $\ff$ is modular. In particular, this implies that if $\ff$ is a spherical fusion category
of  invertible dimension over an algebraic closed field $\kk$, then $\zz(\ff)$ is a modular ribbon fusion category
(this last result was first shown in~\cite{Mueg} using different methods).

Also, this description of the coend of $\zz(\ff)$ leads to an explicit algorithm (involving Hopf diagrams \cite{BV1})
for computing the Reshetikhin-Turaev invariants defined with $\zz(\ff)$. Moreover, this approach allows one to  define these invariants over an arbitrary base ring, without assumption on the dimension of $\ff$ (if the dimension of $\ff$ is not invertible, this yields `non-semisimple' invariants). See~\cite{BV5} for details.

\end{document}